\documentclass[a4paper,english,oneside]{smfcustom}

\newcounter{Version}
\ifnum \arabic{Version} < 2
\usepackage[notref,notcite]{showkeys}
\fi

\ifnum \arabic{Version} = 2
\fi

\usepackage{mathrsfs,mathtools}
\usepackage[english]{babel}
\usepackage[utf8]{inputenc}
\usepackage[T1]{fontenc}
\usepackage{pxfonts,microtype,enumitem}
\usepackage[dvipsnames]{xcolor}
\usepackage[all]{xy}
\usepackage{bm,framed}
\usepackage{etex}
\usepackage[labelfont=normal]{subfig}
\usepackage{tikz,forest}
\usetikzlibrary{patterns,decorations.markings,cd,arrows}
\usetikzlibrary{matrix,arrows,decorations.markings}
\usetikzlibrary{patterns,arrows,trees,shapes}
\usepackage{young}
\newdimen\dikte     \dikte=0.5pt        
\newdimen\hoogte    \hoogte=16pt        
\newdimen\breedte   \breedte=1.2\hoogte 
\newcommand{\eps}{\varepsilon}
\newcommand{\dd}{\mathrm d}
\newcommand{\id}{\mathrm {id} }
\renewcommand{\labelitemi}{$\bullet$} 
\usepackage{todonotes}

\usepackage{navigator}
\embeddedfile{Printed version of the companion CoCalc worksheet}{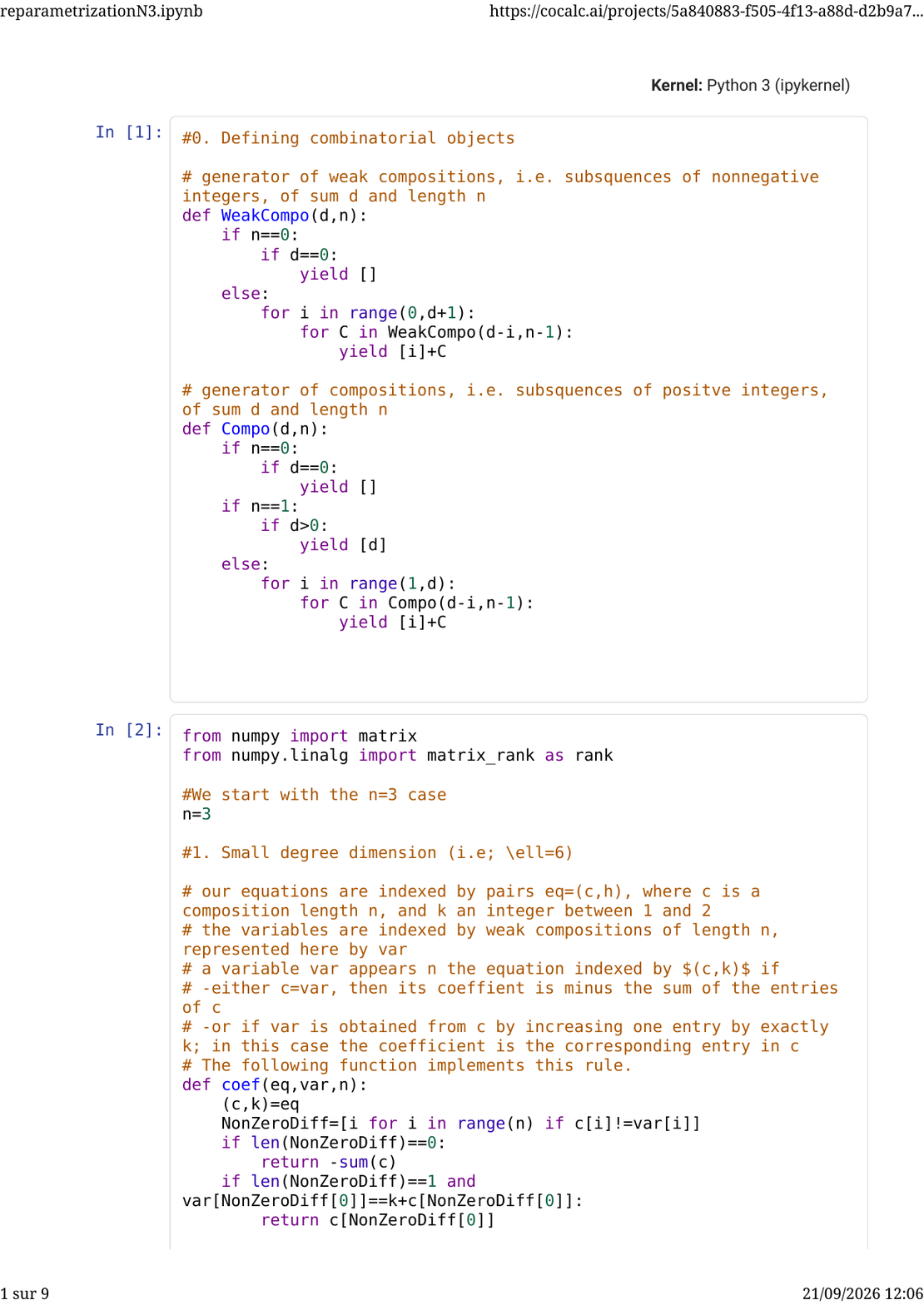}
\newcommand{\valentin}[1]{
\ifnum \arabic{Version} = 0
\todo[color=cyan!70,inline]{Valentin: #1}
\fi
}
\newcommand{\vf}[1]{{#1}}
\newcommand{\vupe}[1]{{#1}}
\newcommand{\pe}[1]{{#1}}
\newcommand{\pebox}[1]{
\ifnum \arabic{Version} = 0
\todo[color=Melon,inline]{PE: #1}
\fi
}

\newtheorem{thmx}{Theorem}
\newtheorem{theorem}{Theorem}[section]
\newtheorem{proposition}[theorem]{Proposition}
\newtheorem{conjecture}[theorem]{Conjecture}

\newtheorem{lemma}[theorem]{Lemma}
\newtheorem{corollary}[theorem]{Corollary}
\newtheorem{definition}[theorem]{Definition}
\newtheorem{notation}[theorem]{Notation}

\newtheorem{remark}[theorem]{Remark}

\newcommand{\T}{\S\kern .15em\relax }

\newcommand{\ev}{\mathrm{ev}}

\newcommand{\cO}{\mathcal{O}}

\newcommand{\Qnp}{{\Q^n_{>0}}}
\newcommand{\Qr}{{\widetilde \Q}}
\newcommand{\Qrp}{{\widetilde \Q_{>0}}} \newcommand{\Zrp}{{\Z_{1,>0}}}
\newcommand{\Rnp}{{\R^n_{>0}}}
\newcommand{\Rnmp}{{\R^{n-1}_{>0}}}
\newcommand{\beps}{{\bm \eps}}
\newcommand{\set}[1]{\{1,\ldots,#1\}}
\newcommand{\Ze}{\Z^n_{\beps,>0}}
\newcommand{\Qe}{\widetilde \Q^n_{\beps,>0}}
\newcommand{\Fe}{\cF_\beps}
\newcommand{\Se}{\cS_\beps}
\newcommand{\ST}[1]{\cS_{#1}}
\newcommand{\wvar}{{\widetilde \varphi}}
\newcommand{\scal}[1]{\langle #1 \rangle}
\newcommand{\iset}[1]{\{1,\ldots,#1\}}

\newcommand{\R}{\mathbb{R}}
\newcommand{\Q}{\mathbb{Q}}
\newcommand{\N}{\mathbb{N}}
\newcommand{\Z}{\mathbb{Z}}
\newcommand{\Zp}{\mathbb{Z}_{>0}}
\newcommand{\Znp}{\mathbb{Z}^n_{>0}}

\newcommand{\cC}{{\mathcal C}} \newcommand{\cD}{{\mathcal D}}
\newcommand{\cF}{{\mathcal F}} 
 \newcommand{\cI}{{\mathcal I}}
\newcommand{\cR}{{\mathcal R}} \newcommand{\cS}{{\mathcal S}}
 
\newcommand{\CR}{{\mathcal C}^\infty(\R)_0}
\newcommand{\CC}{{\mathbb C}}  

\newcommand{\Supp}{\mathrm{Supp}}

\newcommand{\fS}{{\mathfrak S}}
\newcommand{\im}{\mathrm{Im}}

\DeclareMathOperator{\Id}{Id}
\DeclareMathOperator{\AI}{AI}
\DeclareMathOperator{\LI}{LI}
\DeclareMathOperator{\I}{I}
\DeclareMathOperator{\Sym}{Sym}

\begin{document}

\title{Functionals commuting with reparametrizations}
\author{Pierre-Emmanuel Chaput, David Dos Santos Ferreira and Valentin Féray}
\address{Université de Lorraine, CNRS, IECL, F-54000 Nancy, France}
\email{pierre-emmanuel.chaput@univ-lorraine.fr, david.dos-santos-ferreira@univ-lorraine.fr, valentin.feray@univ-lorraine.fr}
\maketitle

\tableofcontents

\begin{abstract}
Motivated by the study of tangent curve bundles of algebraic varieties, 
we consider the following simple-looking problem: given a space $E$ of functions,
describe functionals from $E^n$ to $E$ which commutes with {\em reparametrization},
i.e.~with the right-composition with an element of $E$.
Our main result is a complete description of the algebraic solutions of this problem
when $E$ is the space of polynomials vanishing in $0$, \vupe{which we denote by $\mathbb C[T]_0$}.
This description involves a simple operation, which we call {\em semi-product},
with interesting algebraic and geometric properties.
The proof involves the study of polyharmonic functions on the positive discrete orthant
$\mathbb Z^n_{>0}$, and a relation with the co-invariant algebra of type $A$.
We also provide a generalization connected to colored permutations, and discuss the
analogue problems
where $E=\mathbb C[T]_0/(T^\ell)$ for some $\ell >0$, and where $E$ consists of infinitely differentiable functions.
Finding a complete description of algebraic solutions in these settings is still open.
\end{abstract}

\section{Introduction}

\subsection{Problem statement and main result}
Let us state our problem right away and discuss later the motivation.
Let $(Z,0)$ be a pointed set and $E$ be a set of
functions \vupe{from $Z$ to $Z$ vanishing at $0$.
Assuming that $E$ is} stable under composition, we aim at
describing multivariate functionals $\varphi : E^n \to E$ ,
which ``commute with reparametrizations'' in the following sense: 
for all $(f_1,\dots,f_n)$ in $E^n$, for all $\gamma$ in $E$,
\begin{equation}\label{e:reparametrage_intro}
 \varphi(f_1\circ \gamma,\dots,f_n \circ \gamma)=\varphi(f_1,\dots,f_n)\circ \gamma.
\end{equation}

%
%

Trivial solutions to this problem are given choosing a function
$H:Z^n \to Z$ mapping $(0,\ldots,0)$ to $0$, and then setting
\begin{equation}\label{eq:space_transfo}
\varphi(f_1,\ldots,f_n)(z)  \coloneqq  H(f_1(z),\ldots,f_n(z)).
\end{equation}
If differentiation and integration make sense on $E$,
a less trivial solution for $n=2$ is given by the following indefinite integral vanishing at $0$:
\begin{equation}
\label{e:intro_semi_product}
[f_1:f_2] \coloneqq \int f_1'f_2.
\end{equation}
That $\varphi(f_1,f_2)$ satisfies \eqref{e:reparametrage_intro} follows
from a simple change of variable formula in the integral.
Since all functions are assumed to vanish at $0$
(\vupe{the importance of this assumption is discussed in Remark~\ref{r:constant} below}),
we have$\int f_1'f_2 + \int f_1f_2' = f_1f_2$.
\vupe{This prompts us to call} $[f_1:f_2]$ a \emph{semi product}.
By induction, we also define, for $m \ge 1$, the \emph{iterated semi product}
\begin{equation}
\label{e:intro_iterated_semi_product}
	[f_1:f_2:\cdots:f_m] = [f_1:[f_2:\cdots:f_m]] = \int f_1' \left( \int f'_2 \left( \cdots \int f_{m-1}' f_m  \right)\right).
\end{equation}
This is also a solution of \eqref{e:reparametrage_intro}.

Our main result concerns the case in which $E$ is the space of polynomials in a single variable $T$ vanishing at $0$, written as $E =\CC[T]_0$.
In this case, we prove that iterated semi products, possibly applied to repeated arguments,
are the only {\em algebraic} functionals that commute with reparametrizations.
Here, {\em algebraic}  means that the coefficient of a given monomial $T^d$ in $\varphi(f_1,\dots,f_n)$
\vupe{depends polynomially on} the coefficients of $f_1$, \dots, $f_n$ (in particular it depends on finitely
many of them).
\begin{thmx}\label{thm:main_intro}
Let $E=\CC[T]_0$. Then the functions
\begin{equation}\label{e:def_varsigma_intro}
\varphi_\sigma(f_1,\dots,f_n) \coloneqq [f_{\sigma(1)} : \cdots : f_{\sigma(m)}],
\end{equation}
indexed by integers $m \ge 1$ and functions $\sigma:[1,\ldots,m] \to [1,\ldots,n]$,
form a linear basis of the space $\AI_n(E)$ of algebraic functionals $\varphi:(\CC[T]_0)^n \to \CC[T]_0$
satisfying \eqref{e:reparametrage_intro}.
\end{thmx}

\begin{remark}
 The semi-product operation is not associative and,
 obviously, $(f_1,f_2,f_3) \mapsto [ [f_1:f_2]:f_3]$ also commutes with reparametrizations.
 However, we have
 \begin{equation}\label{eq:Zinbiel}
  [ [f_1:f_2]:f_3]=\int \left(\int f_1' f_2\right)' f_3=\int f_1' f_2f_3=[f_1:f_2f_3]= [f_1:[f_2:f_3]] + [f_1:[f_3:f_2]]. 
  \end{equation}
Consequently, $(f_1,f_2,f_3) \mapsto [ [f_1:f_2]:f_3]$ is a linear combination of the functions $\varphi_\sigma$, as predicted by Theorem~\ref{thm:main_intro}.
 We note that the relation \eqref{eq:Zinbiel} satisfied by the semi-product is the defining
 relations of the so-called {\em dual Liebniz algebras} or {\em Zinbiel algebras}, introduced by Loday in \cite{LodayZinbiel}. Hence, $E$, equipped with the semi-product,
 is a Zinbiel algebra.
\end{remark}
\begin{remark}
\vupe{The spaces of (algebraic or general) solutions of \eqref{e:reparametrage_intro} also enjoy a nice algebraic structure. Indeed if $\varphi:E^n \to E$ is an (algebraic or general) solution and $\varphi_1$,\dots $\varphi_n$ are also (algebraic or general) solutions of \eqref{e:reparametrage_intro}, repectively on $E^{m_1}$,\dots, $E^{m_n}$, then
\[ (f_{1,1},\dots,f_{1,m_1},\dots,f_{n,1},\dots,f_{n,m_n}) \mapsto \varphi\big( \varphi_1(f_{1,1},\dots,f_{1,m_1}),\dots, \varphi_n(f_{n,1},\dots,f_{n,m_n})\big)\]
is also an (algebraic or general) solution of \eqref{e:reparametrage_intro}.
In other terms}, the spaces of solutions of \eqref{e:reparametrage_intro}
for all $n \ge 1$ form a suboperad of the {\em endomorphism operad} of $E$, i.e.~the operad of functions from $E^n$ to $E$, for all $n \ge 1$; see, e.g.,~\cite{Operads} for an introduction to operads. \vupe{This was used implicitely above to construct further solutions
from the semi-product operation. Our Theorem~\ref{thm:main_intro} means that the
operad of algebraic solutions  of \eqref{e:reparametrage_intro} is generated the semi-product operation.}
\end{remark}
\begin{remark}
Under suitable assumptions on $(f_1,f_2)$, the semi-product $[f_1,f_2]$ also has a geometric interpretation as area of some region of the plane,
which makes both the decomposition $f_1f_2=[f_1:f_2]+[f_2:f_1]$ 
and the commutation with reparametrization straight-forward.
See Section~\ref{sec:geometric_int_semi_prod}.
\end{remark}

\begin{remark}
\label{r:constant}
\vupe{We explain here the assumption that elements of $E$ vanish at $0$.
In the opposite situation where} $E$ contains all the constant functions, then  $\varphi$
must be of the form
$$\varphi(f_1,\ldots,f_n)(x)=H(f_1(x),\ldots,f_n(x))$$ for any $x$ in $Z$, \pe{where for $u \in Z$, $c_u$ is the constant function $x \mapsto u$
and $H$ is the function from $Z^n$ to $Z$ defined
by $H(x_1,\ldots,x_n)=\varphi(c_{x_1},\ldots,c_{x_n})(0)$.} Indeed, this follows by applying \eqref{e:reparametrage} with the constant function $\gamma = c_x$.
This is not very interesting, and we will only consider cases
where all functions vanish at $0$.
\end{remark}
\subsection{Motivation}
\label{sec:intro_motivation}
The initial motivation for this project is connected to the
computation of some space of sections of vector bundles in
the Demailly-Semple tower~\cite{CDM}. 
Assume we have a smooth algebraic variety $X$ over $\CC$,
with $V \subseteq T_X$ a subbundle of the tangent bundle;
concretely this means that for every $x$ in $X$,
we are given a subspace $V_x$ of the tangeant space $T_xX$,
with the map $x \mapsto V_x$ being algebraic.
Let $E$ be the quotient ring
$\CC[T]_0/T^\ell$ for some $\ell >0$.
We are interested in the set $\mathcal C_V$ of order $\ell$ jets of
curves from $(\CC,0)$ to $X$ that are tangent to
$V$, i.e.~such that, for each $t$,
the differential $d\gamma(t)$ lies in $V_{\gamma(t)}$.
\medskip

A simple case, usually refered to as the {\em integrable} case,
is when there exist a smooth algebraic variety $Y$ and
a smooth map $\pi:X \to Y$ such that, for each $x$ in $X$, we have $V_x=\ker(d\pi_x)$.
Then, assuming that a local inversion theorem holds (for example
if we are in the $\cC^\infty$ case), $X$ locally identifies
to a product $U \times \CC^r$, where $r$ is the rank of
$V$, and each $V_x$ to the factor $\{0\} \times \CC^r$ of
$T_x$. Moreover, it can be seen that
the set of tangeant jets of curves mapping $0$ to $x$
is isomorphic to $\mathbb C^r \otimes E$.
We say that $\mathcal C_V$ is locally trivial.
Demailly noticed \cite{demailly} that 
the set $\mathcal C_V$ is in fact always locally trivial,
even if $V$ is not integral.
\medskip

In \cite{CDM}, the critical observation is that some
important constructions on $V$ (such as
the Demailly-Semple tower) in fact only need the geometric
data of $\cC_V$. The set $\cC_V$
is organised as a bundle over
$X$, by mapping tangent curves $\gamma \in \mathcal C_V$
to $\gamma(0)$.
Even more importantly, it gives rise to a
principal bundle with structure group the set of
automorphisms of $\CC^r \otimes E$ (here $r$ is the rank
of $V$) which are equivariant for the action of $E$,
i.e.~with commute with reparametrizations. This leads to
an analysis of constructions such as vector bundles on
the Demailly-Semple tower with two separate ingredients~:
\begin{itemize}
 \item the structure group
 and its linear representations;
 \item the principal bundle
 associated to $\cC_V$, selecting the
 relevant geometric data among all the geometric information
 provided by the subbundle $V \subset T_X$.
\end{itemize}

An element of the structure group is given by $r$
functionals
$E^r \to E$ commuting with reparametrizations.
This is the reason why it is important,
in \cite{CDM}, to have information about those functionals.
In particular, an explicit description of these functions,
such as Theorem~\ref{thm:main_intro},
may facilitate some computations.
\medskip

Note that Demailly's argument above adapts to different geometric
settings. For example, we might consider algebraic or holomorphic
curves $\CC \to X$, or even formal ones (with formal series as
coordinates).
We may also assume that $X$ is a real differential manifold and
consider $\cC^\infty$ curves $\R \to X$. Taking $E$ to be respectively
the set of polynomial, entire, formal functions $\CC \to \CC$,
or $\cC^\infty(\R)$, in all cases vanishing at $0$,
we still have a local isomorphism between
$\cC_V$ and $\CC^r \otimes E$.

In \cite{CDM}, only jets of curves
to some order $\ell$ are considered, so that
the relevant space $E$ is the quotient ring $\CC[T]_0/(T^\ell)$.
Unfortunately, we have not been able to obtain an analogue of Theorem~\ref{thm:main_intro}
in the case $E=\CC[T]_0/(T^\ell)$; see Section~\ref{ssec:intro_TruncatedCase} for details.

\subsection{Proof strategy}
It is standard that a homogeneous algebraic functional of
degree $d$ on some space $E$
can be written as a multilinear function applied to
potentially repeated arguments.
This is known as {\em polarization}. This idea can be
applied to our setting,
and it turns out that Theorem~\ref{thm:main_intro} can be
deduced from an analogous statement
describing the multilinear functions $\varphi$ satisfying \eqref{e:reparametrage_intro}.
This reduction is relatively straight-forward, and the core of the proof of Theorem~\ref{thm:main_intro}
is the description of the multilinear solutions of \eqref{e:reparametrage_intro}.
\smallskip

To this end, we first prove that a multilinear solution $\varphi$ of \eqref{e:reparametrage_intro}
is necessarily of the form
\[\varphi(T^{x_1},\ldots,T^{x_n})=c(x_1,\ldots,x_n) T^{\sum_i x_i},\]
 for some function $c$ from $\Znp$ to the coefficient field $\CC$
  (see Lemma~\ref{l:form_multilinear_solution}).
 Moreover, the function $c$ should be
 \begin{description}
 \item[homogeneous] for all integers $x_1,\ldots,x_n \ge 1$, and $k \ge 1$, we have
 $c(kx_1,\ldots,kx_n)= c(x_1,\ldots,x_n)$;
 \item[polyharmonic] for all integers $x_1,\ldots,x_n \ge 1$, and $h \ge 1$,
we have 
\begin{equation}
\label{e:recurrence_c_intro}
\forall h \geq 1\;,\;\big ( \sum_i x_i \big ) \, c(x_1,\ldots,x_n)
= \sum_i \, x_i \, c(x_1,\ldots,x_i + h,\ldots,x_n).
\end{equation}
 \end{description}
 This translates our original algebraic problem into a more combinatorial one, which consists in solving 
 the above discrete functional equations.
 Recall that we can produce $n!$ solutions, namely the functions $\varphi_\sigma$ defined in 
 \eqref{e:def_varsigma_intro}, where $\sigma$ is a permutation in $S_n$ 
(that $\sigma$ is a permutation ensures that $\varphi_\sigma$ is multilinear, and not only algebraic).
Moreover, proving that these functions are linearly independent is not difficult.
Hence, what we need to show is 
 that the space of homogeneous polyharmonic functions $c:\Znp \to \CC$ 
 has dimension at most $n!$.
 \smallskip
 
A heuristic argument is the following, where we use the standard discrete difference operator $\delta_i$ defined by
\[(\delta_i c)(x_1,\dots,x_n)=c(x_1,\dots,x_i+1,\dots,x_n)-c(x_1,\dots,x_n).\]
Equation \eqref{e:recurrence_c_intro} is equivalent to the
fact that $c$ satisfies
$$\forall h \geq 1\;,\; D_h \cdot c = 0, \ \mbox{ where } \
D_h \coloneqq \sum x_i \delta_i^h.$$
A function $c$ \vupe{on $\Znp$ vanishes if and only} if all polynomials (symmetric or not) in $(\delta_i)_{1 \le i \le n}$ \vupe{applied to $c$ and evaluated at $(1,\dots,1)$} vanish.
\vf{When evaluated at $(1,\dots,1)$,
 the operator $D_h$ yields} the $h$-th power sum
of the discrete difference operators $(\delta_i)_{1 \le i \le n}$.
Since this holds for all $h\ge 1$,
all symmetric polynomials in $(\delta_i)_{1 \le i \le n}$ with zero constant term
applied to $c$ and evaluated at $(1,\dots,1)$ vanish.
Therefore the dimension of the space of polyharmonic functions should be equal to that of
the quotient space of all polynomials by the ideal generated
by symmetric ones with no constant term,
i.e.~the space of {\em symmetric coinvariants}.
This is a well-studied space, known in particular to be of dimension $n!$ \footnote{Much more is known on this coinvariant space: in particular, it is isomorphic as $S_n$-representation to the regular representation, and generalizations of this statement
in many directions have been investigated, see e.g.~\cite{coinvariant}
and references therein.}.
\smallskip

Applying directly this heuristic is however made difficult by the fact that the polyharmonic relation
involves a multiplication by $x_i$ and that $x_i$ and $\delta_i$ do not commute.
We could overcome this difficulty easily if we were working with usual differential operators instead of discrete difference ones, since the relation between $x_i$ and $\partial_i$ is simpler.
 The next idea is therefore to use a deformation argument based
on the following observation. The homogeneity condition
allows to extend any
solution $c$ from $\Zp^n$ to $\Q^n$, and the extended
solution satisfies the poyharmonic relation not only for
positive {\em integers} $h$, but for all positive {rational
numbers} $h$. We introduce an algebra
\[ \cD = \CC[H]\langle X_1,\ldots,X_n,D_1,\ldots,D_n \rangle, \]
that accounts for both the differential operators and the
discrete difference operators. Specializing to $H=0$,
we have the usual differential operators and we know that
the space of solutions has dimension $n!$.
Specializing to $H=h>0$, we have the discrete difference operators
on a lattice with mesh $h$; since we have extended our function $c$
from $\Zp^n$ to $\Q^n$, any rational $h>0$ is convenient for us.
We then argue with topological considerations that for $h$
sufficiently small the space of solutions on $(h\Zp)^n$
also has dimension at most $n!$, which concludes the proof.
\medskip

The above argument shows that the space of polyharmonic homogeneous functions
$c:\Zp^n \to \CC$ \vupe{has dimension at most $n!$}.
The heuristic comparison with the space of coinvariants suggests that
the space of polyharmonic functions $c$
(without assuming homogeneity) also has dimension $n!$.
By reversing some of the step of our proof,
we show that it is indeed the case,
so homogeneity is in fact be a consequence of polyharmonicity
(see Section~\ref{ss:homog-from-induction}).
Finding a direct proof of that seems however challenging,
and our argument crucially relies on the connection with
describing the solution of \eqref{e:reparametrage_intro}.

\begin{remark}\label{rmk:commutators_intro}
Via some commutator computations, it can be proven that a function $c$ is polyharmonic
if and only if it is {\em biharmonic}, in the sense that it satisfies~\eqref{e:recurrence_c_intro} for $h=1$ and $h=2$; see Remark~\ref{rmk:commutators} for details. Nevertheless,
our proof of Theorem~\ref{thm:main_intro} uses relation~\eqref{e:recurrence_c_intro} for all positive $h$.
\end{remark}

\begin{remark}
There is a large literature on positive ``harmonic'' functions on graded graphs;
see, e.g., \cite{HarmonicFunctions}.
In general, such functions form an infinite-dimensional convex set whose extremal points 
can be determined through various methods. 
Our equation~\eqref{e:recurrence_c_intro} for $h=1$ tells us that the function
$c$ is harmonic on the graded graph $\Zp^n$ with edge multiplicity $x_i$ on the edge going from $(x_1,\dots,x_n)$ to $(x_1,\dots,x_i+1,\dots,x_n)$.
However, our problem is fundamentally different
 from the usually considered question in the literature,
since: i) we are not only interested in positive functions,
and ii) the extra condition~\eqref{e:recurrence_c_intro} for $h=2$ reduces the space of solutions to a finite-dimensional space.
\end{remark}

\subsection{Further results and discussions}

\subsubsection{A continuous analogue of the polyharmonicity relations}
Our resolution of the induction relation \eqref{e:recurrence_c_intro}
also gives the solution to the analogue system of partial differential equation.
Namely, we prove that the space of distributions $c$ on $(0,+\infty)^n$ satisfying
\[
\forall k \in \N,\ \sum_{i=1}^n x_i \frac{\partial^k c}{\partial x_i^k} = 0
\]
has dimension $n!$, and we exhibit an explicit basis of solutions which are connected
to the functionals $\varphi_\sigma$. 
See Section~\ref{ss:system_EDP} for details.

\subsubsection{The truncated case}
\label{ssec:intro_TruncatedCase}
As explained in Section~\ref{sec:intro_motivation}, the most relevant case
for geometric applications in that of the space $E=\CC[T]_0/(T^\ell)$,
i.e.~truncated Taylor expansions up to degree $\ell$.
In this case, we believe that, as in the case $E=\CC[T]_0$,
the iterated semi-products are the only algebraic functionals commuting with reparametrizations.
We prove it for functionals of total degree 2 and 3 in $f_1,\dots,f_n$, 
but are not able to deal with the general case.

Again in this case, one can reduce the problem through polarization to the study of multilinear solutions, and the original algebraic problem can be recasted into discrete functional equations.
Namely one should now study polyharmonic functions $c$ defined on
the truncated space\footnote{The polyharmonic condition is only required when both sides of the identity are well-defined, i.e.~when $x_1+\dots+x_n+h \le \ell$.}
\[\Omega^n_\ell \coloneqq \{(x_1,\dots,x_n) \in \Zp^n: x_1+\dots+x_n < \ell\} \subset \Zp^n.\]
We have tried several approaches to bound from above the dimension of 
such polyharmonic functions.
A first attempt consists in identifying invertible subsystems, and works well for $n=2$
and $n=3$.
In a second approach, we look for extra annihilating operators for $c$ from which
we easily see that $c$ is determined by a finite number of initial conditions --
the (non-commutative) algebra generated by the multiplications by $x_i$ and the discrete difference operators $\delta_i$ is known as an ``Ore algebra'' and there are some efficient algorithms, such as Gröbner basis computation, to perform computation in such algebras, see, e.g.,~\cite{Ore,HolonomicFunctions}.
These approaches do not lead to a full solution of this truncated case,
but give evidence towards our conjecture that iterated semi-products are the only algebraic functionals commuting with reparametrizations.
Details are given in Section~\ref{sec:truncated}.

\subsubsection{A generalization to colored permutation}
As explained in the proof strategy section, the key point of the proof of Theorem~\ref{thm:main_intro}
is to establish that the space of multilinear solutions
of \eqref{e:reparametrage_intro}  has dimension $n!$.
Moreover, the basis elements are naturally indexed by permutations,
i.e.~by elements of Coxeter groups of type A.
A standard question in the field, when this kind of picture arises, is to construct
generalizations for other Coxeter groups, 
or more generally for complex reflection groups.

In Section~\ref{sec:typeB}, we provide a generalization of the main results to colored permutations.
(i.e.~to the complex reflection group $G(r,1,n)$).
This needs to introduce a slight variation of the general problem,
where \vupe{\eqref{e:reparametrage_intro} is only required
for reparametrization functions $\gamma$ living in a subspace $F$ of $E$}.
Fix $r\ge 1$, and let $\lambda$ be a primitive $r$-th root of unity.
Then take $E=\CC[T]_0$ as before, but let $F$ be the subspace of complex-valued 
polynomials $P$ satisfying $P(\lambda T)=\lambda P(T)$.
With this modification, multiplying  the variable by $\lambda$,
namely, setting $\varphi(f(t))=f( \lambda \, t)$, gives a valid solution for $n=1$.
Combining with the previous construction, 
we can construct multilinear solutions
\[\varphi_{\sigma,\bm \eps}(f_1,\dots,f_n)=[ f_{\sigma(1)}(\lambda^\eps_{1} T)\, :\, \dots\, :\, f_{\sigma(n)}(\lambda^\eps_{n} T)],\]
which are indexed by pairs $(\sigma,\bm \eps)$ where $\sigma$ is a permutation
and $\bm \eps$ a list in $(\Z/r\Z)^n$, i.e.~by $r$-colored permutations.
We shall prove (see Theorem~\ref{thm:type_B}) that these form a basis of the space of
multilinear functionals from $E^n$ to $E$, commuting with reparametrizations by elements
of $F$. The polarization argument goes through easily, giving a description of algebraic 
functionals on $E^n$, commuting with reparametrizations by elements of $F$.

This is a natural generalization of our main theorem to colored permutation.
We have not been able to find analogues for other complex reflection groups/Coxeter types.

\subsubsection{The case of infinitely differentiable functions through polynomial approximations}
In Section~\ref{sec:CInfty}, we look at the case where $E=C^\infty(\mathbb R)_0$ is
 the space of infinitely differentiable  functions
$\mathbb R\to\mathbb R$ vanishing at $0$.
Then the functions $\varphi_\sigma$ are still well-defined and solutions of \eqref{e:reparametrage_intro}, 
and the question is again whether they span the whole
space of algebraic/multilinear solutions.

Our results depend on the continuity assumption that we make on $\varphi$.
Without any continuity assumption, we cannot prove that $\varphi$ 
is a linear combination of $\varphi_\sigma$'s, but it must share some properties with such functions, e.g.~for $t_0>0$, the value $\varphi(f_1,\dots,f_n)(t_0)$ \pe{only} depends on the
restriction of $f_1$,\dots, $f_n$ to $[0,t_0]$ (Proposition~\ref{p:past}).
Under some weak continuity assumption (which we call path-continuity, see Section~\ref{ss:path-cont}), we can approximate any multilinear solution $\varphi$
by a linear combination $\psi=\sum a_\sigma \varphi_\sigma$, in the sense that 
 $\varphi(f_1,\dots,f_n)$ and $\psi(f_1,\dots,f_n)$ have the same Taylor expansion in $0$
 at all orders for all $f_1,\dots,f_n$. Finally under a stronger continuity assumption 
 (which we call global continuity, see Section~\ref{ss:global-continuity}), we prove that any multilinear  solution $\varphi$
is a linear combination $\psi=\sum a_\sigma \varphi_\sigma$.
All these results are based on Theorem~\ref{e:def_varsigma_intro}.

\subsubsection{Functions with finite range}
Since our approach relies on derivation (see Lemma~\ref{lem:derivative} below),
we decided to consider a simple setting with non-smooth functions,
namely the set $E$ of functions from $\R$ to $\R$ with finite range vanishing at $0$.
In this case the derivative, and hence the semi-product, are ill-defined (they are equal to $0$ on step functions, but general functions with finite range may not have finite variations).
We prove in Section~\ref{sec:piecewise} that the only {\em multilinear} solutions of~\eqref{e:reparametrage_intro} are the functionals $\varphi(f_1,\dots,f_n)=\kappa\, f_1 \cdots f_n$, where $\kappa$ is in $\R$.
Consequently (using a polarization argument as in Section~\ref{sec:algebraic}), the only {\em algebraic} solutions of~\eqref{e:reparametrage_intro} in this case are of the form 
{$\varphi(f_1,\dots,f_n)(z)=H(f_1(z),\dots,f_n(z))$ for some algebraic function $H$.}

\subsubsection{Non-algebraic solutions}
Theorem~\ref{thm:main_intro} describes the algebraic solutions of~\eqref{e:reparametrage_intro},
and we can ask about the existence and a description of non-algebraic solutions.
For $E=\CC[T]_0$, a simple non-algebraic solution, which is not of the form~\eqref{eq:space_transfo}
is given by
\[ \varphi(f_1,\dots,f_n) = \bm 1 [f_1=f_2= \dots =f_n] \, f_n,\]
where the condition $f_1=f_2=\dots=f_n$ in the indicator should be understood as equalities of functions
 (hence, $\varphi(f_1,\dots,f_n)(z)$ depends on the full functions $f_1$,\dots,$f_n$,
 and not only on $f_1(z)$, \dots, $f_n(z)$).

A simple case, where we can completely describe non-algebraic solutions,
is when $E^*$ is a group for composition. This happens, e.g.~if we consider the set 
 $E= \bigcup_{k \text{ odd}} \mathbb R^o[[T^{1/k}]]$, where $\mathbb R^o[[X]]$ is the set
of formal power series with only odd powers.
Then the values of $\varphi$ on vectors whose last nonzero coordinates is $1$ can be chosen freely,
and if $f_i \ne 0$, we have
\[ \varphi(f_1,\dots,f_i,0,\dots,0)= \varphi(f_1 \circ f_i^{<-1>},\dots,f_{i-1} \circ f_i^{<-1>},1,0,\dots,0) \circ f_i,\]
where $f^{<-1>}$ denotes the compositional inverse of $\varphi$. Details are given in Section~\ref{sec:everything_invertible}. This shows that, at least in this case,
we can build many nontrivial non-algebraic solutions of~\eqref{e:reparametrage_intro}.


{\bf Acknowledgement :}
Antoine Lemenant was initially part of the project.
We thank him for generously sharing his ideas.

{\bf Use of AI:} 
\vf{AI helped us to find the proof that polyharmonicity implies homogeneity,
given in Section~\ref{ss:homog-from-induction} in the Appendix.
It was also used to debug the cocalc worksheet in Section~\ref{sec:truncated}.
AI was {\bf not} used in any of the sections of the core of this paper, neither for proof development, nor for proofreading. The second author, DDSF, belongs to the AI-free caucus of
the \href{https://www.ahmath.org/members}{Association for Human Mathematics}.}

\section{Notations and Preliminaries}

\subsection{Statement of the problem and definition of some solutions}
As in the introduction, we consider a monoid $(E,\circ)$ that we think about
as a space of functions from some pointed set $(Z,0)$ to itself. For example:
\begin{itemize}
\item[$\bullet$] $E=\CC[T]_0 \coloneqq \{P \in \CC[T] \;|\;P(0)=0\}$;
\item[$\bullet$] $E=\CC[T]_0/(T^\ell) \coloneqq \{f \in \CC[T]/(T^\ell) \;|\;f(0)=0\}$;
\item[$\bullet$] $E=\CC[[T]]_0 \coloneqq \{f \in \CC[[T]] \;|\;f(0)=0\}$;
\item[$\bullet$] $E=\cC^\infty(\R,\CC)_0 \coloneqq \cC^\infty(\R,\CC)\cap \{f(0)=0\}.$
\end{itemize}

Recall that we consider the set of functionals $\varphi : E^n \to E$ commuting with reparametrizations, i.e.~such that
\begin{equation}
\label{e:reparametrage}
\forall F \in E^n, \quad \forall \gamma \in E, \quad \varphi(F\circ \gamma)=\varphi(F)\circ \gamma.
\end{equation}
The set of $\varphi : E^n \to E$ satisfying \eqref{e:reparametrage} will be denoted  by $I_n(E)$.
 If $E$ is a $\CC$-vector space, then the set of elements in $I_n(E)$
which are moreover $\CC$-multilinear will be denotedby  $\LI_n(E)$. 
Finally, in the case $E=\CC[T]_0$, the set of elements $\varphi$
of $I_n(E)$ such that
the coefficients of $\varphi(P_1,\ldots,P_n)$ are themselves polynomial expressions in the coefficients of the $P_i$'s will be
denoted by $\AI_n(E)$ and called \emph{algebraic}.

From now on, unless specified otherwise, 
we assume that we are in one of the four settings above.
In particular, derivative and anti-derivative (vanishing at $0$) are well-defined.
Recall from the introduction that we define
\begin{equation}
\label{e:semi-product}
[f_1:f_2] \coloneqq \int f_1'f_2 \quad \mbox{ and } \quad
	[f_1:f_2:\cdots:f_m] = [f_1:[f_2:\cdots:f_m]] = \int f_1' \left( \int f'_2 \left( \cdots \int f_{m-1}' f_m  \right)\right).
\end{equation}

More generally, we introduce the following functions:
\begin{definition}
\label{d:varphi-sigma}
Let $n,m>0$ and let $\sigma : \{1,\ldots,m\} \to \{1,\ldots,n\}$ be any function. We
let $\varphi_\sigma:E^n \to E$ be the map
\begin{equation}
\label{e:varphi-sigma}
(f_1,\ldots,f_n) \mapsto [f_{\sigma(1)}:f_{\sigma(2)}:\cdots:f_{\sigma(m)}].
\end{equation}
\end{definition}

\begin{proposition}
\label{prop:f1f2p}
For all $n,m \in \N$ and for all $\sigma : \{1,\ldots,m\} \to \{1,\ldots,n\}$,
the function $\varphi_\sigma$ is in $AI_n(E)$. Moreover, if $n=m$ and $\sigma$ is
bijective, then $\varphi_\sigma \in \LI_n(E)$.
\end{proposition}
\begin{proof}
We first prove that $\varphi_\sigma$ satisfies \eqref{e:reparametrage}.
Adding repetitions or changing the order of the arguments of a functional $\varphi$
preserves the commutation with reparametrization, so that
it suffices to consider the case $n=m$, and $\sigma$ is the identity.
Using induction, it is enough to check \eqref{e:reparametrage} for $n=m=2$ and $\sigma=\id$, namely for
$$
\varphi : E^2 \to E, (f_1,f_2) \mapsto [f_1:f_2]
$$
\vupe{But this can be done by a simple
 change of variables formula in an integral:
\[[f_1 \circ \gamma:f_2 \circ \gamma](t)=\int_0^t f'_1(\gamma(u))\, \gamma'(u)\, f_2(\gamma(u))\, \dd u=\int_0^{\gamma(t)}
f'_1(v) f_2(v) \dd v= [f_1:f_2](\gamma(t)).\]
Hence  $\varphi_\sigma$ satisfies \eqref{e:reparametrage} for all $n,m \ge 1$ and $\sigma : \{1,\ldots,m\} \to \{1,\ldots,n\}$.}
Moreover, it is clear that $\varphi_\sigma$ is algebraic and furthermore, if $\sigma$ is bijective, then $\varphi_\sigma$ is actually $\CC$-multilinear.
\end{proof}

\subsection{On the codimension of some linear subspaces in the polynomial algebra}

We consider the polynomial algebra $\CC[Y_1,\ldots,Y_n]$. For a non-negative integer $d$,
the subspace of polynomials of degree
at most $d$ will be denoted $\CC_{\leq d}[Y_1,\ldots,Y_n]$.
Given a linear subspace $V$, we say that the quotient $\CC[Y_1,\ldots,Y_n]/V$ has degree at most $d$ if 
\[V + \CC_{\leq d}[Y_1,\ldots,Y_n] = \CC[Y_1,\ldots,Y_n].\]
In this section,
 we provide bounds on the dimension and the degree of $\CC[Y_1,\ldots,Y_n]/V$
  for some linear subspaces $V \subset \CC[Y_1,\ldots,Y_n]$. 
We start with algebraic results 
  (Lemmas~\ref{l:CI-regular} and~\ref{l:dim-coinv}), of which we include proofs for completeness.

Recall that for a ring $A$ and $P_1,\ldots,P_r \in A$, we say that $(P_1,\ldots,P_r)$
is a regular sequence (or an $A$-sequence) when for each $i < r$,
the element $P_{i+1}$ is not a zero divisor in $A / \langle P_1, \dots, P_{i} \rangle$.

\begin{lemma}
\label{l:CI-regular}
Let $P_1,\dots,P_r\in\mathbb{C}[Y_1,\dots,Y_n]$ be homogeneous polynomials of positive
degree defining a complete intersection, i.e. such that
$\dim\; \{x \, | \,P_1(x)=\cdots=P_r(x)=0\}=n-r$. Then $(P_1,\dots,P_r)$ is a regular sequence.
\end{lemma}
\vupe{Informally, the idea is that adding a new polynomial in the family $P_i$ decreases the dimension of the common vanishing locus by at most $1$, and it does decrease if only if $P_i$ is not a zero diviser in $\mathbb{C}[Y_1,\dots,Y_n] / \langle P_1, \dots, P_{i-1} \rangle$.
The difficulty is that this is only true {\em locally}, so we need to work in {\em a local ring} to run this argument.
We recall that, if $A$ is a ring, $M$ an $A$ module, and $\mathfrak{m}$ a prime ideal of $A$, 
then the local ring $(M)_{\mathfrak m}$ is defined as the set of fraction $p/q$ with $p$ in $M$, $q \in A\setminus \mathfrak m$, 
with the identifications $p/q=p'/q'$ if there exists $\alpha \in A \setminus \mathfrak m$ such that $\alpha(q'p-qp')=0$. }

\begin{proof}
Write $R=\mathbb{C}[Y_1,\dots,Y_n]$, and let $\mathfrak{m}$ and $I$
be the ideals of $R$ generated by $(Y_1,\dots,Y_n)$ and $(P_1,\dots,P_r)$, respectively.
We consider the local ring $A=R_{\mathfrak{m}}$ of $R$ at $\mathfrak{m}$,
which is Cohen--Macaulay and has dimension $n$.
Since $I$ is homogeneous, its vanishing locus is a cone and
all its irreducible components pass through the origin, so
\[\dim A/IA=\dim R/I\vupe{=\dim \{x \, | \,P_1(x)=\cdots=P_r(x)=0\}}=n-r.\]
 By \cite[Thm.~17.4(i)]{Matsumura}, $\operatorname{ht}(IA)=r$. As the
$P_i$ are homogeneous of positive degree, they lie in $\mathfrak{m}A$, and
\cite[Thm.~17.4(iii), $(3)\Rightarrow(1)$]{Matsumura} shows that $(P_1,\dots,P_r)$ is an
$A$-sequence.

It remains to deduce that it is an $R$-sequence. For each $i$, let $I_i$ be the ideal generated by
$(P_1,\dots,P_i)$ in $R$. Multiplication by $P_{i+1}$
is a homogeneous endomorphism of the graded module $M_i:=R/I_i$, so its kernel
$K_i$ is a graded submodule of $M_i$. By exactness of localization, $(K_i)_{\mathfrak{m}}$
is the kernel of the multiplication by $P_{i+1}$ on $(M_i)_{\mathfrak{m}}$, which vanishes since
$(P_1,\dots,P_r)$ is an $A$-sequence.
\vupe{Assume for the sake of contradiction that $K_i \ne \{0\}$.
Since $K_i$ is a graded module, it must contain a nonzero homogeneous element $x$.
The image of $x$ in the localization $(K_i)_{\mathfrak{m}}$ must be $0$ (since $(K_i)_{\mathfrak{m}}=\{0\}$),
which means that there exists $\alpha \in R \setminus \mathfrak{m}$ such that $\alpha\, x=0$ in $K_i \subset M_i:=R/I_i$, i.e.~$\alpha\, x$ lies in $I_i$.
Since $I_i$ is a graded ideal, this implies that all homogenous components of $\alpha\, x$
lie in $I_i$; in particular, letting $\alpha_0$ be the constant term in $\alpha$,
the lowest degree component $\alpha_0 x$ of $\alpha\, x$ lies in $I_i$.
Since $\alpha$ is not in $\mathfrak{m}$, its constant term $\alpha_0$ is nonzero,
and we conclude that $x$ is in $I_i$. But we assumed that $x$ was nonzero in $K_i \subseteq M_i=R/I_i$,
reaching a contradiction. Thus $K_i =\{0\}$, i.e.~$P_{i+1}$ is not a zero diviser in $M_i$.
Since this holds for all $i \le r-1$, the sequence $(P_1,\dots,P_r)$ is regular.}
\end{proof}

\begin{remark}\label{rmk:CI-nonhomog}
The homogeneity assumption in Lemma~\ref{l:CI-regular} is essential: a non-homogeneous
complete intersection need not be a regular sequence. In $R=\mathbb{C}[x,y,z]$, the polynomials
\[
  P_1=y(1-x),\qquad P_2=z(1-x),\qquad P_3=x
\]
have common zero locus $\{(0,0,0)\}$, hence define a complete
intersection of dimension $0$; yet, in this order, they do not form a regular sequence, since
$P_2$ is a zero-divisor modulo $(P_1)$: indeed $z(1-x)\cdot y=y(1-x)\cdot z\in(P_1)$ while
$y\notin(P_1)$. After localizing at $\mathfrak{m}=(x,y,z)$ the factor $1-x$ becomes a unit,
so the sequence is regular in $R_{\mathfrak{m}}$; homogeneity allows one to
descend regularity from $R_{\mathfrak{m}}$ to $R$.
\end{remark}

\begin{lemma}
\label{l:dim-coinv}
For $k=1,\ldots,n$, let $P_k \in \CC[Y_1,\ldots,Y_n]$ be a homogeneous
polynomial of positive degree $d_k$. Assume that $(0,\ldots,0)$ is the
only vector $(y_1,\ldots,y_n)$ in $\CC^n$ such that for all $k$,
$P_k(y_1,\ldots,y_n)=0$. Let $I \subset \CC[Y_1,\ldots,Y_n]$ be the
ideal generated by the polynomials $P_k$. We have
$ \CC[Y_1,\ldots,Y_n]/I $
has dimension $\prod_k d_k$ and maximal degree $\sum_k (d_k-1)$. 
\end{lemma}
In the setting of the Lemma, the vanishing locus of $P_1,\ldots,P_n$ is a cone,
and therefore it has dimension $0$ if and only if it is restricted to the origin.
Therefore, the assumption that $(0,\ldots,0)$ is the
only vector in $\CC^n$ on which all $P_k$ vanish simultaneously
is equivalent to saying that $P_1,\dots,P_n$ define a complete intersection.
By Lemma \ref{l:CI-regular}, this implies that $(P_1,\dots,P_n)$ is a regular sequence.

To prove the lemma, we use the notion of Hilbert series. 
Recall that the Hilbert series of a graded vector space
$V=\bigoplus_{d \ge 0} V_d$ is $H_V(t)=\sum_{d \ge 0} \dim(V_d) t^d$.
The lemma can also be seen as a consequence of the Bézout theorem~\cite[8.4(3)]{fulton}.
\begin{proof}
For each $k \le n$, we let $R_k$ be a complement subspace of the ideal $\langle P_1, \dots, P_k \rangle$ in $\CC[Y_1,\ldots,Y_n]$.
Since $P_1,\dots,P_n$ define a complete intersection, they form a regular sequence.
 This implies that
\[ \langle P_1, \dots, P_{k} \rangle = \langle P_1, \dots, P_{k-1} \rangle\ \oplus\ P_k \cdot R_{k-1}. \]
Looking at the corresponding Hilbert series, we have
\[H_{\CC[Y_1,\ldots,Y_n]}(t) - H_{R_k}(t) =H_{\CC[Y_1,\ldots,Y_n]}(t) - H_{R_{k-1}}(t)+t^{d_k}\,  H_{R_{k-1}}(t),\]
i.e.~$H_{R_{k}}(t)=(1-t^{d_k})H_{R_{k-1}}(t)$. Since 
\[H_{R_{0}}(t)=H_{\CC[Y_1,\ldots,Y_n]}(t)=\frac1{(1-t)^n},\]
we get by an immediate induction
\[H_{R_{n}}(t)=\frac{\prod_{k \le n} (1-t^{d_k})}{(1-t)^n}=\prod_{k \le n} (1+t+ \dots + t^{d_k-1}).\] 
 Setting $t=1$ gives the dimension of $\CC[Y_1,\ldots,Y_n]/I$
 and looking at the maximal exponent of $t$ gives us its maximal degree.
\end{proof}

In the next section we will need to apply this lemma to the family $P_k=\sum_{i \le n} a_i \, Y_i^k$,
where the $a_i$'s are {\em positive} real parameters. In the next lemma, we prove that
they indeed define a complete intersection.
\begin{lemma}
\label{l:sol-triviale}
Let $(a_i) \in \Rnp$. The only solution in $\CC^n$ to the system of equations
\begin{equation}
\label{equa:coinv}
\sum_i a_i \, y_i^k = 0
\end{equation}
for $k=1,\ldots,n$ is $(y_1,\ldots,y_n)=(0,\ldots,0)$.
\end{lemma}
\begin{proof}
This proof was suggested to us by Matei Toma, to whom we are grateful.
Given $y_1,\ldots,y_n \in \CC$, we see \eqref{equa:coinv} as a system
of linear
equations in $(a_i)$, and we show inductively on $n$
that it has no solution in $\Rnp$ unless $(y_1,\ldots,y_n)=(0,\ldots,0)$.
It holds trivially for $n=1$.
For $n \ge 2$, note that
the determinant of the system is the Vandermonde determinant
$y_1 \cdots y_n \prod_{i<j} (y_i-y_j)$. If it
does not vanish, then the only solution is $a_1=\cdots=a_n=0$, proving the induction step.
If, say, $y_1=0$, then we get
the same equations in $a_2,\ldots,a_n$, and by induction they have no
solution in $\Rnmp$ \vupe{unless all $y_i$'s vanish}. If, say, $y_1=y_2$, then we get the
equations $(a_1+a_2)y_2^k+a_3y_3^k+\cdots+a_ny_n^k=0$, and again these
equations do not have solutions in $\Rnp$ by induction unless $y_2=\dots=y_n=0$.
This concludes the proof.
\end{proof}

Lemmas \ref{l:dim-coinv} and \ref{l:sol-triviale} imply:

\begin{proposition}
\label{p:coinv-A}
Let $a_1,\ldots,a_n \in \Rnp$.
The ideal $I \subset \CC[Y_1,\ldots,Y_n]$ generated by the polynomials
$P_k \coloneqq \sum_i a_i \, Y_i^k$ is such that
$\CC[Y_1,\ldots,Y_n]/I$ has dimension $n!$ and maximal degree $\binom{n}2 -1$.
\end{proposition}

\begin{remark}
In the special case where $P_k = \sum_i Y_i^k$, the ideal $I$ is generated
by symmetric functions in $Y_1,\dots,Y_n$ with zero constant coefficients.
Then the quotient $\CC[Y_1,\ldots,Y_n]/I$ is known as the
symmetric coinvariant ring, and
its maximal degree and dimension are well-known.
\end{remark}

We conclude this section by a short lemma allowing to replace homogeneous polynomials
by some non-homogeneous deformations.
It might be well-known, but we \vupe{could not locate it in the literature}.

\begin{lemma}
\label{lemm:smaller-deg-perturbation}
Let $S$ be an index set and for each $s \in S$, let $A_s,B_s$ be
polynomials in $\CC[Y_1,\ldots,Y_n]$. Assume that
\begin{itemize}
 \item each polynomial $A_s$ is homogeneous;
 \item each polynomial $B_s$ satisfies that $B_s-A_s$ is a polynomial
 of degree strictly less than $\deg(A_s)$;
 \item $\dim \CC[Y_1,\ldots,Y_n]/\sum_{s \in S} \CC A_s < \infty$.
\end{itemize}
Then, $\CC[Y_1,\ldots,Y_n]/\sum_{s \in S} \CC B_s$ has dimension and maximal degree
smaller or equal to that of $\CC[Y_1,\ldots,Y_n]/\sum_{s \in S} \CC A_s$.
\end{lemma}
\begin{proof}
Denote by $\CC_d[Y_1,\ldots,Y_n]$ the space of homogeneous polynomials
of degree $d$. We choose a subspace $V \subset \CC[Y_1,\ldots,Y_n]$
such that $\dim V=\dim \CC[Y_1,\ldots,Y_n]/\sum_s \CC A_i$ and
$V \oplus \sum_{s} \CC A_s = \CC[Y_1,\ldots,Y_n]$. Moreover
we assume that $V = \oplus_d V_d$ with
$V_d \subset \CC_d[Y_1,\ldots,Y_n]$.
We let $S_d \subset S$ be
such that $\deg(A_s)=d$ for all $s \in S_d$, and we let
$S_{\leq d}=S_0 \cup S_1 \cup \cdots \cup S_d$. We use the
similar notations
$V_{\leq d}=V_0 \oplus \cdots \oplus V_d$ and
$\CC_{\leq d}[Y_1,\ldots,Y_n]$.
We prove by
induction on $d$ that
$V_{\leq d} + \sum_{s \in S_{\leq d}} \CC B_s = \CC_{\leq d}[Y_1,\ldots,Y_n]$, starting with
$d=0$.

Recall that by assumption, for $s \in S_{d+1}$,
$A_s$ is equal to $B_s$ modulo
$\CC_{\leq d}[Y_1,\ldots,Y_n]$. Therefore,
assuming
$V_{\leq d} + \sum_{s \in S_{\leq d}} \CC B_s = \CC_{\leq d}[Y_1,\ldots,Y_n]$
for some integer $d$, we get
\begin{multline}\label{eq:sum_V}
V_{\leq d+1} + \sum_{s \in S_{\leq d+1}} \CC B_s =
\Big (V_{\leq d} + \sum_{s \in S_{\leq d}} \CC B_s \Big) +
\Big (V_{d+1}+ \sum_{s \in S_{d+1}} \CC B_s \Big ) \\
= \CC_{\leq d}[Y_1,\ldots,Y_n] + \Big (V_{d+1}+
\sum_{s \in S_{d+1}} \CC A_s \Big ) 
= \CC_{\leq d+1}[Y_1,\ldots,Y_n].
\end{multline}
Since we proved that $V + \sum_{s \in S} \CC B_s=\CC[Y_1,\ldots,Y_n]$,
the codimension of $\sum_{s \in S} \CC B_s$ in $\CC[Y_1,\ldots,Y_n]$
is at most $\dim(V)$.

For the degree bound, if $d_0$ is the maximal degree of $\CC[Y_1,\ldots,Y_n]/\sum_{s \in S} \CC A_s$, one has $V_d=\{0\}$ for $d>d_0$.
Summing~\eqref{eq:sum_V} for all $d$ gives us
\[V_{\leq d_0} + \sum_{s \in S} \CC B_s =   \CC[Y_1,\ldots,Y_n],  \]
which shows that $\CC[Y_1,\ldots,Y_n]/\sum_{s \in S} \CC B_s$ has maximal degree at most $d_0$.
\end{proof}

\subsection{Geometric interpretation of iterated semi-products}
\label{sec:geometric_int_semi_prod}
In this section, we show that under extra assumptions of $f_1,\dots,f_n$,
the (iterated) semi-products $[f_1 : \dots: f_n]$ correspond to the 
volume (or area for $n=2$) of some region of the space.
The commutation with reparametrizations is straight-forward from this description.
This is independent from the rest of the paper.

\begin{figure}[htbp]
    \centering
    \includegraphics[width=0.75\linewidth]{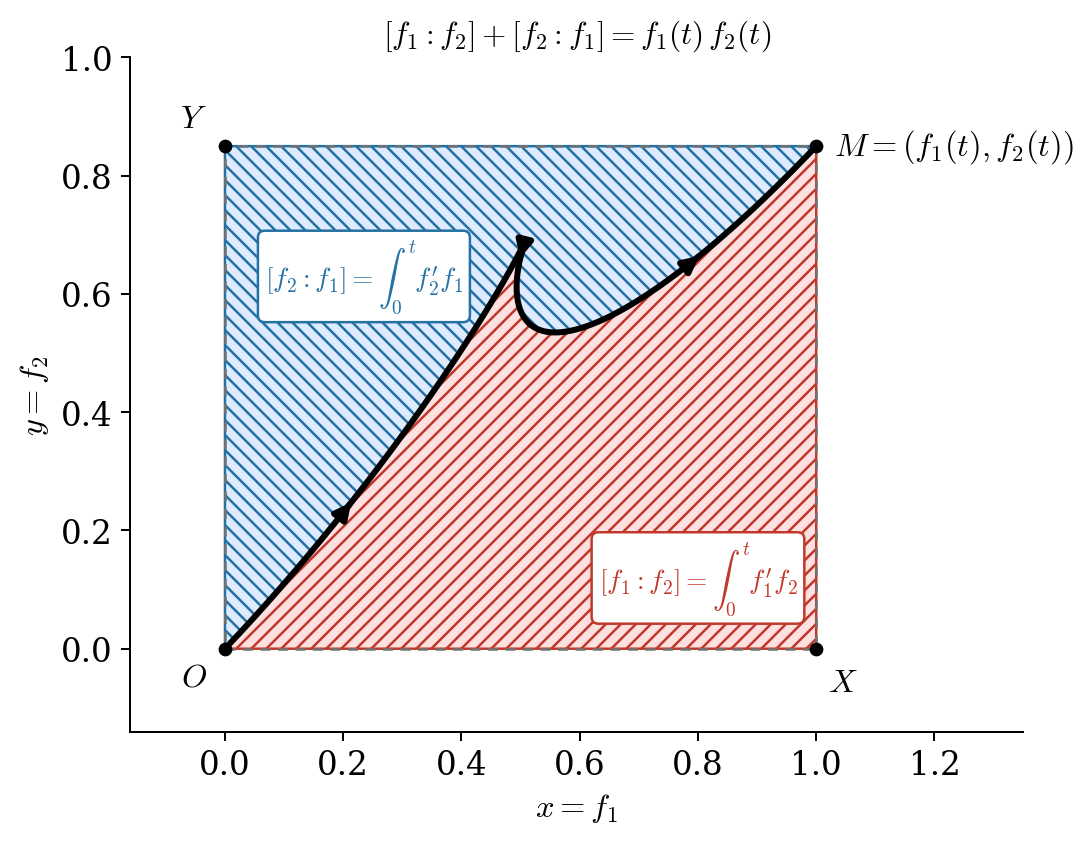}
    \caption{$[f_1{:}f_2]+[f_2{:}f_1]=f_1(t)f_2(t)$, for a path with neither
$f_1$ nor $f_2$ monotone. Red $=\int_0^t f_1'f_2$ (below--right), blue
$=\int_0^t f_2'f_1$ (above--left).}
    \label{fig:semiproduct_area}
\end{figure}

We start with the case of a single semi-product, i.e.~$n=2$.
The following discussion is illustrated on Figure \ref{fig:semiproduct_area}.
Let $t>0$ be fixed. Consider the planar path $\Gamma\colon s\mapsto\bigl(f_1(s),f_2(s)\bigr)$,
$s\in[0,t]$, from the origin $O=(0,0)$ to $M=\bigl(f_1(t),f_2(t)\bigr)$ (recall
$f_i(0)=0$). Suppose $\Gamma$ stays in the interior of the rectangle
$\mathcal R=[0,f_1(t)]\times[0,f_2(t)]$, i.e. $0 < f_i(s) < f_i(t)$ on $]0,t[$, and that
the map $\Gamma$ is injective;
\emph{no monotonicity of $f_1$ or $f_2$ is assumed}. Writing $x=f_1$, $y=f_2$,
\[
  [f_1:f_2](t)=\int_0^t f_1'f_2\,ds=\int_\Gamma y\,dx .
\]
Close $\Gamma$ into a loop by the two segments $M\to X$ and $X\to O$, where
$X=\bigl(f_1(t),0\bigr)$: along $M\to X$ one has $x\equiv f_1(t)$, hence $dx=0$;
along $X\to O$ one has $y\equiv 0$; so both contribute $0$ and
$\oint y\,dx=\int_\Gamma y\,dx$. By Green's theorem $\oint y\,dx$ is the
area enclosed, whence
\[
  [f_1:f_2](t)=\text{area of the region of }\mathcal R\text{ lying below--right of }\Gamma
  \quad(\text{the side of the corner }X;\ \text{red}).
\]
Symmetrically, $[f_2:f_1](t)=\int_0^t f_2'f_1\,ds=\int_\Gamma x\,dy$, and closing
through $Y=\bigl(0,f_2(t)\bigr)$ (where $dy=0$ on $M\to Y$ and $x=0$ on $Y\to O$)
gives
\[
  [f_2:f_1](t)=\text{area of the region above--left of }\Gamma
  \quad(\text{the side of the corner }Y;\ \text{blue}).
\]
The two regions tile $\mathcal R$, \vupe{which gives a geometric interpretation
of the decomposition of a product as semi-products:}
\[
  [f_1:f_2](t)+[f_2:f_1](t)=\operatorname{area}(\mathcal R)=f_1(t)\,f_2(t).
\]
\vupe{The invariance by increasing reparametrizations is also straight-forward on this interpretation,
since the image of $\Gamma$ is unchanged when both $f_1$ and $f_2$ are composed with some $\gamma$.}


\vupe{For general $n$, we do not have such a geometric interpretation,
but provide an integral expression of the iterated semi-product,
which makes apparent both the decomposition of products and the invariance by increasing reparametrizations.}
\begin{lemma}\label{l:iterated-simplex}
  For all $f_1,\ldots,f_n\in E$ and all $t$,
  \begin{equation}\label{e:iterated-simplex}
    [f_1:\cdots:f_n](t)
    =\int_{t\geq s_1\geq\cdots\geq s_n\geq 0}
      f_1'(s_1)\,f_2'(s_2)\cdots f_n'(s_n)\,ds_1\cdots ds_n,
  \end{equation}
  the integral being taken over the simplex $\Delta_t  \coloneqq  \{t\geq s_1\geq\cdots\geq s_n\geq 0\}$.
\end{lemma}
\begin{proof}
  We argue by induction on $n$, the case $n=1$ being the identity
  $[f_1](t)=f_1(t)=\int_0^t f_1'(s_1)\,ds_1$ (recall $f_1(0)=0$). For $n\geq 2$,
  the definition~\eqref{e:semi-product} of the iterated semi-product gives
  \[
    [f_1:\cdots:f_n](t)=\bigl[f_1:[f_2:\cdots:f_n]\bigr](t)
    =\int_0^t f_1'(s_1)\,[f_2:\cdots:f_n](s_1)\,ds_1 .
  \]
  By the induction hypothesis applied at the endpoint $s_1$,
  \[
    [f_2:\cdots:f_n](s_1)
    =\int_{s_1\geq s_2\geq\cdots\geq s_n\geq 0}
      f_2'(s_2)\cdots f_n'(s_n)\,ds_2\cdots ds_n .
  \]
  Substituting this into the previous integral and observing that the two
  conditions $0\leq s_1\leq t$ and $s_1\geq s_2\geq\cdots\geq s_n\geq 0$ describe
  exactly the simplex $\Delta_t$ yields
  \eqref{e:iterated-simplex}.
\end{proof}

\section{The dimension of $\LI_n(\CC[T]_0)$}
\label{sec:dimension_n!}

In this section, we describe the mutlilinear solutions of \eqref{e:reparametrage} 
on the space $\CC[T]_0$ of polynomials without constant terms.
In this case, it turns out that the solutions $\varphi_\sigma$
exhibited in Proposition~\ref{prop:f1f2p}
are linearly independent and span the space of solutions.

\subsection{Independence of iterated semi-products}
\label{s:independence_semi_products}
The following functions $c_\sigma$ will be useful in the proof of
\vupe{the independence of the $\varphi_\sigma$'s}.
\begin{definition}
\label{d:c-sigma}
Let $\sigma:\set{n} \to \set{n}$ be a bijection. For $x_1,\ldots,x_n \in \R$, we let
\begin{equation}
\label{e:c-sigma}
c_\sigma(x_1,\ldots,x_n) \coloneqq \frac{x_{\sigma(1)} \cdots x_{\sigma(n)}}{(x_{\sigma(1)}+\cdots+x_{\sigma(n)})\cdots
(x_{\sigma(n-1)}+x_{\sigma(n)})x_{\sigma(n)}}.
\end{equation}
\end{definition}

\begin{lemma}
\label{l:c-sigma_ind}
The family $(c_\sigma)_{\sigma \in \fS_n}$ of rational functions on $x_1,\ldots,x_n$ is linearly independent.
\end{lemma}
\begin{proof}
Let $(a_\sigma) \in \CC^{\fS_n}$ be such that $\sum_\sigma a_\sigma c_\sigma=0$. Then
\[
\sum_\sigma \frac{a_\sigma}{(x_{\sigma(1)}+\cdots+x_{\sigma(n)})\cdots(x_{\sigma(n-1)}
+x_{\sigma(n)})x_{\sigma(n)}}=0.
\]
Multiplying by $x_k$ and setting $x_k=0$ we get the identities
$$
\forall k \in \{1,\ldots,n\}, \quad \sum_{\sigma:\sigma(n)=k} \frac{a_\sigma}
{(x_{\sigma(1)}+\cdots+x_{\sigma(n-1)}) \cdots (x_{\sigma(n-2)}+x_{\sigma(n-1)})x_{\sigma(n-1)}}=0.
$$

By induction, we get $a_\sigma=0$ for all $\sigma \in \fS_n$.
\end{proof}

\begin{lemma}
\label{l:phi-sigma_ind}
The set of elements $\varphi_\sigma$ for $\sigma \in \fS_n$ form a free family of elements of the $\CC$-vector space
$\LI_n(E)$.
\end{lemma}
\begin{proof}
A direct computation shows that for $\sigma:\set{n} \to \set{n}$,
\begin{equation}
\label{value-bracket}
\varphi_\sigma(T^{i_1},\dots,T^{i_n})=[T^{\sigma(i_1)}:\cdots:T^{\sigma(i_n)}] = c_\sigma(i_1,\ldots,i_n)\, T^{\sum i_{\sigma(k)}}.
\end{equation}

Assuming we have a linear relation $\sum_{\sigma \in \fS_n} a_\sigma \varphi_\sigma = 0$,
applying it to $(T^{i_1},\ldots,T^{i_n})$ we find
$$
\sum_\sigma a_\sigma \, c_\sigma(i_1,\ldots,i_n)=0\,.
$$
This rational function vanishes on $\N^n$, thus it also vanishes on $\R^n$. We get
$\sum a_\sigma c_\sigma \equiv 0$, and thus $a_\sigma=0$ for all $\sigma \in \fS_n$, by Lemma
\ref{l:c-sigma_ind}.
\end{proof}

\subsection{Reduction to a system of discrete partial differential equations}
\vupe{We now take $\varphi$ in $\LI_n(\CC[T]_0)$ and want to prove that it is a linear combination of $\varphi_\sigma$'s.}
We start by an easy lemma on the general form of $\varphi$.
\begin{lemma}\label{l:form_multilinear_solution}
Let $\varphi \in \LI_n(\CC[T]_0)$. Then there exists a function $c:\Znp \to \CC$ such that
\begin{equation}
\varphi(T^{x_1},\ldots,T^{x_n})=c(x_1,\ldots,x_n) T^{\sum_i x_i}.
\end{equation}
\end{lemma}
\begin{proof} We use the condition \eqref{e:reparametrage} with $F=(T^{x_1},\ldots,T^{x_n})$ and $\gamma(t)=\tau t$, $\tau \in \R$. We obtain
$$
\varphi(T^{x_1},\ldots,T^{x_n})(\tau T) =\varphi(\tau^{x_1} T^{x_1},\ldots,\tau^{x_n} T^{x_n})
=\tau^{\sum x_k} \varphi(T^{x_1},\ldots,T^{x_n}).
$$
Therefore, $\varphi(T^{x_1},\ldots,T^{x_n})$ is homogeneous of degree $\sum x_i$, proving the result.
\end{proof}
The next step is to translate condition \eqref{e:reparametrage} as constraints
on the function $c$ associated with $\varphi$.
To this end, it is convenient to differentiate \eqref{e:reparametrage} with respect to $\gamma$.
\begin{lemma}[Derivative]
\label{lem:derivative}
Let $\varphi \in \LI_n(\CC[T]_0)$. Then, for any
$F=(f_1,\ldots,f_n) \in {(\CC[T]_0)}^n$, and for any $\eta \in \CC[T]_0$,
we have
\begin{eqnarray}
\label{identiteDer}
\eta \cdot [\varphi(F)]'=\sum_{i=1}^n \varphi(f_1,\ldots,\eta f_i',\ldots,f_n).
\end{eqnarray}
\end{lemma}

\begin{proof}
Let us take $F=(f_1,\ldots,f_n)$ in ${(\CC[T]_0)}^n$ and expand it as
$(\sum_{k=1}^d a_{k,i} T^k)_{1 \le i \le n}$, where the $a_{k,i}$'s are complex numbers. 
Given any $h \in \CC$ and $\eta \in \CC[T]_0$, we have
$$
F(h\eta(T)+T)=\left( \sum_{k=0}^d a_{k,i} \sum_{j=0}^k h^{j} \binom{k}{j} \eta(T)^{j}T^{k-j}\right)_{1  \le i \le n}.
$$
Multilinearity of $\varphi$ yields that
$$
\varphi(F(h\eta(T)+T)) = \sum_{(k_1,\ldots,k_n)} \sum_{(j_1,\ldots,j_n)} \prod_{i=1}^n a_{k_i,i} \,  h^{\sum j_i}
\prod_i \binom{k_i}{j_i} \varphi \big (\eta(T)^{j_1} T^{k_1-j_1},\ldots,\eta(T)^{j_n} T^{k_n-j_n} \big ).
$$
Here $k_1,\ldots,k_n$ vary from $0$ to $d$, and $j_i$ varies from $0$ to $k_i$.
We want to extract the coefficient of $h$ in the above equality.
In the right-hand side, only summands where one of the $j_i$'s, say $j_u$, is equal to $1$
and the other are $0$, contribute.
Hence we get
$$
\frac{\partial}{\partial h}_{|h=0} \varphi(F(T+h\eta(T))) =
\sum_{(k_1,\ldots,k_n)} \prod_i a_{k_i,i} \sum_{u=1}^n
k_u \varphi \big (T^{k_1},\ldots,\eta(T) T^{k_u-1},\ldots,T^{k_n} \big ).
$$
Using again multilinearity, this rewrites as
\begin{align*}
\frac{\partial}{\partial h}_{|h=0} \varphi(F(T+h\eta(T))) &=
\sum_{u=1}^n \varphi \big (
\sum_{k_1} a_{k_1,1}T^{k_1},\ldots,\eta(T) \sum_{k_u} a_{k_u,u} k_uT^{k_u-1},\ldots,\sum_{k_n} a_{k_n,n} T^{k_n} \big ),\\
& = \sum_{u=1}^n \varphi \big (f_1,\ldots,\eta(T) f_u',\ldots,f_n \big ).
\end{align*}

On the other hand, \eqref{e:reparametrage} implies that
$\varphi(F(T+h\eta(T))) = \varphi(F)(T+h\eta(T))$. The linear term in $h$ in this quantity is
$\eta(T)[\varphi(F)]'$, yielding \eqref{identiteDer}.
\end{proof}

\begin{proposition}[Polyharmonicity and homogeneity]
\label{p:induction}
Let $\varphi : \CC[T]_0^n \to \CC[T]_0$ be a fonctional of the form
\begin{equation}
\label{e:def-c}
\varphi(T^{x_1},\ldots,T^{x_n})=c(x_1,\ldots,x_n)\, T^{\sum_i x_i}
\mbox{ for }x_1,\ldots,x_n \in \Zp,
\end{equation}
where $c$ is a function from $\Zp^n$ to $\CC$.
If $\varphi \in \LI_n(\CC[T]_0)$, then, for all $x_1,\ldots,x_n \in \N$ and all $k \in \Zp$
we have the relations:
\begin{align}
\label{recurrence-c}
\big ( \sum_i x_i \big ) \, c(x_1,\ldots,x_n)
&= \sum_i \, x_i \, c(x_1,\ldots,x_i + k,\ldots,x_n);\\
\label{c-homogene}
c(k x_1,\ldots,k x_n)&=c(x_1,\ldots,x_n).
\end{align}
\end{proposition}
\begin{proof}
Let $F=(f_1,\dots,f_n)=(T^{x_1},\ldots,T^{x_n})$ and $\eta=T^{k+1}$.
Relation \eqref{e:def-c} implies
\begin{equation}
\label{amontre1}
 \eta \cdot \varphi(F)'
=\big ( \sum x_i \big )\,  c(x_1,\ldots,x_n)T^{k+\sum x_i}.
\end{equation}
On the other hand, using also \eqref{e:def-c}, we have
\begin{align}
 \sum_{i=1}^n \varphi(f_1,\ldots,\eta f_i',\ldots,f_n)=&\sum_i \varphi(T^{x_1},\ldots,T^{k+1} x_i T^{x_i-1},\ldots,T^{x_n})
\nonumber \\
=& \sum_i \, x_i \, c(x_1,\ldots,x_i + k,\ldots,x_n) T^{k+\sum_i x_i }
\label{identiteDerivee2}
\end{align}
Then relation \eqref{identiteDer} gives the equality of the terms computed in \eqref{amontre1} and in
\eqref{identiteDerivee2}, which implies \eqref{recurrence-c}.

To prove \eqref{c-homogene}, we directly apply \eqref{e:reparametrage}
with $F=(f_1,\dots,f_n)=(T^{x_1},\ldots,T^{x_n})$ and $\gamma=T^{k}$. This gives
\[ c(k x_1,\ldots,k x_n)\, T^{k\sum_i x_i} = \varphi(T^{k\, x_1},\ldots,T^{k\, x_n})
= \varphi(T^{x_1},\ldots,T^{ x_n}) \circ T^k
= c(x_1,\ldots,x_n)\, T^{k\sum_i x_i},
\]
implying \eqref{c-homogene}, as claimed.
\end{proof}

\begin{remark}\label{rmk:commutators}
Recall that $\delta_i$ denotes the discrete difference operator defined by
\[(\delta_i c)(x_1,\dots,x_n)=c(x_1,\dots,x_i+1,\dots,x_n)-c(x_1,\dots,x_n).\]
If we consider the operator $\Delta_k=\sum_{i=1}^n x_i ((\delta_i+\Id)^k - \Id)$, then \eqref{recurrence-c}
rewrites as $\Delta_kc \equiv 0$.
A short computation with commutators using $[\delta_i,x_i]=\delta_i+\Id$ shows that
\[ [\Delta_k,\Delta_\ell]=(k-\ell)\Delta_{k+\ell} -k\Delta_k + \ell \Delta_\ell.\]
In particular for $k \ge 2$, we have $[\Delta_k,\Delta_1]=(k-1) \Delta_{k+1} -k\Delta_k + \Delta_1$.
 Consequently, if $c$ is such that $\Delta_1 c =\Delta_2 c=0$,
an immediate induction proves that $\Delta_k c=0$ for all $k \ge 1$.

In other terms, asking that Equation~\eqref{recurrence-c} is valid for all $k \ge 1$
is equivalent to asking it only for $k=1$ and $k=2$.
Nevertheless, the version with all $k\ge 1$ will be useful below to find an upper bound
on the dimension of the space of solutions.
\end{remark}

\subsection{Generalization and difference operator formulation}

We are thus led to study the homogeneous solutions of \eqref{recurrence-c}. We state a result that
bounds the dimension of the space of such solutions, with a slightly higher degree of
generality, that will be useful in Section \ref{sec:typeB}. We choose a positive
integer $r$.
To prove our main Theorem \ref{t:linear}, the choice $r=1$ is enough, but for
Theorem \ref{thm:type_B}, all possible values of $r$ are needed.
For $x \in \Z$, $[x]$ denotes its class in $\Z/r\Z$.
Let $\beps = (\eps_i) \in (\Z/r\Z)^n$ and let $\Ze \subset \Znp$ be the set of $n$-uples
$(x_i)$ such that $[x_i]=\eps_i$.
Let $\Zrp$ denote the set of positive integers $n$ with $[n]=1$.
We also introduce
a weak version of homogeneity and polyharmonicity:

\begin{definition}
\label{d:r-homogene}
A function $c:\Ze \to \CC$:
\begin{itemize} \item is called $r$-homogeneous when for all  $(x_i) \in \Ze$, and all $\lambda \in \Zrp$,
\begin{equation}
\label{e:homo-r}
c(\lambda x_1,\ldots,\lambda x_n)=c(x_1,\ldots,x_n);
\end{equation}
\item is called $r$-polyharmonic when for all $(x_i) \in \Ze$, and all $\forall k \in \Zp$,
\begin{equation}
\label{recurrence-cr}
 \big ( \sum_i x_i \big ) \, c(x_1,\ldots,x_n)
= \sum_i \, x_i \, c(x_1,\ldots,x_i + kr,\ldots,x_n)
\end{equation}
\end{itemize}
\end{definition}
Note that for $(x_i) \in \Ze$, $\lambda \in \Zrp$ and $k \in \Zp$, 
we have $[\lambda x_i]=[x_i+kr]=\eps_i$, so that \eqref{e:homo-r} and \eqref{recurrence-cr} make sense.
Our goal, that will be achieved in Section \ref{s:conclusion}, is the following:
\begin{proposition}
\label{p:bound}
Let $r$ and $\beps$ be fixed.
The vector space of $r$-homogeneous and $r$-polyharmonic functions $c:\Ze \to \CC$
has dimension at most $n!$\,.
\end{proposition}

To prove this Proposition, we consider a $r$-homogeneous and $r$-polyharmonic function $c$.
First, we use $r$-homogeneity to extend $c$ to a function
defined over rational numbers.

\begin{definition}
\label{d:qr}
Let $\Qr$ be the set of fractions $p/q$ where $[q]=1$. It is a subring of $\Q$. Let $\Qrp \subset \Qr$
be the subset of positive elements. Let $\Qe \subset \Qnp$ be the set of $n$-uples
$(x_i=p_i/q_i)_{1 \leq i \leq n}$
with  $[p_i] = \eps_i$ and $[q_i]=1$.
\end{definition}
In case $r=1$, $\Ze = \Znp$ and $\Qe=\Qnp$.
In general,
our function $c$ defines a function $f:\Qe \to \CC$ as follows.
Given $(x_1,\ldots,x_n) \in \Qe$, there exist $Q$ in $\Zrp$ such that $(Qx_1,\ldots,Q x_n) \in \Ze$.
Then we set
\begin{equation}
\label{def-f}
f\left (x_1,\ldots,x_n \right )=c(Qx_1,\ldots,Q x_n)
\end{equation}
By \eqref{e:homo-r}, the right-hand side does not depend on the choice of $Q$, and $f$ is well-defined.
The function $f$ defined this way inherits the polyharmonicity property of $c$.
To state this formally, let us introduce some generalized discrete difference operators.
\begin{notation}
\label{delta_ih}
Let $\Fe$ denote the set of functions $\Qe \to \CC$.
Let $h \in r\,\Qrp$.
We define $\delta_{i,h}:\Fe \to \Fe$ by
\[(\delta_{i,h} f)(x_1,\dots,x_n)=
\frac{f(x_1,\dots,x_i+h,\dots,x_n)-f(x_1,\dots,x_n)}{h/r}.\]
We let $\delta_i \coloneqq \delta_{i,r}$, which can also be applied to functions $\Ze \to \CC$.
\end{notation}

\begin{lemma}
\label{lem:discrete_diff_cf}
Let $c:\Ze \to \CC$ and $f:\Qe \to \CC$ be as above. Then we have
\begin{align}
\label{equa:discrete_diff_c}
\forall h \in r\Zp, \forall k \in \Zp\ ,\quad & \sum_{i=1}^n x_i \, \delta_{i,h}^k c = 0;\\
\label{e:f}
\forall h \in r\,\Qrp,\ \forall k \in \Zp,\ \quad &
\sum_{i=1}^n x_i \, \delta_{i,h}^kf = 0 \in \Fe\,.
\end{align}
\end{lemma}
\begin{proof}
We start with \eqref{equa:discrete_diff_c}. For this, fix a positive multiple $h$ of $r$ and
 observe that
  \[c(x_1,\ldots,x_i + kh,\ldots,x_n)= (r/h\, \delta_{i,h}+\Id)^k c(x_1,\ldots,x_k,\ldots,x_n).\]
  Thus \eqref{recurrence-cr} with $k$ substituted by $kh/r$ rewrites as
  \[\sum_{i=1}^n x_i \left( (r/h\, \delta_{i,h}+\Id)^k - \Id \right)\, c(x_1,\ldots,x_k,\ldots,x_n) =0.\]
  An immediate induction from here shows \eqref{equa:discrete_diff_c} for all $k>0$ and all $h$ in  $r\Zp$.
  Equation \eqref{e:f} follows by homogeneity.
\end{proof}

\subsection{Discrete differential equations defining a complete intersection}
\label{sec:upper_bound_general}

To deal with equation \eqref{e:f}, we consider a more
general setting. Namely, we introduce the algebra $\cD$ defined as follows:
\[ \cD = \CC[H]\langle X_1,\ldots,X_n,D_1,\ldots,D_n \rangle, \]
where 
\begin{itemize}
\item the $X_i$'s commute with each other;
\item the $D_i$'s commute with each other;
\item the $X_i$'s and  the $D_i$'s satisfy the following relation
\begin{equation}
\label{commutationXD}
\left \{
\begin{array}{l}
D_i X_j = X_j D_i \quad \mbox{if} \ \ i \neq j ; \\
D_i X_i = (X_i+H) D_i + r.
\end{array}
\right.
\end{equation}
\end{itemize}

Specializing $H$ to $h \in r \, \Qrp$,
the algebra $\cD$ is isomorphic to the algebra
generated by the operators
$\Fe \to \Fe$ of multiplication by $X_i$ and of
divided difference $\delta_{i,h}$ defined in Notation \ref{delta_ih}.
Indeed, the relation $\delta_{i,h} X_i = (X_i+h) \delta_{i,h}  + r$ is straightforward to check
(when $r=h=1$, it is sometimes referred to as Leibniz rule for divided differences).
Specializing $H$ to $0$ \vf{and $r$ to $1$}, we get the standard algebra of differential operators on $n$
variables with the usual Leibniz rule.

Note that each element $P \in \cD$ has a unique representation as
\begin{equation}
\label{conv-P}
P = \sum_{a_1,\dots,a_n \ge 0} Q_{a_1,\dots,a_n}(X_1,\ldots,X_n,H)\, D_1^{a_1} \cdots D_n^{a_n}.
\end{equation}
 We define the degree of such an element $P$ as $a_1+\dots+a_n$, i.e.~we only consider
the degree in the $D_i$'s in this canonical form throughout the section.
With this notation, given $h \in r\,\Qrp$, we set
\begin{equation}
P(\delta_{1,h},\ldots,\delta_{n,h})=
\sum_{a_1,\dots,a_n \ge 0}  Q_{a_1,\dots,a_n} (X_1,\ldots,X_n,h)\, \delta_{1,h}^{a_1}\cdots \delta_{n,h}^{a_n}.
\end{equation}
It is an operator $\Fe \to \Fe$.
Furthermore,
for $x=(x_1,\ldots,x_n) \in \Q^n$ and $h \in \Q$, we set
\begin{equation}
\ev_{x,h}(P) = \sum_{a_1,\dots,a_n \ge 0} Q_{a_1,\dots,a_n}(x_1,\ldots,x_n,h)
D_1^{a_1} \dots D_n^{a_n} \in \CC[D_1,\ldots,D_n].
\end{equation}
The map $\ev_{x,h}$ is a linear operator from $\cD$ to $\CC[D_1,\ldots,D_n]$,
but not an algebra morphism.
Directly from the definition, we observe that the following compatibility relation holds for 
all $P \in \cD$, $x \in \Qe$ and $h \in r\,\Qrp$:
\begin{equation}\label{equa-ev}
P(\delta_{1,h},\ldots,\delta_{n,h})(f)(x) =
\ev_{x,h}(P)(\delta_{1,h},\ldots,\delta_{n,h})(f)(x).
\end{equation}

From now on, we fix a family $P_1,\ldots,P_n$ of $n$ elements in $\cD$.
Recall that  homogeneous polynomials 
define a complete intersection if and only if the only vector
on which they all vanish is the origin.
We  make the following assumption on $P_1,\ldots,P_n$:
\begin{equation}
\label{e:CI}
\tag{CI}
\begin{array}{l}
\mbox{For all }x \in \Rnp \mbox{, the elements }
\ev_{x,0}(P_k) \in \CC[D_1,\ldots,D_n] \mbox{ are homogeneous of the same degree }\\
\mbox{and define a complete intersection in }\CC^n.\mbox{ We let }d_1, \dots, d_n
\mbox{ be the degree of }\\ev_{x,0}(P_1), \dots, \ev_{x,0}(P_n)
\mbox{ and we let }\pi=\prod_i d_i \mbox{ and }\rho=\sum_i (d_i-1).
\end{array}
\end{equation}
Our goal is to study the space $\ST{\bm \eps}$ of functions $\Qe \to \CC$ which are {\em annihilated by the $P_k$'s}, i.e.
such that for every $h \in r\,\Qrp$ and $k \le n$, we have
$P_k(\delta_{1,h},\ldots,\delta_{n,h})(f)\equiv 0$, as function on $\Qe$.
\medskip

As an intermediate step, we introduce a finite variant.
For $x=(x_i) \in \Qe$, $h \in r\,\Qrp$ and $N \in \N$, we denote by
$S(x,h,N)$ the set of points $x+hy$ with $y=(y_i) \in \N^n$ such that
$\sum y_i \leq N$. We always assume that for all $k$, we have $N \geq d_k$.
The set of functions on $S(x,h,N)$ annihilated by the $P_k$'s will be denoted $\ST{S(x,h,N)}$.
More formally, let $\ST{S(x,h,N)}$ be the set of functions
$f:S(x,h,N) \to \CC$ such that, for every $k \le n$, we have
$
P_k(\delta_{1,h},\ldots,\delta_{n,h})(f)
\equiv 0,
$
as functions on  $S(x,h,N-d_k)$.\medskip

Finally, let $\cI \subset \cD$ be the left $\cD$-module linearly
generated by the elements $P_k$ for $k \in \{1,\ldots,n\}$.
We have the following simple characterization of elements of $\ST{S(x,h,N)}$.
\begin{lemma}
\label{equa-point}
Let $x \in \Qe ,h \in r\,\Qrp,N\in \N$ be fixed and let $f:S(x,h,N) \to \CC$ be fixed.
The following are equivalent:
\begin{enumerate}
\item $f \in \ST{S(x,h,N)}$;
\item for all $Q \in \cI$ of degree $d$
with $d \leq N$, we have
$Q(\delta_{1,h},\ldots,\delta_{n,h})(f)(x)=0 \in \CC$
\item for all $Q \in \ev_{x,h}(\cI)$ of degree $d$
with $d \leq N$, we have
$Q(\delta_{1,h},\ldots,\delta_{n,h})(f)(x)=0 \in \CC$..
\end{enumerate}
\end{lemma}
\begin{proof}
The equivalence of (2) and (3) follows from \eqref{equa-ev}, so let us only prove (1) $\Leftrightarrow$ (2).

Assume (1), i.e.~that $f$ is in $\ST{S(x,h,N)}$. Then for $k \leq n$,
$g \coloneqq P_k(\delta_{1,h},\ldots,\delta_{n,h})(f)$ is the zero function in
$\cF(S(x,h,N-d_k),\CC)$.
Applying \vf{further} divided difference operators $\delta_{i,h}$ to $g$ or multiplying
$g$ by some $x_i$ will still
yield the zero function, so that condition (2) is satisfied.

Conversely, assume (2). Then, given $k \leq n$,
$g \coloneqq P_k(\delta_{1,h},\ldots,\delta_{n,h})(f)$ is an element of
$\cF(S(x,h,N-d_k),\CC)$ satisfying
$D(g)(x)=0$ for all products $D$ of $d$ elements $\delta_{i,h}$
with $d\leq N-d_k$. This implies that $g \equiv 0$, proving (1).
\end{proof}

For our next Lemma, we recall Assumption \eqref{e:CI}.

\begin{lemma}
\label{lemm:kappa}
Let $N \in \N$ with $N>\pi$,
let $K \subset \Rnp$ compact. There exists $\kappa >0$ such
that for all $h$ in $\R$ with $|h| \leq \kappa$, and for all
$x \in K$, the quotient 
$ \CC_{\leq N}[D_1,\ldots,D_n]\big / \ev_{x,h} ( \cI ) $
has dimension at most $\pi$ and maximal degree at most $\rho$.
\end{lemma}
\begin{proof}
Let us first consider the case $h=0$ (and $x \in \Rnp$ arbitrary).
We let $\mathcal J$ be  the ideal 
generated by  $\ev_{x,0}(P_k)$ in $\CC[D_1,\ldots,D_n]$.
By Assumption \eqref{e:CI} and Lemma \ref{l:dim-coinv}, the quotient
 $\CC[D_1,\ldots,D_n] / \mathcal J$ 
 has dimension $\pi$ and maximal degree $\rho$.

However, since $\ev_{x,0}$ is not an algebra morphism, for $Q$ in $\CC[D_1,\ldots,D_n]$,
in general we have 
\[Q(D_1,\ldots,D_n) \ev_{x,0}(P_k) \ne \ev_{x,0}(Q(D_1,\ldots,D_n) P_k).\] 
Consequently $\mathcal J \ne \ev_{x,0}(\mathcal I)$.
But, for $h=0$, the second relation in~\eqref{commutationXD} simplifies to $X_iD_i=D_iX_i +r$,
 i.e.~$X_i$ and $D_i$ commute up to smaller degree terms.
Thus, if $Q$ is a homogeneous
polynomial in $n$ variables of degree $d$, then
$\ev_{x,0}(Q(D_1,\ldots,D_n) P_k)$ is equal to
the sum of $Q(D_1,\ldots,D_n) \ev_{x,0}(P_k)$ and terms of degree
strictly less that $d+d_k$. 
Thus it follows from Lemma
\ref{lemm:smaller-deg-perturbation} that
$\CC[D_1,\ldots,D_n] \big / \ev_{x,0} ( \cI )$
has dimension at most $\pi$ and maximal degree at most $\rho$.

Let us call $G \subset \CC^n \times \CC$ the set of pairs $(x,h)$ such that $\CC[D_1,\ldots,D_n] \big / \ev_{x,h} ( \cI )$ has dimension  at most $\pi$ and degree at most $\rho$.
The above argument shows that $G$ contains $\Rnp \times \{0\}$, and hence $K\times \{0\}$.
We claim that $G$ is open. Indeed, let $(x_0,h_0)$ be in $G$,
and let $(R_j)_{j \le \pi_0}$ be representatives of a basis of the quotient $\CC_{\leq N}[D_1,\ldots,D_n] \big / \ev_{x_0,h_0} ( \cI )$, with  $\pi_0 \le \pi$ elements.
 Since $(x_0,h_0)$ be in $G$, one can assume that each $R_j$
has degree at most $\rho$.
Also, we can find a finite family of polynomials $(Q_{j,k})_{j\le M}$
such that $\ev_{x_0,h_0}(Q_{j,k}(D_1,\ldots,D_n) P_k)$ form a basis of 
$\ev_{x_0,h_0} ( \cI ) \cap  \CC_{\leq N}[D_1,\ldots,D_n]$.
Hence $(R_j)_{j \le \pi_0}$ and \vupe{$\ev_{x_0,h_0}(Q_{j,k}(D_1,\ldots,D_n) P_k)$} together from a basis
of $ \CC_{\leq N}[D_1,\ldots,D_n]$.
But, for general $(x,h)$ in $\CC^n \times \CC$  the condition that
$R_j$ and $\ev_{x,h}(Q_{j,k}(D_1,\ldots,D_n) P_k)$
 form a basis of $\CC_{\leq N}[D_1,\ldots,D_n]$ can be expressed by
the non-vanishing of the determinant of a matrix whose coefficients are
polynomials in $x,h$. 
This determinant does not vanish at $(x_0,h_0)$,
and when it does not vanish at some $(x,h)$, we have that $(x,h)$ belongs to $G$
(the converse might not be true). Hence $G$
contains an open neighbourhood of $(x_0,h_0)$, and
it follows that $G$ is open.

Now we have proved that $G$ is an open subset of $\CC^n \times \CC$ which contains $K\times \{0\}$, hence it contains $K\times [-\kappa,\kappa]$ for $\kappa$ small enough.
This proves the lemma.
\end{proof}


\begin{lemma}
\label{res-inj}
Let $x \in \Qe,h \in r\,\Qrp,N\in \N$ be fixed with $N\ge \rho$. Assume $K \subset \Rnp$ is a compact
containing $x$ and $h \leq \kappa$, where $\kappa$ is provided by Lemma~\ref{lemm:kappa}.
Then, restriction yields
an injective map $\ST{S(x,h,N)} \to \ST{S(x,h,\rho)}$ and
$\dim \ST{S(x,h,N)} \leq \pi$.
\end{lemma}
\begin{proof}
Let $f \in \ST{S(x,h,N)}$. Considering a product $D$ of $D_i$'s  of degree $d \le N$,
we can write
\[D = Q + \sum_{j =1}^{\pi_0} \alpha_j R_j,\]
 where $Q$ is in $\ev_{x,h}(\cI)$, $\alpha_j \in \CC$, and the elements
 $R_j$ for $1\le j\le \pi_0$ are representatives of a basis of the quotient
 $\CC[D_1,\ldots,D_n] / \ev_{x,h}(\cI)$. Moreover, by Lemma~\ref{lemm:kappa},  one has $\pi_0 \le \pi$ and one can assume that the $R_j$'s have degree at most $\rho$.
By Lemma~\ref{equa-point}, one has $Q (\delta_{1,h},\ldots,\delta_{n,h}) f=0$, so that
 \begin{equation}\label{eq:Dfx}
D(\delta_{1,h},\ldots,\delta_{n,h}) f(x) = \sum_{j =1}^{\pi_0} \alpha_j R_j(\delta_{1,h},\ldots,\delta_{n,h}) f (x).
 \end{equation}
From here, we easily prove the two statements in the lemma:
\begin{itemize}
\item Assume that $f_{|S(x,h,\rho)}=0$. Then it follows that $P(\delta_{1,h},\ldots,\delta_{n,h})f(x)=0$ for all polynomial $P$ in $\CC[D_1,\ldots,D_n]$ of degree at most $\rho$. In particular $R_j(\delta_{1,h},\ldots,\delta_{n,h})f(x)=0$ for all $j$,
and using \eqref{eq:Dfx}, we see that
\[D(\delta_{1,h},\ldots,\delta_{n,h}) f(x) =0.\]
Since $D$ can be any product of $D_i$'s  of degree $d \le N$, this implies that
$f \equiv 0$ on $S(x,h,N)$. Hence restriction is an injective map from
$\ST{S(x,h,N)}$ to $\ST{S(x,h,\rho)}$, as claimed.
\item On the other hand, a function on $S(x,h,N)$ is determined by the values
$D(\delta_{1,h},\ldots,\delta_{n,h}) f(x)$, where $D$ runs over all possible products of 
$D_i$'s  of degree $d \le N$.
 By \eqref{eq:Dfx}, if $f$ is in  $\ST{S(x,h,N)}$, such values are themselves determined
 by the $\pi_0$ scalar values $R_j(\delta_{1,h},\ldots,\delta_{n,h})(x)$.
 This proves $\dim(\ST{S(x,h,N)}) \le \pi_0 \le \pi$.\qedhere
\end{itemize}
\end{proof}

We now bound the dimension of the space $\ST{\bm \eps}$ of functions on  the whole space $\Qe$ that are annihilated by the $P_k$.

\begin{figure}
\centering
\includegraphics{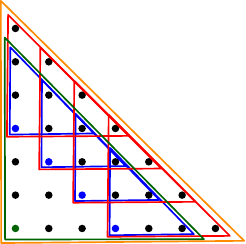}
\caption{Illustration of the proof of injectivity of the restriction from $\ST{S(x_0,h,M)}$ to $\ST{S(x_0,h,M-1)}$ in the case $r=1$, $n=2$, $\rho=2$ and $M=6$. The point $x_0$ is in green, while points $y$ in $S(x_0,h,M-\rho-1)$
and not in $S(x_0,h,M-\rho-2)$ are in blue.
Lemma~\ref{res-inj} implies that each restriction map from a red triangle to the corresponding blue triangle is injective. We deduce immediately that the restriction map from the orange triangle to the green one is injective as well.}
\label{fig:induction-on-triangles}
\end{figure}

\begin{proposition}\label{p:bound_dim_general}
$\dim(\ST{\bm \eps}) \le \pi$.
\end{proposition}
\begin{proof}
Assume for the sake of contradiction $\dim\, \Se > \pi$.
Then we can find a {\em finite} subset
$Z \subset \Qe$ such that
\[\dim \{ f_{|Z} :\, f\in \Se \}>\pi.\]
Given $Z$ and $\bm \eps$, we introduce a number of objects.
\begin{itemize}
\item
First, let us choose $d \in \N$ such that $[d]=1$ and $Z \subset \frac 1d \Ze$.
We denote by $x_0$ the smallest element in $\frac 1d\Ze$, i.e.~$x_0  \coloneqq  \frac 1d(\tilde \epsilon_1,\ldots,\tilde \epsilon_n)$, where $\tilde \eps_i \in \iset{r}$ is such that \vupe{$[\tilde \eps_i]=\eps_i$}.
By construction, we have $Z \subset x_0 + \frac rd \mathbb Z^n_{\ge 0}$.
\item
Second, we can find $A>0$, such that $Z$ and $x_0$ are both contained in the compact set
\[  K_{\tau,A} \coloneqq \Big\{(x_1,\dots,x_n) \in \Rnp: \sum_{i \le n} x_i \le A \text{ and, for all }i \le n, \ x_i \ge \tfrac1d \Big\} .\]
\item Third, choosing $N=\rho+1$ and $K$ in Lemma~\ref{lemm:kappa}
provides us some $\kappa>0$. We can then take a multiple $D$ of $d$ with $[D]=1$ such that
$h \coloneqq \frac rD <\kappa$. Since $D$ is a multiple of $d$, we have $Z \subset x_0 + h \mathbb Z^n_{\ge 0}$.
\end{itemize}

We now start bounding dimensions of space of functions annihilated by the $P_k$'s recursively.
Since $x_0$ is in $K_{\tau,A}$ and $h < \kappa$, Lemma~\ref{res-inj}
ensures that 
$\ST{S(x_0,h,\rho\vf{+1})}$ has dimension at most $\pi$. 

Now, we claim that for all $M > \rho$,
the restriction map from $\ST{S(x_0,h,M)}$ to $\ST{S(x_0,h,M-1)}$ is injective,
as long as $S(x_0,h,M-\rho-1) \subset K$.
Note that $M$ is typically bigger than $N=\rho+1$, so that we cannot use directly
Lemma~\ref{res-inj} (our value of $\kappa$ was fixed so that Lemmas~\ref{lemm:kappa} and~\ref{res-inj} hold for $N=\rho+1$, not for larger $N$).

However, we can apply Lemma~\ref{res-inj} to any $y \in S(x_0,h,M-\rho-1)$ and $N=\rho+1$. This implies that the restriction map from $\ST{S(y,h,\rho+1)}$ to 
$\ST{S(y,h,\rho)}$ is injective. But the set $S(x_0,h,M)$ is covered by
\[S(x_0,h,M-\rho-2) \cup \bigcup_{y \in S(x_0,h,M-\rho-1)} S(y,h,\rho+1).\]
Hence the values of a function $f$ in $\ST{S(x_0,h,M)}$ are determined by
its values on 
\[S(x_0,h,M-\rho-2) \cup \bigcup_{y \in S(x_0,h,M-\rho-1)} S(y,h,\rho),\]
which is equal to $S(x_0,h,M-1)$.
We have proved that the restriction from $\ST{S(x_0,h,M)}$ to $\ST{S(x_0,h,M-1)}$ is indeed injective (the argument is illustrated on Fig.~\ref{fig:induction-on-triangles}).
An immediate induction implies that $\ST{S(x_0,h,M)}$ has dimension at most $\pi$,
as long as we have the inclusion $S(x_0,h,M-\rho-1) \subset K$.

Recalling that $Z \subset K \cap (x_0+h\Znp)$.
It implies that,
for $M$ sufficiently large, we have $Z \subset S(x_0,h,M)$.
Moreover since $Z \subset K$ and $K$ is of the form $K_{\tau,A}$ for some $\tau,A>0$,
this $M$ can be chosen so that $S(x_0,h,M-\rho-1) \subset K$.
 Then, since every element $f$ in $\Se$ gives an element of $\ST{ S(x_0,h,M)}$
 by restriction, we have
 \[ \{ f_{|Z} :\, f\in \Se \} \subseteq \{ f_{|Z} : f \in \ST{ S(x_0,h,M)} \}.\]
 But the former space has dimension bigger than $\pi$, while we have proved that the latter
 has dimension at most $\pi$. We have reached a contradiction, proving $\dim\, \Se \le \pi$,
 as wanted.
\end{proof}


\subsection{Conclusion}
\label{s:conclusion}
We can now prove the main results of the section.
First Proposition~\ref{p:bound} is simply a particular case of Proposition~\ref{p:bound_dim_general}
\begin{proof}[Proof of Proposition~\ref{p:bound}]
We set $P_k=\sum_i X_iD_i^k \in \cD$, thus $\pi=n!$\,.
For $x \in \Rnp$, we have $\ev_{x,0}(P_k)=\sum x_i D_i^k$, 
and Assumption~\eqref{e:CI} is satisfied, by Lemma \ref{l:sol-triviale}.
Now, by Lemma \ref{lem:discrete_diff_cf}, any $r$-homogeneous and $r$-polyharmonic
function $c$ on $\Ze$ can be extended to a function $f$ satisfying \eqref{e:f}, i.e.~living in $\Se$. This extension is clearly an injective map, so that the dimension of the space of $r$-homogeneous and $r$-polyharmonic
function is at most that of $\Se$, which does not exceed $n!$ by Proposition~\ref{p:bound_dim_general}.
\end{proof}
Combining this upper bound on dimension with the explicit independent family
of solutions exhibited in Section~\ref{s:independence_semi_products}, we get the following result.

\begin{theorem}
\label{t:linear}
\vupe{The functions $c_\sigma$,
where $\sigma$ runs \pe{over} $\fS_n$, form a basis of the space of homogeneous polyharmonic functions $c:\Znp \to \CC$.
Consequently,}
a basis of $\LI_n(\CC[T]_0)$ is provided by the elements $\varphi_\sigma$, for $\sigma \in \fS_n$.
\end{theorem}
\begin{proof}
\vupe{
By \eqref{value-bracket} and Proposition~\ref{p:induction}, the functions $c_\sigma$ are
homogeneous and polyharmonic. Moreover, they are linearly independent by Lemma~\ref{l:c-sigma_ind}.
On the other hand Proposition~\ref{p:bound} (for $r=1$ and $\bm \eps=(0,\dots,0)$) states
that the space of homogeneous and polyharmonic functions $c:\Znp \to \CC$ has dimension at most $n!$.
\pe{We conclude that $(c_\sigma)_{\sigma \in \fS_n}$ is a basis of this space.}}

\vupe{
Regarding $LI_n(\CC[T]_0)$, Lemma~\ref{l:form_multilinear_solution} and Proposition~\ref{p:induction}
shows that any element $\varphi$ of $LI_n(\CC[T]_0)$ writes as 
\[\varphi(T^{x_1},\ldots,T^{x_n})=c(x_1,\ldots,x_n) T^{\sum_i x_i},\]
for some homogeneous and polyharmonic functions $c:\Znp \to \CC$. By the first part of the theorem,
we see that $\varphi$ must be a linear combination of $(\varphi_\sigma)_{\sigma \in \fS_n}$.
Since the $\varphi_\sigma$ are indeed in $LI_n(\CC[T]_0)$ (Proposition~\ref{prop:f1f2p})
 and are linearly independent (Lemma~\ref{l:phi-sigma_ind}),
they form a basis of $LI_n(\CC[T]_0)$, as asserted.}
\end{proof}

\subsection{A continuous analogue of polyharmonicity}
\label{ss:system_EDP}
In this section, we consider a system of PDEs which is a continous analogue of the 
 polyharmonicity relations \eqref{recurrence-c} ; namely we look for distributions 
 $c \in \mathscr{D}'((0,+\infty)^n)$ satisfying
 \begin{equation}
\label{equa-diff-c}
\forall k \in \N,\ \sum_{i=1}^n x_i \frac{\partial^k c}{\partial x_i^k} = 0
\end{equation} 
 It turns out the results of this section provide us with an explicit basis of solutions for this system of PDEs.
 \begin{theorem}
 The functions
$c_\sigma:\;(0,+\infty)^n \to \CC$ defined by
$$
c_\sigma(x_1,\ldots,x_n)=\frac{\prod_{i=1}^n x_i}
{(x_{\sigma(1)}+\cdots+x_{\sigma(n)})\cdots(x_{\sigma(n-1)}+x_{\sigma(n)})x_{\sigma(n)}}
$$
indexed by permutations $\sigma$ of $n$
form a basis of the space of distribution solutions of the system \eqref{equa-diff-c} of PDEs.
 \end{theorem}
\begin{proof}
Note first that distribution solutions  of the system \eqref{equa-diff-c} of PDEs are in fact analytic functions.  The differential operator
\[  \sum_{i=1}^n x_i \frac{\partial^2 }{\partial x_i^2} = 0 \]
has real analytic coefficients in $(0,+\infty)^n$ and its characteristic set
\[ \bigg\{ (x,\xi) \in \pe{(0,+\infty) \times (\R_*)^n} :  \sum_{i=1}^n x_i \xi_i^2 = 0 \bigg\} \]
is empty.
 Therefore, by~\cite[Theorem 8.6.1]{Hor}, the analytic wave front set $\mathop{\rm WF}_A(c)$ of \vupe{any solution $c$ of \eqref{equa-diff-c}}
 is empty.
Since the projection of the analytic wave front set on $(0,\infty)^n$ is the analytic singular support $\mathop{\rm sing supp}_A(c)$ of $c$ (see~\cite[Theorem 8.4.5]{Hor}), the distribution $c$ coincides with an analytic function $c : (0,\infty)^n \to \CC$.
In the rest of the proof we assume therefore that $c$ is real-analytic.

We next prove that for each permutation $\sigma$ of $n$,
the function $c_\sigma$ indeed solves \eqref{equa-diff-c}.
Recall by Proposition~\ref{prop:f1f2p} that $\varphi_\sigma \in \LI_n(E)$.
By \eqref{value-bracket} and Proposition~\ref{p:induction}, it follows that $c_\sigma$ satisfies 
\[ \big ( \sum_{i=1}^n x_i \big ) \, c_\sigma(x_1,\ldots,x_n)
= \sum_{i=1}^n \, x_i \, c_\sigma(x_1,\ldots,x_i + u,\ldots,x_n),\]
 whenever
$x_1,\ldots,x_n,u$ belong to $\Zp$. 
Since this is an equality of rational functions, it holds for all
$x_1,\ldots,x_n,u \in \R$. 
Performing a Taylor expansion in the variable $u$,
and extracting the coefficient of $u^k$, we get
\[0= \frac{1}{k!}\,\sum_{i=1}^n x_i \frac{\partial^k }{\partial x_i^k}c_\sigma.\]
Hence $c_\sigma$ satisfies \eqref{equa-diff-c}, as claimed.

Next, we know by Lemma~\ref{l:c-sigma_ind} that the $c_\sigma$ (indexed
by permutations $\sigma$ of $n$) are linearly independent as rational functions,
hence also as functions on $(0,+\infty)^n$.
It remains to show that they span the whole space of analytic solutions or,
 equivalently, that the space of analytic solutions has dimension at most $n!$. 
 
 Recall that two analytic functions on $(0,+\infty)^n$ coincide if and only if they have the same partial derivatives of all orders at $x_0=(1,\dots,1)$.
 Therefore we focus at partial derivatives at this point.
 To this end, let us use the notation of Section~\ref{sec:upper_bound_general}
 with $h=0$, $\delta_i=\delta_{i,0}$ being the differential with repect to $x_i$,
 and $P_k=\sum x_i \delta_i^k$.
 By a straightforward analogue of Lemma~\ref{equa-point},
 a function $c$ satisfies \eqref{equa-diff-c} at $x_0$ if and only if
 \[Q(\delta_1,\dots,\delta_n)(c)\, (x_0)=0,\]
 for all $Q$ in $\ev_{x_0,0}(\mathcal I)$, 
 where $\mathcal I$ is the left ideal generated by the $P_k$'s.
 But it follows from the proof of Lemma~\ref{lemm:kappa},
  that $\CC[D_1,\dots,D_n]/\ev_{x_0,0}(\mathcal I)$ has dimension at most $n!$.
  Hence it is spaned by $n!$ differential operators $(R_j)_{j \le n!}$.
  Now, considering any product $D$ of partial derivatives $\delta_1,\dots,\delta_n$,
  we can decompose $D$ as $D=Q+\sum_{j \le n!} \alpha_j R_j$, where $Q$ is in $\ev_{x_0,0}(\mathcal I)$ and $\alpha_j$ in $\CC$. It follows that
  \[ D(\delta_1,\dots,\delta_n)(c)\, (x_0) = \sum_{j \le n!} \alpha_j R_j(\delta_1,\dots,\delta_n)(c)\, (x_0).\]
  Hence the values of all partial derivatives of $c$ in $x_0$ depend linearly 
  on the $n!$ parameters $(\alpha_j)_{j \le n!}$.
  Since analytic functions are completely determined by such derivatives,
  this proves that the space of  solutions of \eqref{equa-diff-c} has dimension
  at most $n!$. The theorem is proved.
\end{proof}

\begin{remark}
Again, an easy commutator computation shows that writing $P_k=\sum x_i \delta_i^k$ as in the above proof, we have $[P_k,P_j]=(k-j)P_{k+j-1}$. Hence satisfying \eqref{equa-diff-c} for $\bm{k \le 3}$ is equivalent to satisfying it for all $k$.
\end{remark}

\section{Algebraic solutions}
\label{sec:algebraic}
We now study $\AI_n(E)$ for $E=\CC[T]_0$, 
i.e.~the solutions $\varphi$ of \eqref{e:reparametrage} such that
the coefficients of $\varphi(f_1,\dots,f_n)$ depend polynomially 
on the coefficients of $f_1$,\dots,$f_n$. 
It turns out that this problem can be reduced to the problem of finding
multilinear solutions, thanks to a standard trick called {\em polarization}.

\subsection{Polarization of multihomogeneous polynomials}
\label{s:polarization}

When $E$ is a vector space, it is well-known that
a polynomial homogeneous function $\varphi : E \to \CC$ of degree $d$
admits a unique polarization, namely a symmetric $d$-linear function
$\wvar : E^d \to \CC$ such that for all $x \in E$, we have
$\wvar(x,\ldots,x)=\varphi(x)$. We recall this result in the more general context
of multilinear functions.

\vf{Let $E_1,\ldots,E_n,F$ be $\CC$-vector spaces
(possibly of infinite dimension).
A polynomial function from $E_1 \times \dots \times E_n$ to $F$
is an element of $\big(\bigotimes_i S E_i^* \big) \otimes F$,
where $S E^*$ is the symmetric algebra on the dual space $E^*$ of $E$.
Concretely, if we fix a basis of $F$
and a basis $(\alpha_{i,b})_{b \in B_i}$ of each of the dual spaces $E_i^*$
(for $i$ in $\{1,\dots,n\}$), then 
a function $\varphi:E_1 \times \dots \times E_n \to F$
is called {\em polynomial} if each coordinate of
$\varphi(-,\ldots,-)$ is a polynomial in $(\alpha_{i,b}(-))_{i\le n, b\in B_i}$ (thus involving only a finite number of those); this notion does not depend on the chosen bases.}

Let $\tau_1,\ldots,\tau_n \in \CC$ and $\varphi:E_1 \times \dots \times E_n \to F$ be a polynomial
function. We may write
$$
\varphi(\tau_1 x_1,\ldots,\tau_n x_n) = \sum_{d_1,\ldots,d_n}
\tau_1^{d_1} \cdots \tau_n^{d_n} \varphi_{(d_1,\ldots,d_n)}
(x_1,\ldots,x_n).
$$
The element $\varphi_{(d_1,\ldots,d_n)}$ is called the multihomogeneous
part of $\varphi$ of degree $(d_1,\ldots,d_n)$.

Now, given a multihomogeneous polynomial of degree $(d_1,\ldots,d_n)$
with $d=\sum d_i$,
there is a unique $d$-multilinear and
$(\fS_{d_1} \times \fS_{d_2} \times \dots \times \fS_{d_n})$-invariant
function
$\wvar:E_1^{d_1} \times \cdots \times E_n^{d_n} \to F$
such that for all $x_1,\ldots,x_n \in E_1 \times \cdots \times E_n$, we
have
\begin{equation}
\label{e:def-polar}
\wvar(x_1,\ldots,x_1,x_2,\ldots,x_2,\ldots,x_n,\ldots,x_n)=
\varphi(x_1,\ldots,x_n)
\end{equation}
(each $x_i$ is repeated $d_i$ times).
Indeed, existence can be proved by first considering a monomial $\varphi$ --
in which case $\wvar$ is simply the product of variables --
and then extending to all polynomials by linearity.
Uniqueness follows from the fact that \eqref{e:def-polar} implies that
$d_1!\cdots d_n!\; \wvar(y_1,y_2,\ldots,y_d)$ is the coefficient of
$\alpha_1 \cdots \alpha_d$ in
\begin{equation}
\label{e:unic-polar}
\varphi(\alpha_1 y_1+\cdots+\alpha_{d_1} y_{d_1},\ldots,
\alpha_{d-d_n+1} y_{d-d_n+1} + \alpha_d y_d).
\end{equation}
This function $\wvar$ is called the polarization of $\varphi$.
\vupe{Both as an illustration and as preparation for the next section,
let us compute the polarization of $\varphi_\sigma$.}

\begin{lemma}\label{lem:pol_phi}
\vupe{Let $\sigma$ be a function from $\{1,\ldots,m\}$ to $\{1,\ldots,n\}$, for $n,m \ge 1$.
For $i \le n$, let $d_i=|\sigma^{-1}(i)|$ and $\tau:\{1,\ldots,m\} \to \{1,\ldots,n\}$
be the unique nondecreasing function taking value $i$ exactly $d_i$ times (for each $i \le n$). Then
 \begin{equation}\label{eq:pol_phi}
  \widetilde{\varphi_\sigma} = \frac{1}{d_1!\, \cdots\, d_n!} \sum_{\pi \in \fS_m:\, \sigma=\tau \circ \pi} \varphi_\pi.
  \end{equation}}
\end{lemma}
\begin{proof}
\vupe{Let us call $\psi_\sigma$ the right-hand side of~\eqref{eq:pol_phi}.
For $\pi_0$ in $\fS_{d_1} \times \fS_{d_2} \times \dots \times \fS_{d_n}$, we have
$\tau \circ \pi_0=\tau$, and hence, left-composition of $\pi$ by $\pi_0$ permutes the elements in the sum in~\eqref{eq:pol_phi}.
Since 
\[\varphi_{\pi_0 \circ \pi}(g_1, \dots,g_m)=
[g_{\pi_0(\pi(1))}:\dots:g_{\pi_0(\pi(m))}]=\varphi_\pi(g_{\pi_0(1)}, \dots,g_{\pi_0(m)}),\]
we get
\begin{align*}
\psi_\sigma(g_1, \dots,g_m)&=\frac{1}{d_1!\, \cdots\, d_n!} \sum_{\pi \in \fS_m:\, \sigma=\tau \circ \pi} \varphi_{\pi_0 \circ \pi}\\&=\frac{1}{d_1!\, \cdots\, d_n!} \sum_{\pi \in \fS_m:\, \sigma=\tau \circ \pi}
\varphi_\pi(g_{\pi_0(1)}, \dots,g_{\pi_0(m)})=\psi_\sigma(g_{\pi_0(1)}, \dots,g_{\pi_0(m)}).
\end{align*}
This proves that $\psi_\sigma$ is $\fS_{d_1} \times \fS_{d_2} \times \dots \times \fS_{d_n}$-invariant. It remains to prove that specializing to repeated arguments, i.e.~to
$(f_{\tau(1)},\dots,f_{\tau(m)})$, gives $\varphi_\sigma(f_1,\dots,f_n)$.
But this is clear since, for $\pi \in \fS_m$ satisfying $\sigma=\tau \circ \pi$, we have
\[\varphi_\pi(f_{\tau(1)},\dots,f_{\tau(m)})= \varphi_{\tau \circ \pi}(f_1,\dots,f_n)
=\varphi_\sigma(f_1,\dots,f_n).\]
By uniqueness of polarization, we get $\widetilde{\varphi_\sigma} =\psi_\sigma$, as wanted.}
\end{proof}

\subsection{Generating set of $\AI_n(\CC[T]_0)$}
Recall from Definition~\ref{d:varphi-sigma} that for a function $\sigma$ from $\{1,\ldots,m\}$ to $\{1,\ldots,n\}$, we det
\[\varphi_\sigma(f_1,\ldots,f_n) = [f_{\sigma(1)}:f_{\sigma(2)}:\cdots:f_{\sigma(m)}].\]
\begin{theorem}
\label{t:algebraic}
A basis of $\AI_n,m(\CC[T]_0)$ is provided by the elements
$\varphi_\sigma$, for $\sigma:[1,\ldots,m] \to [1,\ldots,n]$.
\end{theorem}
\begin{proof}
First consider $\varphi \in \AI_n(E)$. 
Being algebraic, the map $\varphi$ can be considered as a polynomial
function $E^n \to E$, so that the results of Section \ref{s:polarization}
apply. We can assume that $\varphi$ is multihomogeneous of degree
$(d_1,\ldots,d_n)$, and consider its polarization $\wvar$.
Let $d=\sum d_i$, let $f_1,\ldots,f_d,\gamma \in E$ and let
$\alpha_1,\ldots,\alpha_d \in \CC$.

We have $(\alpha_1f_1+\cdots+\alpha_{d_1}f_{d_1}) \circ \gamma
=\alpha_1f_1 \circ \gamma+\cdots+\alpha_{d_1}f_{d_1}\circ \gamma$, and
similarly for other linear combinations. Condition
\eqref{e:reparametrage} implies that
\begin{multline*}
\varphi(\alpha_1 f_1 \circ \gamma+\cdots+\alpha_{d_1}
f_{d_1} \circ \gamma,
\ldots, \alpha_{d-d_n+1} f_{d-d_n+1} \circ \gamma  + \dots+ \alpha_d f_d \circ \gamma) \\
=
\varphi(\alpha_1 f_1 +\cdots+\alpha_{d_1} f_{d_1},
\ldots, \alpha_{d-d_n+1} f_{d-d_n+1} + \cdots + \alpha_d f_d) \circ \gamma
\end{multline*}
The coefficient of $\alpha_1 \cdots \alpha_d$ is
$d_1! \cdots d_n!$ times
$\wvar(f_1 \circ \gamma,\ldots,f_d \circ \gamma)$ in the left hand side
and $\wvar(f_1,\ldots,f_d) \circ \gamma$ in the right hand side, so that
$\wvar$ commutes with reparametrizations and is an element of $\LI_d(E)$.

By Theorem~\ref{t:linear}, we may write $\wvar = \sum_{\pi \in \fS_d} a_\pi \varphi_\pi$ for some scalars $a_\pi \in \CC$.
Let $\tau:\{1,\ldots,m\} \to \{1,\ldots,n\}$ be the \vupe{unique} nondecreasing function \vupe{taking} each value $i \in \{1,\ldots,n\}$ \vupe{exactly} $d_i$
times. We have
\begin{align*}
\varphi(f_1,\ldots,f_n) & =  \wvar(f_{\tau(1)},\ldots,f_{\tau(m)}) 
= \sum_{\pi \in \fS_m} a_\pi \varphi_\pi(f_{\tau(1)},\ldots,f_{\tau(m)}) \\
&=\sum_{\pi \in \fS_m} a_\pi [f_{\tau \circ \pi(1)}: \cdots : f_{\tau \circ \pi(m)}] 
= \sum_{\pi \in \fS_m} a_\pi \varphi_{\tau \circ \pi}(f_1,\ldots,f_n).
\end{align*}
This shows that $(\varphi_\sigma)_{m \ge 0, \sigma:\{1,\ldots,m\} \to \{1,\ldots,n\}}$
is a generating family of $\AI_n(E)$.

It remains to show that it is free. Assume that we have a linear combination $\sum_{\sigma} a_\sigma \varphi_\sigma$ that vanishes. Up to extracting the multihomogenous components,
we can assume that it is multi-homogeneous of degree $(d_1,\dots,d_n)$, i.e.~that
all functions $\sigma$ involved satisfies $|\sigma^{-1}(i)|=d_i$ for all $i \le n$.
Looking at the relation $\sum_{\sigma} a_\sigma \varphi_\sigma=0$, 
and using the uniqueness of polarization, we get $\sum_{\sigma} a_\sigma \widetilde{\varphi_\sigma}=0$.
\vupe{Writing $m=d_1+\dots+d_n$
and using Lemma~\ref{lem:pol_phi},
this simplifies to}
\[ \sum_{\sigma: \{1,\ldots,m\} \to \{1,\ldots,n\} \atop |\sigma^{-1}(i)|=d_i} a_\sigma \left( \sum_{\pi \in \fS_m:\, \sigma=\tau \circ \pi} \varphi_\pi \right)
=\sum_{\pi \in \fS_m} a_{\tau \circ \pi} \varphi_\pi =0.\]
Now recall from Lemma~\ref{l:phi-sigma_ind} that the multilinear functions
$\varphi_\pi$ indexed by {\em permutations} $\pi$ are independent,
so that the above relation implies $a_{\tau \circ \pi}=0$ for all $\pi$ in $S_m$.
Clearly, each $\sigma:\{1,\ldots,m\} \to \{1,\ldots,n\}$ 
satisfying $|\sigma^{-1}(i)|=d_i$ for all $i \le n$ decomposes as $\tau \circ \pi$ for some $\pi$, so that $a_\sigma=0$ for all such $\sigma$.
Since this holds for all $d_1,\dots,d_n\ge 0$, the family $(\varphi_\sigma)$ is a free family,
as claimed.
\end{proof}

\section{The truncated case $E=\CC[T]_0/(T^\ell)$}
\label{sec:truncated}

Let us now consider the case $E=\CC[T]_0/(T^\ell)$.
Similarly to the case $E=\CC[T]$, all multilinear solutions of \eqref{e:reparametrage} are of the form
\[\varphi(T^{x_1},\dots,T^{x_n})=c(x_1,\dots,x_n) T^{x_1+\dots+x_n},\]
but now the complex-valued function $c$ is only defined on
\[\Omega^n_\ell \coloneqq \{(x_1,\dots,x_n) \in \Zp^n: x_1+\dots+x_n < \ell\}.\]
The function $c$ must also necessarily satisfy the same induction relation as in Proposition~\ref{p:induction},
namely, for all $k \ge 1$ and $x_1+\cdots+x_n \le \ell-k$, we have
\begin{equation}\label{eq:recurrence_truncated}
\big ( \sum_i x_i \big )\, c(x_1,\ldots,x_n)
= \sum_i \, x_i \, c(x_1,\ldots,x_i + k,\ldots,x_n).
\end{equation}
\vupe{Throughout this section,
 let us denote by $\cS^n_\ell$ the space of solutions of \eqref{eq:recurrence_truncated} on $\Omega^n_\ell$.}
In this setting, since $c$ is not defined for large $x_1,\dots,x_n$,
we cannot extend solutions of~\eqref{eq:recurrence_truncated} to the whole space $\Qnp:w$ and apply the deformation techniques
from Section~\ref{sec:dimension_n!}.
\vupe{We believe nevertheless that $\cS^n_\ell$
has again dimension $n!$, provided that $\ell$ is large enough.
We state this as a conjecture}
\begin{conjecture}\label{conj:truncated_solutions}
  \vupe{For $\ell\ge \binom{n+1}{2}$, we have $\dim(\cS^n_\ell)=n!$.}
\end{conjecture}
\vupe{Note that this conjecture implies that the solutions of~\eqref{eq:recurrence_truncated}
on the whole space $\Znp$ has dimension at most $n!$, which is the most difficult part of Theorem~\ref{t:linear}.
In this section, we provide evidence for Conjecture~\ref{conj:truncated_solutions}. 
In particular, our first approach proves} the conjecture for $n=2$ and $n=3$,
and giving strong evidence for $n=4$.

\subsection{First approach: fixed size subsystems}
The idea here is to compute the coefficient $c(x_1,\cdots,x_n)$ recursively on the degree $x_1+\cdots+x_n$.
If we know the coefficients up to some degree $d$, then the coefficients $c(x_1,\cdots,x_n)$ of degree $d+1$
are solutions of a system of equations
\[ \begin{cases}
 \sum_i \, x_i \, c(x_1,\ldots,x_i + 1,\ldots,x_n)
 = \big ( \sum_i x_i \big ) \, c(x_1,\ldots,x_n), & \text{ for $x_1+\cdots+x_n=d$;}\\
\sum_i \, x_i \, c(x_1,\ldots,x_i + 2,\ldots,x_n)
= \big ( \sum_i x_i \big ) \, c(x_1,\ldots,x_n), & \text{ for $x_1+\cdots+x_n=d-1$.}
\end{cases}\]
Note that taking more equations corresponding to $k \ge 3$ is not useful because of Remark~\ref{rmk:commutators}.
From Conjecture~\ref{conj:truncated_solutions}, we expect that, for $d$ large, this system of equations
has at most one solution, so that the dimension does not increase when we increase $\ell$.
Note that the size of the system and the number of variables increase with $d$. 
However, it is possible to extract some subsystems of fixed size (with respect to $d$, but the size still depends on $n$),
and to show that these systems have unique solutions for small values of $n$.

\subsubsection{The $n=2$ case}
\begin{proposition}\label{prop:n2}
Conjecture~\ref{conj:truncated_solutions} holds for $n=2$ and all $\ell \ge 3$.
\end{proposition}
\begin{proof}
\vupe{
For $\ell=3$, we have three variables $c(1,1)$, $c(1,2)$ and $c(2,1)$,
and a single equation $2c(1,1)=c(1,2)+c(2,1)$, hence the space of solution has indeed dimension $2$.}

\vupe{Restricting solutions from $\Omega^2_{\ell+1}$ to $\Omega^2_{\ell}$
defines a map from $\cS^2_{\ell+1}$ to $\cS^2_{\ell}$.
We will show that this map is injective.
To this end, let $c$ be a function on $\Omega^2_{\ell+1}$ solving~\eqref{eq:recurrence_truncated},
such that $c(i,j)=0$ as soon as $i+j \le \ell$.}
Let $(i_0,j_0)$ be positive integers with $i_0+j_0=\ell-1$.
Equations \eqref{eq:recurrence_truncated} with $(x_1,x_2,k)$ specialized to $(i_0,j_0,2)$, $(i_0+1,j_0,1)$
and $(i_0,j_0+1,1)$ can be rewritten in matrix form as
$$
\begin{pmatrix}
 j_0 & 0 & i_0 \\
 j_0+1 & i_0 & 0 \\
 0 & j_0 & i_0+1
\end{pmatrix}
\begin{pmatrix}
 c(i_0,j_0+2) \\
 c(i_0+1,j_0+1) \\
 c(i_0+2,j_0)
\end{pmatrix}
=
\begin{pmatrix}
 (i_0+j_0) c(i_0,j_0) \\
 (i_0+j_0+1) c(i_0,j_0+1) \\
 (i_0+j_0+1) c(i_0+1,j_0)
\end{pmatrix}.
$$
The right-hand side is the zero vector by hypothesis.
Moreover, we can compute
$$\det \begin{pmatrix}
 j_0 & 0 & i_0 \\
 j_0+1 & i_0 & 0 \\
 0 & j_0 & i_0+1
\end{pmatrix} = -i_0j_0(i_0+j_0+2),$$
showing that this matrix is invertible for all pairs $(i_0,j_0)$ of positive integers.
Hence $c(i_0,j_0+2)=c(i_0+1,j_0+1)=c(i_0+2,j_0)$. Since this holds for all $i_0,j_0$ with $i_0+j_0=\ell-1$,
this proves that all coefficients of degree $\ell+1$ are zero.
Hence $c$ is identically zero, and \vupe{the restriction map is an injection from $\cS^2_{\ell+1}$ to $\cS^2_{\ell}$.
Thus $\dim(\cS^2_{\ell+1}) \le \dim(\cS^2_{\ell})$, and an immediate induction yield $\dim(\cS^2_{\ell}) \le 2$
for all $\ell \ge 3$.}

\vupe{It remains to prove that $\dim(\cS^2_{\ell}) \ge 2$. If, for some $\ell_0 \ge 3$,
we had $\dim(\cS^2_{\ell}) < 2$, then, by the above argument, we would have $\dim(\cS^2_{\ell}) < 2$ for all $\ell \ge \ell_0$.
Consequelty, the space of solutions of polyharmonic functions on $\mathbb Z_{>0}^2$ would have dimension less than 2.
This is in contradiction with the indendence of $c_\sigma$'s (Lemma~\ref{l:c-sigma_ind}) -- recall that the $c_\sigma$'s
are polyharmonic by Propositions~\ref{prop:f1f2p} and~\ref{p:induction}.
This concludes the proof of the proposition.}
\end{proof}

\subsubsection{The $n=3$ case}
\begin{proposition}\label{prop:n3}
Conjecture~\ref{conj:truncated_solutions} holds for $n=3$ and all $\ell \ge 6$.
\end{proposition}
\begin{proof}
We start by looking at the case $\ell=6$.
Then \eqref{eq:recurrence_truncated} defines a linear system whose variables are the $c(x_1,x_2,x_3)$, indexed by sequences $(x_1,x_2,x_3)$ of positive integers, whose sum
is at most $6$. There are 20 of them.
On the other hand, the equations are indexed by pairs $((x_1,x_2,x_3),k)$,
where  $(x_1,x_2,x_3)$ is a sequence of positive integers and $k$ is in $\{1,2\}$
such that $x_1+ x_2+x_3+k \le 6$. There are 14 such pairs,
and a computer assisted computation shows that these equations indeed define
a space of solutions of dimension $6$; see the companion CoCalc worksheet\footnote{Available at \href{https://cocalc.com/share/public_paths/9e9de75dd492cec9240a750f91869a3059fd7661}{\texttt{https://cocalc.com/share/public\_paths/9e9de75dd492cec9240a750f91869a3059fd7661}}; alternatively, a pdf version has been embedded in this file (\openfilelink{Cocalc_Reparametrization.pdf}{\underline{\textcolor{blue}{clickable link}}}; might not work with all pdf readers, but works, e.g., if your pdf is open in firefox).}.

\vupe{It remains to prove that the restriction map from  $\cS^3_{\ell+1}$ to $\cS^3_{\ell}$
is injective for any $\ell \ge 6$ (the last argument of the proof of Proposition~\ref{prop:n2} is valid for any $n$).}
For this, let us consider a solution $c$ of \eqref{eq:recurrence_truncated}
such that $c(x_1,x_2,x_3)=0$ on all triples $(x_1,x_2,x_3)$ of sum at most $\ell$.

We fix a triple of indices $\bm i=(i_1,i_2,i_3)$ of sum $\ell-A$, where $A$ is a ``depth'' parameter, to be fixed later.
The idea is to keep in the system \eqref{eq:recurrence_truncated} the equations involving only variables $(x_1,x_2,x_3)$ of degree $\ell+1$ which are above (coordinate-wise) the triple $(i_1,i_2,i_3)$.
These equations are precisely those indexed by triples $\bm x=(x_1,x_2,x_3)=
 \bm i +\bm a$, where $\bm a$ is a triple $(a_1,a_2,a_3)$ of nonnegative integers
 of sum $A-1$ (if $k=2$) or $A$ (if $k=1$).
 On the other hand, the variables are indexed by triples $\bm x=(x_1,x_2,x_3)=\bm i +\bm a$, where $\bm a$ is a triple $(a_1,a_2,a_3)$ of nonnegative integers
 of sum $A+1$.
 If the corresponding system of equations is invertible for all $(i_1,i_2,i_3)$,
  then vanishing of coefficients $c(x_1,x_2,x_3)$
 for $x_1+x_2+x_3 \le \ell$ will imply the vanishing of coefficients $c(x_1,x_2,x_3)$
 for $x_1+x_2+x_3 = \ell+1$. This will prove that the restriction from  $\cS^3_{\ell+1}$ to $\cS^3_{\ell}$
is an injective map and prove Proposition~\ref{prop:n3}.
 So we only need to check that the system of equations is indeed invertible.
 
Taking $A=2$ gives only 9 equations for only 10 variables, 
so we have no chance to get an invertible system.
On the opposite, taking $A=3$ gives 16 equations for 15 variables, so it is enough to check,
that for all positive $(i_1,i_2,i_3)$, at least one of the principal minors is nonzero.
With the help of a computer, we observe that all principal minors of the system can be factorized into degree $1$ factor with only positive coefficients in $i_1$, $i_2$ and $i_3$, so that in fact all minors never vanish for positive $i_1$, $i_2$ and $i_3$; see again the companion Cocalc worksheet (where we computed only one minor to save some computation time).
This proves the proposition.
\end{proof}

\subsubsection{Case $n=4$}
We cannot prove Conjecture~\ref{conj:truncated_solutions} for $n=4$, but we collected some evidence in favour of it.
First the conjecture holds for $\ell=10$.
In this case, we have a system of 210 variables and 196 equation,
and the computer confirms that the space of solutions has dimension 24.

For the recursive step, we proceed as for $n=3$. 
Taking $A=6$ gives 140 equations for 120 variables.
Computing a minor of size 120,
keeping $i_1$, $i_2$, $i_3$ and $i_4$ as formal variables, seems too demanding for the computer,
but we can look at particular values of $i_1$, $i_2$, $i_3$ and $i_4$.
In particular, while the first minor (i.e.~taking the first 120 rows of the matrix) is sometimes singular\footnote{This is different from the case $n=3$, where all minors are invertible for all positive vectors $(i_1,i_2,i_3)$.},
we have checked that the sixth minor (taking lines 6 to 125 of the matrix) is invertible for all $(i_1,\dots,i_4)$ summing to at most 17, and for all samples of random integers between 1 and 1000 that we have generated.
Again these computations are shown in our Cocalc worksheet (though only up to sum 10, to save computation time and output space).
This gives a strong evidence that the corresponding system has a unique solution,
which would imply Conjecture~\ref{conj:truncated_solutions} for $n=4$.

\subsection{Second approach: operator calculus}
We propose a second approach to Conjecture~\ref{conj:truncated_solutions}.
Let $\mathcal F$ be the space of function from $\Zp^n$ to $\CC$ and $S_i$ be the shift operator in the $i$-th variable, which is an operator from $\mathcal F$ to $\mathcal F$ defined by
\[ S_i c(x_1,\dots x_n) =c (x_1,\dots, x_{i-1}, x_i +1, x_{i+1}, \dots, x_n).\]
The multiplication by $x_i$ is also an operator from $\mathcal F$ to $\mathcal F$,
which we will simply denote by $x_i$.
Note that $x_i$ and $S_j$ do not commute: in fact, for any function $c$, we have
\begin{multline*}
[S_i,x_i] c= S_i x_i c - x_i S_i c  = (x_i+1) c (x_1,\dots, x_{i-1}, x_i +1, x_{i+1}, \dots, x_n)
- x_i c (x_1,\dots, x_{i-1}, x_i +1, x_{i+1}, \dots, x_n)\\
= c (x_1,\dots, x_{i-1}, x_i +1, x_{i+1}, \dots, x_n)=S_i c,
\end{multline*}
i.e.~$[S_i,x_i] =S_i$. On the other hand, operators $S_i$ and/or $x_j$ associated with different indices do commute. The (non-commutative) subalgebra of operators from $\mathcal F$ to $\mathcal F$ spanned by $S_1,\dots,S_n$ and $x_1,\dots,x_n$ is sometimes referred to as the Ore algebra \cite{Ore}. 

Now, eq.~\eqref{eq:recurrence_truncated} means that $R_k c=0$ on $\Omega^n_{\ell-k}$, where
$R_k = \sum x_i - \sum x_i S_i^k$. This implies that $D c=0$ on $\Omega^n_{\ell-d}$,
where $D$ is any operator of degree $d$ living in the left ideal $\mathcal I$ generated by $\{R_k, k \ge 1\}$ (in fact, because of Remark~\ref{rmk:commutators}, this is the same as the left ideal generated by only $R_1$ and $R_2$).
The strategy to analyze the solutions of~\eqref{eq:recurrence_truncated} consists in finding
nice elements in $\mathcal I$. We illustrate it in the case $n=2$.
\begin{lemma}\label{lem:OreCalculus}
There exists polynomials $P_0$, $P_1$, $P_2$ in $x_1$ and $x_2$ such that $P_2(x_1,x_2) \ne 0$ for any $x_1,x_2$ in $\Zp$, and the element
\[D \coloneqq P_2(x_1,x_2) S_2^2 + P_1(x_1,x_2) S_2 + P_0(x_1,x_2)\]
is in $\mathcal I$.
\end{lemma}
\begin{proof}
Setting $\overline{R_1} = x_1 (S_1-1) - x_2 (S_2-1)$, we can check, either by hand or using the \texttt{OreAlgebra} function of the \texttt{HolonomicFunctions} mathematica package~\cite{HolonomicFunctions},
that
\[\overline{R_1} R_1 -(x_1+1) R_2 -R_1=P_2(x_1,x_2) S_2 + P_1(x_1,x_2) S_1 + P_0(x_1,x_2),\]
with
\[P_0(x_1,x_2)=-x_1\, x_2-x_2^2, \quad P_1(x_1,x_2)=2x_2 +2x_1\, x_2+ 2x_2^2, \quad 
 P_2(x_1,x_2)=-2x_2 -x_1\, x_2- x_2^2.\]
 Since $\overline{R_1} R_1 -(x_1+1) R_2 -R_1$ is in $\mathcal I$,
 and $P_2(x_1,x_2)$ is indeed nonzero for positive $x_1$ and $x_2$,
 this proves the lemma.
\end{proof}

\vupe{As a consequence, we get an alternative proof of the case $n=2$ in Conjecture~\ref{conj:truncated_solutions}:}
\begin{corollary}
Let $n=2$ and $\ell \ge 3$. Assume that $c:\Omega^n_\ell \to \CC$ satisfies \eqref{eq:recurrence_truncated}
and that $c(1,1)=c(1,2)=0$ are equal to $0$.
Then $c$ is identically zero.
\end{corollary}
\begin{proof}
Take $P_0$, $P_1$ and $P_2$ be as in Lemma~\ref{lem:OreCalculus}.
Since $D$ is in $\mathcal I$, we have $Dc\equiv 0$ for any solution $c$ of \eqref{eq:recurrence_truncated}, i.e.~
\[P_2(x_1,x_2) c(x_1,x_2+2) + P_1(x_1,x_2) c(x_1,x_2+1) + P_0(x_1,x_2)c(x_1,x_2)=0.\]
In particular, specializing to $x_1=1$, an immediate recurrence starting from $c(1,1)=c(1,2)=0$ shows that $c(1,x_2)=0$ 
\vupe{for any $x_2 \le \ell-1$ (for larger $x_2$, the pair $(1,x_2)$ is not in $\Omega^n_\ell$) }
-- note that the condition $P_2(x_1,x_2) \ne 0$ is crucial here.

We now proceed by induction on $x_1$. Assume that for a given $x_1\ge 1$,
 $c(x_1,x_2)=0$ \vupe{for any $x_2$ in $\{1,\dots,\ell-x_1\}$}. This holds for $x_1=1$, as proved above.
 Then writing
 \[x_1 c(x_1+1,x_2)= (x_1+x_2) c(x_1,x_2) - x_2 c(x_1,x_2+1)\]
 proves that $c(x_1+1,x_2)=0$ for any \vupe{$x_2$ in $\{1,\dots,\ell-x_1-1\}$}.
 Hence $c(x_1,x_2)=0$ for all $(x_1,x_2)$ in $\Omega^n_\ell$, as claimed.
\end{proof}

Unfortunately, we have not been able to find appropriate elements of $\mathcal I$ to generalize this approach to any $n \ge 2$,
\vupe{even with the help of the computer algebra package \texttt{HolonomicFunctions}~\cite{HolonomicFunctions} computing Gröbner bases
in the non-commutative algebra generated spanned by $S_1,\dots,S_n$ and multiplication by $x_1,\dots,x_n$.}



\section{A generalization to colored permutations}
\label{sec:typeB}
We consider here a variant of our main problem. Let $r \geq 2$ be an integer and let $\lambda \in \CC$
be a primitive $r$-th root of $1$. Recall that we use $[n]$ for the residue of $n$ mod $r$
($r$ being fixed throughout the section).
Let $E=\CC[T]_0$ as in Sections~\ref{sec:dimension_n!} and \ref{sec:algebraic},
and additionally set $F=\CC^\lambda[T]$ be the set of polynomials $P$ such that
$P(\lambda T)=\lambda P(T)$ (this means that $P$ is a linear combination of monomials
$T^n$ with $[n]=1$). Observe that $F$ is not a subalgebra of $E$, but is stable under compositions.
In case $r=2$ and $\lambda=-1$, $F$ is the set of odd polynomials.

We are interested in the space $\LI_n(E,F)$ of multilinear functions $\varphi$
satisfying the following property: 
for all $(f_1,\dots,f_n)$ in $E^n$, for all $\gamma$ in $F$,
\begin{equation}\label{e:reparametrage_odd}
 \varphi(f_1\circ \gamma,\dots,f_n \circ \gamma)=\varphi(f_1,\dots,f_n)\circ \gamma.
\end{equation}
Note that the difference with \eqref{e:reparametrage} is that we only require
here this identity to hold for $\gamma$ in the smaller space $F$, and not for all $\gamma$ in $E$.

Naturally, we have $\LI_n(E) \subseteq \LI_n(E,F)$. 
The inclusion turns out to be strict: for example for $n=1$,
the functional $\varphi$ defined by $\varphi(f(T))=f(\lambda T)$ is in $\LI_n(E,F)$ but not in $\LI_n(E)$.
Building on this, we can construct the following elements $\varphi_{\sigma,\bm \eps}$ of $\LI_n(E,F)$, indexed by pairs $(\sigma,\bm \eps)$,
where $\sigma$ is in $\Sym_n$ and $\bm \eps=(\eps_1,\dots,\eps_n)$ is a tuple of elements in $\Z/r\Z$:
\begin{equation}
\label{eq:def_var_sigma_eps}
\varphi_{\sigma,\bm \eps}(f_1,\dots,f_n)=[ f_{\sigma(1)}(\lambda^{\eps_{1}} T)\, :\, \dots\, :\, f_{\sigma(n)}(\lambda^{\eps_{n}} T)].
\end{equation}

The main theorem of this section is the following analog of Theorem~\ref{t:linear}.
\begin{theorem}\label{thm:type_B}
A basis of $\LI_n(E,F)$ is provided by the elements $\varphi_{\sigma,\bm \eps}$, for 
$\sigma$ in $\Sym_n$ and $\bm \eps$ in $(\Z/r\Z)^n$.
\end{theorem}
In particular, $\LI_n(E,F)$ has dimension $r^n n!$.
The polarization argument goes through without modification,
and we also get a description of the set $\AI_n(E,F)$ of 
algebraic solutions of \eqref{e:reparametrage_odd}.
\begin{corollary}\label{cor:type_B}
A basis of $\AI_n(E,F)$ is provided by the elements $\varphi_{\sigma,\bm \eps}$
defined in~\eqref{eq:def_var_sigma_eps}, where \vupe{ $m$ is a positive integer,
$\sigma$ any function} from $\{1,\dots,m\}$ to $\{1,\dots,n\}$
and $\bm \eps$ in $(\Z/r\Z)^m$.
\end{corollary}

The proof of Theorem~\ref{thm:type_B} follows the same line as that of Theorem~\ref{t:linear}. First we need to prove the linear independence of the functions $\varphi_{\sigma,\bm \eps}$.
\begin{lemma}
The functionals $(\varphi_{\sigma,\bm \eps})_{\sigma,\bm \eps}$, for 
$\sigma$ in $\Sym_n$ and $\bm \eps$ in $(\Z/r\Z)^n$, are linearly independent.
\end{lemma}
\begin{proof}
Again we will look at the images of monomial. We have
\[\varphi_{\sigma,\bm \eps}(T^{i_1},\ldots,T^{i_n}) = \frac{\prod_k \lambda^{\eps_ki_{\sigma(k)}}}{(i_{\sigma(1)}+\cdots+i_{\sigma(n)})\cdots(i_{\sigma(n-1)}+i_{\sigma(n)})i_{\sigma(n)}}
T^{\sum_j i_j}.\]
Hence, 
assuming that we have a linear relation $\sum a_{\sigma,\bm \eps} \varphi_{\sigma,\bm \eps} = 0$, applying it to $(T^{i_1},\ldots,T^{i_n})$ we find
$$
\sum_{\sigma,\bm \eps} \frac{a_{\sigma,\bm \eps} \prod_k \lambda^{\eps_ki_{\sigma(k)}}
\prod_k i_k}{(i_{\sigma(1)}+\cdots+i_{\sigma(n)})\cdots(i_{\sigma(n-1)}+i_{\sigma(n)})i_{\sigma(n)}} = 0\,.
$$

We obtain
\begin{equation}
\label{e:comb-lin}
\sum_{\sigma,\beps} a_{\sigma,\bm \eps} \frac{\lambda^{\sum_k \eps_k i_{\sigma(k)}}}{(i_{\sigma(1)}+\cdots+i_{\sigma(n)})\cdots(i_{\sigma(n-1)}+i_{\sigma(n)})i_{\sigma(n)}} = 0\,.
\end{equation}

Equation \eqref{e:comb-lin} holds for all positive integers $i_1,\dots,i_n$.
For each $k$ such that $1 \leq k \leq r$, we fix $\mu_k \in \Z/r\Z$, and restrict
\eqref{e:comb-lin} to lists such that $[i_k] = \mu_k$. For any such list, we get
\[
\sum_{\sigma} \sum_{\beps} a_{\sigma,\bm \eps} \lambda^{\sum_k \eps_k \mu_{\sigma(k)}}
\frac{1}{(i_{\sigma(1)}+\cdots+i_{\sigma(n)})\cdots(i_{\sigma(n-1)}+i_{\sigma(n)})i_{\sigma(n)}} = 0\,.
\]
As in the proof of Lemma~\ref{l:phi-sigma_ind},
this implies that the rational function
\[
\sum_{\sigma} \sum_{\beps} a_{\sigma,\bm \eps} \lambda^{\sum_k \eps_k \mu_{\sigma(k)}}
\frac{1}{(x_{\sigma(1)}+\cdots+x_{\sigma(n)})\cdots(x_{\sigma(n-1)}+x_{\sigma(n)})x_{\sigma(n)}}
\]
is identically $0$, and thus that, for all $\sigma \in \fS_n$,
\[ \sum_{\beps} a_{\sigma,\bm \eps} \lambda^{\sum_k \eps_k \mu_{\sigma(k)}}=0 .\]
Fix $\sigma$ in $S_n$, and consider the above equation for each $n$-tuple $(\mu_k) \in (\Z/r\Z)^n$.
This yields a linear system
with $r^n$ equations and the $r^n$ variables $(a_{\sigma,\beps})_{\beps \in (\Z/r\Z)^n}$.
Setting $\bm \nu = (\nu_k)_{1 \leq k \leq n}$ with $\nu_k=\mu_{\sigma(k)}$,
the matrix of this system has coefficients
$(\lambda^{\scal{\bm \nu,\beps}})_{\bm \nu,\beps}$. This is the $n$-th tensorial power of the
matrix $(\lambda^{kl})_{1\leq k,l\leq n}$ obtained for $n=1$. The Van der Monde matrix
$(\lambda^{kl})_{1\leq k,l\leq n}$ has a non-vanishing determinant, and so also its $n$-th
tensorial power.

Hence, it is invertible, and we can conclude that $a_{\sigma,\beps}=0$, for all $\beps$ in
$(\Z/r\Z)^n$. Since this holds for all permutations $\sigma$, we have proved that all coefficients
$a_{\sigma,\beps}$ must be $0$, i.e.~$ \varphi_{\sigma,\beps}$ is an independent family.
\end{proof}
It remains to prove that the dimension of $\LI_n(E,F)$ is at most $r^n n!$.
First observe that Lemma~\ref{l:form_multilinear_solution} still holds in our context,
i.e.~any $\varphi$ in $\LI_n(E,F)$ maps $(T_{x_1},\dots,T_{x_n})$ to $c(x_1,\dots,x_n)T^{x_1+\dots+x_n}$. Moreover the underlying function $c$ must be $r$-homogeneous and $r$-polyharmonic, in the sense of Definition~\ref{d:r-homogene}.
The proofs of these properties are straight-forward modifications of the corresponding statements in
Section~\ref{sec:dimension_n!}. 
Now, if we fix $\bm \eps$ in $(\Z/r\Z)^n$, Proposition \ref{p:bound} implies that the restriction of
$c$ to $\Ze$ belongs to a space of dimension at most $n!$\,. Since we have $r^n$ choices for $\bm \eps$,
the dimension of the space of functions $c$ is at most $r^nn!$. Theorem \ref{thm:type_B} is proved.

\section{The case of infinitely differentiable functions}
\label{sec:CInfty}
Throughout this section, we let $E=C^\infty(\mathbb R)_0$ denote the space of infinitely differentiable  functions
$\mathbb R\to\mathbb C$ vanishing at $0$.
Recall that $\I_n(E)$, resp.~$\LI_n(E)$, denote the space of solutions, resp.~multilinear solutions of~\eqref{e:reparametrage} with $n$ arguments. 
We are not able to obtain a full classification of such solutions,
but we get some partial results under various continuity assumptions.
Note that some of the results below apply to all $\varphi \in \I_n(E)$,
while other use the additional multilinearity assumption, i.e.~assume $\varphi \in \LI_n(E)$.

\subsection{General results without continuity assumptions}

\begin{proposition}[Polynomials are stable]
\label{p:stabilisation}
Let $E=\CR$ and $G \coloneqq \CC[X]_0$. Then for all $\varphi \in \LI_n(E)$ it holds
$$\varphi(G^n)\subset G.$$
\end{proposition}

\begin{proof} Since  $\varphi$ is multilinear on $G$,  it is completely determined by the values of
$\varphi(X^{i_1}, \dots, X^{i_n})$. The invariance by $\gamma$ such that $\gamma(t) \coloneqq \tau t$ provides
\[
\varphi((\tau X)^{i_1} , \dots , (\tau X)^{i_n})(t)=\tau^{\sum_{k=1}^n i_k}\varphi(X^{i_1}, \dots, X^{i_n})(t) = \varphi(X^{i_1}, \dots, X^{i_n})(\tau t).
\]
Taking  $t=1$ yields
\[
\tau^{\sum_{k=1}^n i_k}\varphi(X^{i_1}, \dots, X^{i_n})(1)= \varphi(X^{i_1}, \dots, X^{i_n})(\tau),
\]
showing that  $\varphi(X^{i_1}, \dots, X^{i_n})(t)$.
\end{proof}

\begin{proposition}[Dependance only on the past]
\label{p:past}
 Let $f_1,\ldots,f_n,g_1,\ldots,g_n \in \CR$, $\varphi \in \I_n(E)$ and let $t_0>0$. Assume that,
 for all $t \in [0,t_0]$ and all  $i \in \{1,\ldots,n\}$, we have $f_i(t)=g_i(t)$. Then
 $\varphi(f_1,\ldots,f_n)(t) = \varphi(g_1,\ldots,g_n)(t)$ for
 $0 \leq t \leq t_0$.
\end{proposition}
\begin{proof}
Let $\gamma \in \CR$ be such that $\im(\gamma) = [0,t_0]$.
Then the hypothesis implies $f_i \circ \gamma=g_i \circ \gamma$ for $i \in \{1,\ldots,n\}$.
Thus $\varphi(f_1 \circ \gamma,\ldots,f_n \circ \gamma)=\varphi(g_1\circ \gamma,\ldots,g_n \circ \gamma)$. Since
$\varphi \in I_n(E)$, we deduce $$\varphi(f_1,\ldots,f_n) \circ \gamma = \varphi(g_1,\ldots,g_n) \circ \gamma.$$ Since
$\im(\gamma)=[0,t_0]$, we obtain $\varphi(f_1,\ldots,f_n)(t) = \varphi(g_1,\ldots,g_n)(t)$ for $t \in [0,t_0]$.
\end{proof}

\begin{remark}
 Of course, we have a similar result for $t_0<0$: if for all $t \in [t_0,0]$
 and all  $i \in \{1,\ldots,n\}$, we have $f_i(t)=g_i(t)$, then
 $\varphi(f_1,\ldots,f_n)(t) = \varphi(g_1,\ldots,g_n)(t)$ for
 $t_0 \leq t \leq 0$.
\end{remark}

For our last result in this section, we assume that $n=2$.
\begin{proposition}[Disjoint supports]
\label{p:disjoint-support}
 Let $f_1,f_2$ be in $E$ satisfying $\Supp(f_1) \subset ]-\infty,a]$ and
 $\Supp(f_2) \subset [b,+\infty[$ with $a<b$. Then, for any $\varphi$ in $\I_2(E)$, the following holds.
 \begin{enumerate}
  \item If $b>0$, we have $\varphi(f_1,f_2)_{|]-\infty,b]}
  = \varphi(f_1,0)_{|]-\infty,b]}$\,;
  \item If $a<0$, we have $\varphi(f_1,f_2)_{|[a,+\infty[}
  = \varphi(0,f_2)_{|[a,+\infty[}$\,.
 \end{enumerate}
 Consequently, if $\varphi$ is multilinear and if $a<0<b$, then $\varphi(f_1,f_2)=0$.
\end{proposition}

\begin{proof}
We consider the first case, with $b>0$. Let $\epsilon>0$ such that $b-2\epsilon>a$ and $b-2\epsilon>0$. Let $\gamma \in \CR$
be such that
$$
\left \{
\begin{array}{l}
 \gamma(t)=t \mbox{ for } t \leq b-2\epsilon \mbox{, and}\\
 b-2\epsilon \leq \gamma(t) \leq b-\epsilon \mbox{ for } t \geq b-2 \epsilon.
\end{array}
\right .
$$
Then
\begin{itemize}
 \item[$\bullet$] $\gamma(0)=0$.
 \item[$\bullet$] $f_2 \circ \gamma = 0$. Indeed, $\forall t \in \R$, $\gamma(t) \leq b-\epsilon$, and $\Supp(f_2) \subset [b,+\infty[$.
 \item[$\bullet$] $f_1 \circ \gamma = f_1$. Indeed, for $t \leq b-2\epsilon$, $\gamma(t)=t$, so $f_1 \circ \gamma(t)=f_1(t)$, and for
 $t>b-2\epsilon$, $\gamma(t)>a$ so $f_1 \circ \gamma(t)=0=f_1(t)$.
\end{itemize}
\medskip

We obtain $\varphi(f_1,f_2)\circ \gamma=\varphi(f_1 \circ \gamma,f_2 \circ \gamma)=\varphi(f_1,0)$. For $t \leq b-2 \epsilon$,
this implies $\varphi(f_1,f_2)(t)=\varphi(f_1,0)(t)$, since $\gamma(t)=t$. Since $\epsilon$ can be arbitrary small, we obtain
the equality $\varphi(f_1,f_2)(t)=\varphi(f_1,0)(t)$ for $t < b$. By continuity, this holds for $t \leq b$. This proves the first item in the lemma.

The case $a<0$ is symmetric. The last statement, assuming that $\varphi$
is multilinear, follows by noticing that in this case $\varphi(f_1,0)=\varphi(0,f_2)=0$.
\end{proof}


\subsection{Path-continuity}
\label{ss:path-cont}
In this section, we assume that $\varphi$ is \emph{path-continuous},
in the following sense:

\begin{definition}
\label{d:path-cont}
We say that $\varphi$ is \emph{path-continuous} if it satisfies the following continuity
property. Let $f_1,\dots,f_n\colon\mathbb R\times\mathbb R\to\mathbb R$,
$(s,t)\mapsto f_i(s,t)$, be such that
\begin{itemize}
  \item for every $s$, the function $f_i(s,\cdot)$ lies in $E$ (i.e. is $C^\infty$ in
        $t$ and vanishes at $t=0$);
  \item for every $t$, the function $s\mapsto f_i(s,t)$ is continuous.
\end{itemize}
Then for every $t$, the function
$s\mapsto\varphi\big(f_1(s,\cdot),\dots,f_n(s,\cdot)\big)(t)$ is continuous.
\end{definition}
All the Landau symbols $O$ throughout the section are considered around $t=0$.
\begin{lemma}
\label{l:finite-det}
Let $\varphi\in \LI_n(E)$ be path-continuous. If $g\in E$ and $g=O(t^{N+1})$, then
for all $f_2,\dots,f_n\in E$ we have $\varphi(g,f_2,\dots,f_n)=O(t^{N+1})$. The same holds in any of the $n$ entries.
\end{lemma}

\begin{proof}
Write $g(t)=t^{N+1}h(t)$ with $h\in C^\infty(\mathbb R)$. Applying~\eqref{e:reparametrage}
to the reparametrization $\gamma=\tau\,\mathrm{id}\in E$ and evaluating at $t=1$ gives,
for every $\tau\in\mathbb R$,
\[
  \varphi(g,f_2,\dots,f_n)(\tau)
  =\varphi\big(g(\tau\,\cdot),f_2(\tau\,\cdot),\dots,f_n(\tau\,\cdot)\big)(1).
\]
Since $g(\tau t)=\tau^{N+1}\,u_\tau(t)$ with $u_\tau(t) \coloneqq t^{N+1}h(\tau t)\in E$, linearity
in the first entry yields
\[
  \varphi(g,f_2,\dots,f_n)(\tau)
  =\tau^{N+1}\,\varphi\big(u_\tau,\,f_2(\tau\,\cdot),\dots,f_n(\tau\,\cdot)\big)(1).
\]
The family $\tau\mapsto\big(u_\tau,f_2(\tau\,\cdot),\dots,f_n(\tau\,\cdot)\big)$ is
$C^\infty$ in $t$ and continuous in $\tau$ (including at $\tau=0$). By path-continuity,
$\tau\mapsto\varphi(u_\tau,\dots)(1)$ is continuous, hence bounded by some constant $C$
on $[-1,1]$. Therefore $|\varphi(g,f_2,\dots,f_n)(\tau)|\le C\,|\tau|^{N+1}$ for
$|\tau|\le 1$.
\end{proof}

\begin{proposition}\label{p:jet-classification}
Let $\varphi\in \LI_n(E)$ be path-continuous. Then there exist scalars $(a_\sigma)_{\sigma\in S_n}$
such that, for every $F=(f_1,\dots,f_n)\in E^n$, the smooth functions $\varphi(F)$ and
$\sum_\sigma a_\sigma\varphi_\sigma(F)$ have the same Taylor expansion at $0$
\vupe{at all orders}.
\end{proposition}

\begin{proof}
By Proposition \ref{p:stabilisation}, $\varphi$ restrict to a functional
$\CC[T]_0^n \to \CC[T]_0$, and
by Theorem~\ref{t:linear} there exist coefficients $(a_\sigma)_{\sigma \in \fS_n}$ such that this
restriction is equal to $\sum a_\sigma \varphi_\sigma$. Set
$\psi \coloneqq \varphi-\sum a_\sigma \varphi_\sigma \in \LI_n(E)$: it is path-continuous and vanishes on polynomials.

Fix $F\in E^n$ and $N\ge 0$. For $1\le i\le n$ write $P_i \coloneqq P^{(i)}_N$ for the Taylor
polynomial of $f_i$ of degree $N$, and $r_i \coloneqq f_i-P_i$, so that $r_i=O(t^{N+1})$. Changing the arguments one at a time, set
$A_i \coloneqq \psi(f_1,\dots,f_i,P_{i+1},\dots,P_n)$ for $0\le i\le n$, so that
$A_0=\psi(P_1,\dots,P_n)$ and $A_n=\psi(F)$. Since $f_i=P_i+r_i$, linearity in the
$i$-th slot gives
\[
  A_i-A_{i-1}
  =\psi(f_1,\dots,f_{i-1},f_i-P_i,P_{i+1},\dots,P_n)
  =\psi(f_1,\dots,f_{i-1},r_i,P_{i+1},\dots,P_n),
\]
and summing this identity over $i$ yields
\[
  \psi(F)=\psi(P_1,\dots,P_n)
  +\sum_{i=1}^n\psi(f_1,\dots,f_{i-1},\,r_i,\,P_{i+1},\dots,P_n).
\]
The first term vanishes because $\psi$ vanishes on polynomials. For $i\le n$, the $i$-th term
of the sum is $O(t^{N+1})$, by Lemma~\ref{l:finite-det} applied to the $i$-th entry. Hence
$\psi(F)=O(t^{N+1})$ for every $N$, i.e. $\psi(F)$ has vanishing Taylor expansion at $0$.
\end{proof}

\subsection{Characterization of $\LI_n(E)$ under a global continuity assumption}
 \label{ss:global-continuity}
For $f \in \CR$,  $\ell \in \Zp$ and $K\subset \R$ compact, we denote by $\|f\|_{C^{\ell}(K)}$
the usual norm on $C^\ell(K)$, namely
\[
\|f\|_{C^{\ell}(K)} \coloneqq  \sum_{k=0}^{\ell} \max_{x\in K} |f^{(k)}(x)|.
\]

\begin{definition}
We say that $\varphi$ is \emph{globally continuous} when for all $t\in \R$, there exists $C>0$,
$\ell\in \Zp$ (all depending on $t$) such that for all $(f_1,\dots , f_n) \in E^n$,
\begin{eqnarray}
|\varphi(f_1,\dots,f_n)(t)| \leq C \prod_{k=1}^n \|f_k\|_{C^{\ell}([-t,t])}. \label{continuity}
\end{eqnarray}
\end{definition}

We will show that the only globally continuous elements of $LI_n(E)$ are the linear combinations of
$\varphi_\sigma$. We start with a lemma.

\begin{lemma}[Approximation by polynomials]
\label{approximation}
Let $f\in \CR$. Then for all $\ell\in \Zp$ and $t\in \R$, one can find a sequence $(g_k)_{k\geq 0}$ of polynomials such that  $g_k(0)=0$ for all $k$ and
$$\|f-g_k\|_{C^{\ell}([-t,t])}\underset{k\to +\infty}{\longrightarrow} 0.$$ 
\end{lemma}
\begin{proof} We denote   $K \coloneqq [-t,t]$. Since $f^{(\ell)}$ is a continuous function on the compact set $K$, by the famous Weierstrass Theorem we know that there exists a sequence of polynomials $h_k\in \R[T]_0 $  such that $\|h_k - f^{(\ell)}\|_{L^\infty(K)} \to_{k\to +\infty} 0$. Then   for $s\in K$ we define 
$$q_k(s)  \coloneqq  f^{(\ell-1)}(0)+ \int_0^s h_k(x) dx. $$
The sequence  $(q_k)_{k\geq 0}$ is a sequence of polynomials  that converges uniformly on  $K$ to the function
$f^{(\ell)}$, and moreover $q_k(0)=f^{(\ell-1)}(0)$.  If $\ell=1$ then the proof is over. Otherwise, we continue 
integrating and after taking $\ell$ primitives we get a sequence $g_k$ converging to $f$ in the norm
$\|\cdot\|_{C^{\ell}(K)}$ with the property that $g_k(0)=f(0)=0$, as desired.
This proves of the Lemma.
\end{proof}

\begin{proposition} Let  $E=\CR$ and $\varphi \in \LI_n(E)$ be globally continous.
Then $\varphi \in Vect\{\varphi_\sigma, \; \sigma \in \fS_n \}$.
\end{proposition}
\begin{proof} Let $\varphi \in \LI_n(E)$ satisfying the continuity assumption \eqref{continuity}
of the Proposition.
We first restrict $\varphi$ to the subspace $G^n \subset E^n$, where $G \coloneqq \CC[X]_0$ is the space of polynomials.
Applying Proposition \ref{p:stabilisation} we deduce that $\varphi \in \LI_n(G)$. Then by
Theorem~\ref{t:linear} we know that the restriction of  $\varphi$ on $G^n$ is a linear combination of
$\varphi_\sigma$, for $\sigma \in \fS_n$.
In other words there exists some constants $a_\sigma$ such that for all
$(g_1,\cdots, g_n)\in G^n$,
\begin{equation}
\varphi(g_1,\cdots, g_n)=\sum_{\sigma \in \fS_n} a_\sigma \varphi_\sigma(g_1,\cdots , g_n). \label{charact}
\end{equation}

The idea is to extend this on $E^n$ by a density argument.
First, it follows from Lemma~\ref{l:iterated-simplex}
that $\varphi_\sigma$ is globally continuous, where we can choose
uniformly $\ell=1$ and $C=|t^{n}|$ in the definition of global continuity.

Let $t>0$ be fixed and let $C>0$, $\ell\in \Zp$ (both depending on $t$) such that for all $(f_1,\dots , f_n) \in
E^n$  the inequality \eqref{continuity} holds. Let $(f_1,\dots , f_n) \in E^n$ be fixed.  For each $1\leq i \leq n$
we consider a sequence  $(g_{i,k})_{k\geq 0}$ of polynomials given by Lemma~\ref{approximation} such that $g_{i,k}$
converges to $f_i$ for the norm $\|\cdot \|_{C^{\ell}([-t,t])}$.
Then the continuity assumption \eqref{continuity} implies that  
$$\varphi(g_{1,k},\cdots, g_{n,k})(t) \underset{k\to +\infty}{\longrightarrow}  \varphi(f_1,\cdots, f_n)(t).$$
On the other hand, convergence in norm $C^\ell$ implies
convergence holds in norm $C^1$ and it follows that
for any $\sigma \in \fS_n$,
$$\varphi_\sigma(g_{1,k},\dots,g_{n,k})(t) = [g_{\sigma(1),k}:g_{\sigma(2),k}:\cdots:g_{\sigma(n),k}](t) \underset{k\to +\infty}{\longrightarrow}  [f_{\sigma(1)}:f_{\sigma(2)}:\cdots:f_{\sigma(n)}](t)=\varphi_\sigma(f_1,\dots,f_n).$$
Using~\eqref{charact} for $(g_{1,k},\dots,g_{n,k})$, and letting $k$ goes to infinity, we obtain
\[\varphi(f_1,\cdots, f_n)(t)=\sum_{\sigma \in \fS_n} a_\sigma \varphi_\sigma(f_1,\cdots , f_n)(t).\]
Since this holds for all $t$ in $\mathbb R$, this finishes the proof of the proposition.
\end{proof}
 
\section{Other toy models}
  
\subsection{A group-theoretic solution}
\label{sec:everything_invertible}

In this section, we consider the case where $E$ is as follows:

\begin{definition}
\label{d:quasi-group}
We say that $(E,\circ)$ is a \emph{quasi-group} when $E$ has an identity $\mathrm{id}$
and a left-absorbing element $0$
(so $0\circ f=0$ for all $f$), and when every $f\in E\setminus\{0\}$ is invertible
for $\circ$.
\end{definition}

For $1\le i\le n$, let $\xi_i\colon E^{i-1}\to E^n$ be defined by
$\xi_i(g_1,\dots,g_{i-1})=(g_1,\dots,g_{i-1},\mathrm{id},0,\dots,0)$.

\begin{proposition}
\label{p:quasi-group}
Assume that $(E,\circ)$ is a quasi-group with $|E| \geq 3$. Then
\[
  \varphi\longmapsto(\psi_1,\dots,\psi_n),\qquad \psi_i \coloneqq \varphi\circ\xi_i,
\]
is a bijection $\I_n(E)\simeq \prod_{i=1}^n\operatorname{Map}(E^{i-1},E)$,
whose inverse sends $(\psi_1,\dots,\psi_n)$ to the map $\varphi$ given by
$\varphi(0,\dots,0)=0$ and, for $(f_1,\dots,f_n) \ne (0,\dots,0)$,
\begin{equation}
\label{e:quasi-group}
\varphi(f_1,\dots,f_n)=\psi_i\bigl(f_1\circ f_i^{-1},\dots,f_{i-1}\circ f_i^{-1}\bigr)\circ f_i,
\end{equation}
where $i$ is the index of the last nonzero entry of $f$, i.e.~$i=\max\{j:\, f_j\ne 0\}$.
\end{proposition}
\vupe{This statement illustrates that the space of solutions of \eqref{e:reparametrage}
can be huge, if we do not impose the additional algebraicity or multilinearity condition.}
\begin{proof}
First, let us prove that $\varphi(0,\dots,0)=0$ for every $\varphi\in I_n(E)$.
By left-absorption gives $(0,\dots,0)\circ\gamma=(0,\dots,0)$ for all $\gamma\in E\setminus\{0\}$.
Hence, setting $c \coloneqq \varphi(0,\dots,0)$,
equation \eqref{e:reparametrage}
yields $c\circ\gamma=c$. If $c\neq0$ then $c$ is
invertible, and we get $\gamma=\mathrm{id}$.
This should hold for all $\gamma\in E\setminus\{0\}$, contradicting $|E| \geq 3$.
So $c=0$.

Let $\varphi\in I_n(E)$ and put $\psi_i \coloneqq \varphi\circ\xi_i$. Fix
$F=(f_1,\dots,f_n)$ and set $i \coloneqq \max\{j:f_j\neq0\}$. For $j<i$ one has
$f_j=(f_j\circ f_i^{-1})\circ f_i$, while $f_i=\mathrm{id}\circ f_i$ and
$f_j=0=0\circ f_i$ for $j>i$; hence, entrywise,
$F=\xi_i\bigl(f_1\circ f_i^{-1},\dots,f_{i-1}\circ f_i^{-1}\bigr)\circ f_i$. Applying
\eqref{e:reparametrage} with $\gamma=f_i$,
\[
  \varphi(F)
  =\varphi\bigl(\xi_i(f_1\circ f_i^{-1},\dots,f_{i-1}\circ f_i^{-1})\bigr)\circ f_i
  =\psi_i\bigl(f_1\circ f_i^{-1},\dots,f_{i-1}\circ f_i^{-1}\bigr)\circ f_i,
\]
which is the displayed formula; in particular $\varphi$ is determined by
$(\psi_1,\dots,\psi_n)$.

Conversely, let $\psi_i\colon E^{i-1}\to E$ be arbitrary and define $\varphi$ by the
formula \eqref{e:quasi-group}, and $\varphi(0,\dots,0)=0$. Let $\gamma\in E\setminus\{0\}$ and let
$F\neq(0,\dots,0)$ have last nonzero index $i$. As $\gamma$ is invertible,
$f_j\circ\gamma=0\iff f_j=0$, so $F\circ\gamma$ also has last nonzero index $i$.
 Using
$(f_j\circ\gamma)\circ(f_i\circ\gamma)^{-1}=f_j\circ f_i^{-1}$, we can write
\[
  \varphi(F\circ\gamma)
  =\psi_i\bigl(f_1\circ f_i^{-1},\dots,f_{i-1}\circ f_i^{-1}\bigr)\circ(f_i\circ\gamma)
  =\varphi(F)\circ\gamma.
\]
This identity also holds trivial when $F=(0,\dots,0)$; thus $\varphi\in \I_n(E)$.

Finally, let us check that the two assignments are mutually inverse. The formula recovers $\varphi$ from its
$\psi_i$'s by the argument for the first implication. Conversely $\varphi\circ\xi_i$ recovers $\psi_i$,
since $\xi_i(g_1,\dots,g_{i-1})$ has last nonzero index $i$ with $i$-th entry
$\mathrm{id}$, so
$\varphi(\xi_i(g))=\psi_i(g_1\circ\mathrm{id}^{-1},\dots,g_{i-1}\circ\mathrm{id}^{-1})\circ\mathrm{id}
=\psi_i(g_1,\dots,g_{i-1})$.
\end{proof}

\begin{remark}
For $i=1$ one has $\operatorname{Map}(E^0,E)=E$, and \eqref{e:quasi-group} reads
$\varphi(f_1,0,\dots,0)=\psi_1\circ f_1$ with $\psi_1=\varphi(\mathrm{id},0,\dots,0)\in E$.
\end{remark}
\medskip

Now we provide examples of such $E$. One possible example is to take any field $(K,+,\times)$ and
$E=(K,\times)$. It is also possible to take any group $G$, and add to $G$ an absorbing element $0$.
But these examples are quite far from the original motivation, and we will now provide an example
within the framework of Puiseux series.

Let $\Gamma$ be the set of rational numbers with an odd denominator:
\[ \Gamma = \mathbb{Z}_{(2)} = \left\{ \frac{p}{q} \mathrel{\Big|} p, q \in \mathbb{Z}, \, q \text{ is odd} \right\}. \]
$\Gamma$ is a subgroup of $(\mathbb{Q}, +)$ and a totally ordered abelian group.
Let $\mathcal{O}$ be the set of rational numbers where both the numerator and the denominator are odd integers:
\[ \mathcal{O} = \left\{ \frac{p}{q} \in \Gamma \mathrel{\Big|} p \text{ is odd} \right\}. \]

We further consider the field of Puiseux series $K = \bigcup_n \mathbb{R}((T^{1/{(2n+1)}}))$;
note that we only consider reciprocal of odd integers in the generators,
so that all exponents that appear in elements of $K$ are in $\Gamma$.
 Let also $K_{\geq 0}$, resp.~$K_{>0}$, be the subset
 where the {\em exponents} are non-negative, resp.~positive,
and inside $K_{> 0}$ the subset $E$ where all the exponents are in $\cO$.
For $f,g \in K_{\geq 0}$, the composition $f \circ g$ is well-defined
because for an odd integer $q$, any real number has a unique
$q$-th root, as well as any series of the form $1+h$ with $h \in K_{>0}$,
and $(1+h)^{1/q} \in K_{\geq 0}$.

%
%

\begin{proposition}
The set $E$ forms a monoid under composition, and every non-zero element in $E$ is invertible for
the composition.
\end{proposition}
\begin{proof}
As a first step, we formalize the notion of an ``odd'' or ``even'' generalized series by defining an automorphism that mimics
evaluating $f(-T)$.
Define a map $\pi: \Gamma \to \mathbb{Z}_2$ by mapping a fraction $r = p/q$ (in lowest terms, $q$ odd) to its
numerator modulo 2:
\[ \pi\left(\frac{p}{q}\right) = p \pmod 2. \]
Observe that $\pi$ is a group homomorphism.
Now, define a map $\sigma: K \to K$ by:
\[ \sigma\left( \sum_{r \in \Gamma} c_r T^r \right) = \sum_{r \in \Gamma} c_r (-1)^{\pi(r)} T^r. \]
Because $\pi$ is a homomorphism, $\sigma(T^r T^s) = \sigma(T^r)\sigma(T^s)$.
Therefore, $\sigma$ is an $\mathbb{R}$-algebra automorphism of $K$.
Notice that $\pi(r) = 1 \iff r \in \mathcal{O}$. Consequently:
\[ f \in E \iff \sigma(f) = -f .\]
In other words, $E$ is exactly the set of ``odd'' Puiseux series.

As a second step, we show that $E$ is closed under compositions.
Let $f, g \in E$. We want to show that $f \circ g \in E$.
Since composition in $K_{>0}$ is defined, $f(g(T)) = \sum_{r \in \Supp(f)} c_r g(T)^r$.
Apply the automorphism $\sigma$:
\begin{equation}
\label{e:compo-odd}
\sigma(f \circ g) = \sigma\left( \sum_{r \in \Supp(f)} c_r g^r \right)
= \sum_{r \in \Supp(f)} c_r \sigma(g)^r
= \sum_{r \in \Supp(f)} c_r (-g)^r = - \sum_{r \in \Supp(f)} c_r g^r = -(f \circ g)\,.
\end{equation}

Since $\sigma(f \circ g) = -(f \circ g)$, we conclude that $f \circ g \in E$.
Because the identity element is in $E$, $E$ is a monoid.

As the last step, we show that every non-zero element in $E$ is invertible
for the composition.
Let $f \in E \setminus \{0\}$. Let $cT^r$ be the monomial
of lowest degree: we can write $f$ as the following composition of series
in $E$:
\[
  f=(cT) \circ (T+g(T)) \circ (T^r),
\]
where $g(T)$ has valuation at least $2$.
Therefore, using the second step,
it is enough to show that those three series are invertible in
$E$. For $T^r$ and $cT$, this is obvious. For the series $T+g(T)$, we compute
its inverse by observing that $(T+g(T)) \circ h(T)=T$ if and only if
$h(T)=T-g(h(T))$. Thus, we define the inverse $h$ as the limit of the sequence $(h_k)_{k \ge 0}$ 
defined inductively by $h_0(T)=T$ and $h_{k+1}(T)=T-g(h_k(T))$.
By the second step, each series $h_k$ is in $E$, therefore $h$ is also in $E$.
\end{proof}

\subsection{The case of functions with finite range}
\label{sec:piecewise}

 A function $f\colon(\mathbb R,0)\to(\mathbb R,0)$ has finite range if and only if it can
be written as a \emph{finite} sum
\[
  f=\sum_{i=1}^N \alpha_i\,\mathbf 1_{A_i},
\]
where $\alpha_i\neq 0$ are real numbers and the $A_i$ are disjoint subsets of
$\mathbb R$. Such an expression of $f$ is unique. 
The function $f$ vanishes at $0$ if and only if none of the $A_i$'s contains $0$.
The set $E$ of functions with finite range vanishing at $0$ is a subalgebra of
the algebra of functions from $\R$ to $\R$, and it is stable under
composition.

If $\cC$ is a condition, we denote by $\delta_\cC$ the integer $1$ if $\cC$ is true,
and $0$ if $\cC$ is false. We shall only need the following two elementary composition identities.
For a point
$c$ in $\R$ and sets $A,C$ such that $0\notin A$, one has
\begin{equation}\label{e:comp-indic}
  \mathbf 1_A\circ(c\,\mathbf 1_C)=\delta_{c\in A}\,\mathbf 1_C=\mathbf 1_A(c)\,\mathbf 1_C.
\end{equation}
By linearity it follows that, for any $h\in E$,
\begin{equation}\label{e:comp-h}
  h\circ(c\,\mathbf 1_C)=h(c)\,\mathbf 1_C.
\end{equation}
\begin{proposition}\label{p:piecewise}
If $\varphi\in \LI_n(E)$, then there exists $\kappa\in\mathbb R$ such that
\[
  \varphi(f_1,\dots,f_n)=\kappa\, f_1\cdots f_n\qquad\text{for all }f_1,\dots,f_n\in E.
\]
\end{proposition}

\begin{proof}
By multilinearity, $\varphi$ is determined by its values on tuples of indicator
functions, so it suffices to show, that for any subsets $A_1,\dots,A_n$ not containing $0$
 \begin{equation}\label{e:to_prove}
 \varphi(\mathbf 1_{A_1},\dots,\mathbf 1_{A_n}) = \kappa\, \mathbf 1_{\bigcap_{k\le n} A_k}.
\end{equation}

\smallskip
\emph{Step 1: the diagonal.}
Let $A\neq\varnothing$, set $h_A \coloneqq \varphi(\mathbf 1_A,\dots,\mathbf 1_A)$ and pick $c\in A$.
Since $\mathbf 1_A\circ(c\,\mathbf 1_A)=\mathbf 1_A$ by~\eqref{e:comp-indic}, the commutation with reparametrizations
\eqref{e:reparametrage} gives $h_A=h_A\circ(c\,\mathbf 1_A)$.
\vupe{The right-hand side equals $h(c)\,\mathbf 1_A$
by~\eqref{e:comp-h}. Since the left-hand side is independent from $c$ in $A$,
the scalar $ h(c)$ is independent of $c$; let us denote it by $\lambda(A)$.
We have proved that}
\begin{equation}\label{e:diagonal}
  \varphi(\mathbf 1_A,\dots,\mathbf 1_A)=\lambda(A)\,\mathbf 1_A.
\end{equation}
\smallskip

\emph{Step 2: general formula.}
Let now $A_1,\dots,A_n$ be arbitrary, set $h \coloneqq \varphi(\mathbf 1_{A_1},\dots,\mathbf 1_{A_n})$
and $D \coloneqq \bigcap_k A_k$. For any $c$ in $\R$, applying~\eqref{e:reparametrage} to
$\gamma=c\,\mathbf 1_D$ and using~\eqref{e:comp-indic}, \eqref{e:comp-h}, \eqref{e:diagonal} and multilinearity,
\begin{multline*}
  h(c)\,\mathbf 1_D = \varphi(\mathbf 1_{A_1},\dots,\mathbf 1_{A_n}) \circ (c\,\mathbf 1_D)
  =\varphi\big(\mathbf 1_{A_1} \circ (c\,\mathbf 1_D) ,\dots,\mathbf 1_{A_n} \circ (c\,\mathbf 1_D) \big)\\
  =\varphi\big(\delta_{c\in A_1}\,\mathbf 1_D,\dots,\delta_{c\in A_n}\,\mathbf 1_D\big)
  =\Big(\prod_{k=1}^n\delta_{c\in A_k}\Big)\,\varphi(\mathbf 1_D,\dots,\mathbf 1_D)=
  \delta_{c\in D} \lambda(D)\,\mathbf 1_D.
\end{multline*}
This is an equality of functions on $\R$. Specializing it at some element of $D$, we get that
$h(c)= \delta_{c\in D} \lambda(D)$. Since this holds for any $c$ in $\R$,
we get 
\begin{equation}\label{e:general}
h=\varphi(\mathbf 1_{A_1},\dots,\mathbf 1_{A_n})=\lambda(D) \,\mathbf 1_D.
\end{equation}

\emph{Step 3: $\lambda$ is constant.}
Let $A$ and $C$ be subsets of $\R$ containing a common element $c$, we have
\[\lambda(A) \,\mathbf 1_C=\lambda(A) \,\mathbf 1_{A} \circ (c\,\mathbf 1_C)=\varphi(\mathbf 1_{A},\dots,\mathbf 1_{A}) \circ (c \mathbf 1_C)
  =\varphi\big(\mathbf 1_{A} \circ (c\,\mathbf 1_C) ,\dots,\mathbf 1_{A} \circ (c\,\mathbf 1_C) \big)\\
  =\varphi(\mathbf 1_C,\dots,\mathbf 1_C) =\lambda(C) \,\mathbf 1_C.\]
This implies $\lambda(A)=\lambda(C)$. 
For any subsets $A$ and $B$ of $\R$, one can find another set $C$ 
that intersects both $A$ and $B$, and we have $\lambda(A)=\lambda(C)=\lambda(B)$.
Letting $\kappa$ be the common values of all $\lambda(D)$,
Eq.~\eqref{e:general} simplifies to~\eqref{e:to_prove}, which proves the proposition.
\end{proof}

\appendix
\section{Local to global principles}

\subsection{Homogeneity is a consequence of polyharmonicity}
\label{ss:homog-from-induction}
Theorem~\ref{t:linear},
  studies homogeneous polyharmonic functions
$c:\Znp \to \CC$. In this section, we prove that the homogeneity hypothesis is superfluous,
i.e.~it in fact follows from polyharmonicity.
Consequently, any polyharmonic function $c\colon\mathbb Z_{>0}^n\to\mathbb C$
is a linear combination of the functions
 $c_\sigma$ introduced in Definition \ref{d:c-sigma}.

Throughout this subsection, $c\colon\mathbb Z_{>0}^n\to\mathbb C$ is a polyharmonic function, and $\varphi\colon E^n\to E$ ($E=\mathbb C[T]_0$) is the
associated multilinear map, defined by~\eqref{e:def-c}. 
The idea is to prove a kind of converse of Lemma~\ref{lem:derivative}.
For this, it will be convenient to extend our functional $\varphi$ the space
 $\widehat E=\mathbb C[[T]]_0$ of formal power series.
 \vupe{We endow $E$ with the $T$-adic topology, i.e.~$f_n$ converges to $f$ if,
 for any fixed $k$, the sequence $[t^k]f_n$ is eventually equal to $[t^k]f$.}

\subsubsection*{Extension to formal power series.}
For $F=(f_1,\dots,f_n)\in\widehat E^{\,n}$ with
$f_i=\sum_{k\ge1}a_{i,k}T^k$, set
\[
  \widehat\varphi(F) \coloneqq \sum_{d\ge n}\Big(\sum_{\substack{k_1,\dots,k_n\ge1\\ k_1+\cdots+k_n=d}}
  c(k_1,\dots,k_n)\prod_{i=1}^n a_{i,k_i}\Big)\,T^d .
\]
Each coefficient is a finite sum, so $\widehat\varphi\colon\widehat E^{\,n}\to\widehat E$ is a
well-defined multilinear map, continuous for the $T$-adic topology, restricting to $\varphi$ on
$E^n$.

\begin{lemma}\label{l:der-formal}
For all $F\in\widehat E^{\,n}$ and all $\eta\in\widehat E$,
\begin{equation}\label{eq:eta_phihat_prime}
  \eta\,[\widehat\varphi(F)]'=\sum_{i=1}^n\widehat\varphi(f_1,\dots,\eta f_i',\dots,f_n).
\end{equation}
\end{lemma}
\begin{proof}
For monomials $F=(T^{x_1},\dots,T^{x_n})$ and $\eta=T^u$, the identity \eqref{eq:eta_phihat_prime}
reduces exactly to~\eqref{recurrence-c} -- see the computation of~\eqref{amontre1} and~\eqref{identiteDerivee2} in
the proof of Proposition~\ref{p:induction}; \vupe{hence \eqref{eq:eta_phihat_prime} holds in this case}.
By multilinearity in $F$ and linearity in $\eta$ it holds
for all polynomials $F,\eta$. 
Moreover, both sides are $T$-adically continuous in $(F,\eta)$
(the coefficient of $T^d$ is a finite expression in the coefficients of $F$ and $\eta$ of degree at most
$d+1$). 
Since polynomials are $T$-adically dense in $\widehat E$, 
we conclude that \eqref{eq:eta_phihat_prime} holds on all of $\widehat E$.
\end{proof}

For the next result, we will need the following generalization of Leibniz rule from product to any multilinear function. We say that a $s$-dependent element $f(s)$ of $ \widehat E$ is continuously
differentiable in $s$ if  the coefficient of any $T^k$, for fixed $k$, is continuously differentiable in $s$. \vupe{Note that this does not imply that $s \mapsto f(s)$ is continuous for the $T$-adic topology, but only for the weaker topology of convergence of each individual coefficient.}
When $f(s)$ is continuously
differentiable in $s$, we write $\partial_s$ for the coefficient-wise derivative with respect to $s$.
\begin{lemma}[Leibniz rule]\label{lem:Leibniz}
Let $f_1(s),\dots,f_n(s)$ be $s$-dependent elements of $\widehat E$ continuously differentiable in $s$, and $\psi$ be a $T$-adically continuous multilinear function from $\widehat  E^n$ to $\widehat E$. Then we have
\[ \partial_s \big(\psi(f_1(s),\dots,f_n(s)) \big)=\sum_{i=1}^n
\psi(f_1(s),\dots, \partial_s f_i(s), \dots, f_n(s)).\]
\end{lemma}

\begin{proof}
 \vupe{\pe{Let $\delta \neq 0$.} Denoting $F(s)=(f_1(s),\dots, f_n(s))$, by multilinearity,
\begin{multline*}
  \psi(F(s+\delta))-\psi(F(s))
  =\sum_{i=1}^n \psi\Big(f_1(s+\delta),\dots,
     f_i(s+\delta)-f_i(s),\dots,f_n(s)\Big)\\
 =\delta \left(\sum_{i=1}^n \psi\Big(f_1(s+\delta),\dots,
       \partial_s f_i(s),\dots,f_n(s)\Big) +
       \psi\Big(f_1(s+\delta),\dots,
       R_i(\delta),\dots,f_n(s)\Big) \right),
\end{multline*}
where $R_i(\delta)=\delta^{-1}(f_i(s+\delta)-f_i(s))-\partial_s f_i(s)$ tends coefficient-wise to $0$ as $\delta$ tends to $0$.
Since $\psi$ is $T$-adically continuous, the coefficient of $t^k$
in $\psi\Big(f_1,\dots,f_n\Big)$ depends only on finitely many
coefficient of $(f_1,\dots, f_n)$, and hence can be seen as a \pe{multilinear} map on a vector space of finite
dimension. Hence it is continuous and we get that, for each $k$,
\begin{align*}
\lim_{\delta \to 0} [t^k]   \psi\Big(f_1(s+\delta),\dots,
       \partial_s f_i(s),\dots,f_n(s)\Big) &= \psi\Big(f_1(s),\dots,
       \partial_s f_i(s),\dots,f_n(s)\Big),\\
       \lim_{\delta \to 0} [t^k]  \psi\Big(f_1(s+\delta),\dots,
       R_i(\delta),\dots,f_n(s)\Big) &=0.
       \end{align*}
Therefore, $[t^k]\psi(F(s))$ is continuously differentiable in $s$ with derivatives given
by
\[\partial_s( [t^k]\psi(F(s)))=[t^k]\, \sum_{i=1}^n
\psi(f_1(s),\dots, \partial_s f_i(s), \dots, f_n(s)).\]
Since this holds for every fixed $k \ge 1$, this proves the lemma.}
\end{proof}
We can now prove the main result of this section.
\begin{proposition}
\label{p:homog-from-induction}
Let $c\colon\mathbb Z_{>0}^n\to\mathbb C$ satisfy~\eqref{recurrence-c}. Then $c$ satisfies
\[
c(\lambda x_1,\ldots,\lambda x_n)=c(x_1,\ldots,x_n).
\]
for all $\lambda,x_1,\ldots,x_n \in \Zp$,
and consequently $c=\sum_{\sigma\in S_n}a_\sigma c_\sigma$ for some scalars $a_\sigma$.
\end{proposition}

\begin{proof}
Fix $\lambda\in\mathbb Z_{>0}$; we want to prove that
\begin{equation}\label{e:commute-power}
  \widehat\varphi(F\circ T^\lambda)=\widehat\varphi(F)\circ T^\lambda\qquad(F\in\widehat E^{\,n}).
\end{equation}
The case $\lambda=1$ is trivial, so assume $\lambda\ge2$. Set $\eta \coloneqq T^\lambda-T$ and, for
$s\in[0,1]$, $\gamma_s \coloneqq \mathrm{id}+s\eta=(1-s)T+sT^\lambda\in\widehat E$, so that $\gamma_0=\mathrm{id}$ and
$\gamma_1=T^\lambda$. We set
\[
  \Psi(s) \coloneqq \widehat\varphi(F\circ\gamma_s),\qquad \Theta(s) \coloneqq \widehat\varphi(F)\circ\gamma_s.
\]
\vupe{Trivially
 $\Psi(0)=\Theta(0)=\widehat\varphi(F)$, and we would like to establish that $\Psi(1)=\Theta(1)$.}
 By construction, $\Psi$ and $\Theta$ are $s$-dependent power series in $T$, 
\vupe{continuously differentiable in $s$},
 so they can be seen as functions of two variables, and we will use the standard
 notation $\partial_s$ and $\partial_t$ for the corresponding partial derivatives.
 We also continue to use the notation $f'$ for the derivative with respect to $t$
 of power series in $\widehat E$.
 
We now apply the Leibniz rule (Lemma~\ref{lem:Leibniz}) to
the multilinear functional $\widehat\varphi$ and to $f_i(s)=f_i\circ\gamma_s$.
 Using that $\partial_s(f_i\circ\gamma_s)=\eta \cdot (f_i'\circ\gamma_s)$,
we get
\[
  \partial_s\Psi=\sum_{i=1}^n\widehat\varphi\big(f_1\circ\gamma_s,\dots,\eta\,(f_i'\circ\gamma_s),\dots,f_n\circ\gamma_s\big).
\]
On the other hand, for $s\in[0,1)$ one has $\gamma_s'(0)=1-s\neq0$, so $\gamma_s'$ is invertible in $\mathbb C[[T]]$ and
$\eta/\gamma_s'\in\widehat E$. Applying Lemma~\ref{l:der-formal} to the tuple $F\circ\gamma_s$ and
the  function $\eta/\gamma_s'$, 
\begin{equation}\label{e:dtPsi}
  \frac{\eta}{\gamma_s'}\,\partial_t\Psi
  =\sum_{i=1}^n\widehat\varphi\Big(f_1\circ\gamma_s,\dots,\tfrac{\eta}{\gamma_s'}(f_i\circ\gamma_s)',\dots,f_n\circ\gamma_s\Big)
  =\sum_{i=1}^n\widehat\varphi\big(f_1\circ\gamma_s,\dots,\eta\,(f_i'\circ\gamma_s),\dots,f_n\circ\gamma_s\big).
\end{equation}
Comparing the two above equations gives $ \partial_s\Psi= \frac{\eta}{\gamma_s'}\,\partial_t\Psi$.
As for $\Theta$, the chain rule gives $\partial_s\Theta=(\widehat\varphi(F)'\circ\gamma_s)\,\eta$ and
$\partial_t\Theta=(\widehat\varphi(F)'\circ\gamma_s)\,\gamma_s'$, hence
$\partial_s\Theta=\tfrac{\eta}{\gamma_s'}\partial_t\Theta$ as well. Therefore $D \coloneqq \Psi-\Theta$
satisfies, for $s \in [0,1)$,
\begin{equation}\label{e:transport}
  \partial_sD=\frac{\eta}{\gamma_s'}\,\partial_tD,\qquad D(0)=0.
\end{equation}

We claim that $D\equiv0$. Write $a \coloneqq \eta/\gamma_s'=\sum_{j\ge1}a_j(s)T^j$ (it has $T$-valuation $\ge1$,
as $\eta$ does) and $D=\sum_{d\ge1}D_d(s)T^d$. Comparing the coefficients of $T^d$ in
$\partial_sD=a\,\partial_tD$ gives
\[
  D_d'(s)=a_1(s)\,d\,D_d(s)+\sum_{j=2}^{d}a_j(s)\,(d-j+1)\,D_{d-j+1}(s),
\]
the terms with $j\ge2$ involving only $D_e$ with $e<d$. By induction on $d\ge1$ the lower terms
vanish, so $D_d$ solves the linear homogeneous ODE $D_d'=a_1\,d\,D_d$ with $D_d(0)=0$; hence
$D_d\equiv0$ on $[0,1)$. As $D_d$ is polynomial in $s$, it vanishes at $s=1$ too. Thus $D(1)=0$,
which is~\eqref{e:commute-power}.

Applying~\eqref{e:commute-power} to $F=(T^{x_1},\dots,T^{x_n})$ gives
$c(\lambda x_1,\dots,\lambda x_n)\,T^{\lambda\sum x_i}=c(x_1,\dots,x_n)\,T^{\lambda\sum x_i}$, i.e.
$c(\lambda x)=c(x)$: this is~\eqref{c-homogene}. \vupe{Finally, since $c$ is homogeneous and
polyharmonic, from Theorem~\ref{t:linear}, it is a linear combination of $c_\sigma$'s, concluding the proof of the proposition.}
\end{proof}

\begin{remark}
 In the above proof, note that the function $\eta/\gamma_s'$ to which we applied Lemma~\ref{l:der-formal} is a power series but not a polynomial.
 This explains the need of extending our function $\varphi$ to 
 $\widehat E=\mathbb C[[T]]_0$ to run the argument 
 (in particular, since $\varphi$ does not commute with the multiplication by $\gamma'_s$,
 applying Lemma~\ref{l:der-formal} to $\eta$ does not lead to \eqref{e:dtPsi}).
 \end{remark}

\begin{remark}
\label{r:local-to-global}
Let us give an informal interpretation of Proposition~\ref{p:homog-from-induction} in terms
of invariance through reparametrizations.
That $c$ is polyharmonic (i.e.~Eq.\eqref{recurrence-c}) means that the associated functional $\varphi$
is invariant at first order at $\gamma=\Id$.
The conclusion (homogeneity) and the proof method tells us that this is enough to guarantee the full
invariance by reparemetrization.
This local to global implication is made possible 
by establishing the transport equation~\eqref{e:transport}.
\end{remark}

\subsection{Local-to-global principle in the infinitely differentiable case}

\vf{We now present a similar result
in the case $E=C^\infty(\mathbb R)_0$ under the path-continuity hypothesis of 
Section~\ref{ss:path-cont}.
 Again we shall see}
\vupe{that invariance by any nondecreasing reparametrizations follows
from the invariance by reparametrizations at first order at the identity.
Because of the similarity with Proposition~\ref{p:homog-from-induction} and since we will not re-use this result later, we do not give a full proof but only emphasize 
differences with the framework of formal power series.}
\begin{proposition}[\vupe{Equivalence between local} and global invariance]
\label{p:derivation}
Let $\varphi\colon E^n\to E$ be multilinear and path-continuous.
The following are equivalent.
\begin{enumerate}
  \item For every $F=(f_1,\dots,f_n)\in E^n$ and every $\gamma\in E$ with $\forall t \in \R, \gamma'(t) \geq 0$,
        \[
          \varphi(F\circ\gamma)=\varphi(F)\circ\gamma .
        \]
  \item For every $F=(f_1,\dots,f_n)\in E^n$ and every $\eta\in E$,
        \[
          \eta\,(\varphi(F))'=\sum_{i=1}^n\varphi(f_1,\dots,\eta f_i',\dots,f_n).
        \]
\end{enumerate}
\end{proposition}

\begin{proof}[Sketch of proof]

\emph{$(1)\Rightarrow(2)$.}
First we explain that it is enough to show this implication for $\eta$
such that $\eta'$ is bounded. Indeed, let $t_0 \in \R$ be fixed (we will
assume $t_0>0$, the other case being similar). Let $\overline{\eta} \in E$
be such that $\overline{\eta}_{|[0,t_0]} = \eta_{|[0,t_0]}$ and
$\overline{\eta}'$ is bounded. We have by assumption
\[
          \overline{\eta}\,(\varphi(F))'=\sum_{i=1}^n\varphi(f_1,\dots,\overline{\eta} f_i',\dots,f_n).
\]
By Proposition \ref{p:past},
$\sum_{i=1}^n\varphi(f_1,\dots,\overline{\eta} f_i',\dots,f_n)(t_0)
=\sum_{i=1}^n\varphi(f_1,\dots,\eta f_i',\dots,f_n)(t_0)$.
We deduce
\[
 \overline{\eta}(t_0) \,\varphi(F)'(t_0) = \eta(t_0) \,\varphi(F)'(t_0)
= \sum_{i=1}^n\varphi(f_1,\dots,\eta f_i',\dots,f_n)(t_0),
\]
which is what we wanted to prove.

\vupe{Now, the proof when $\eta$ is bounded is similar to the proof of Lemma~\ref{lem:derivative},
after noticing that if $\gamma_h(t)=t+h\eta(t)$, then $\gamma'_h$ is nonnegative on all $\mathbb R$ for $t$ small enough.}

\medskip
\emph{$(2)\Rightarrow(1)$.}
Fix $F$ and $\gamma$ in $E$ with $\gamma'(t) \geq 0$ for all $t$.
Put $\eta \coloneqq \gamma-\mathrm{id}$ and, for
$s\in[0,1]$, define $\gamma_s=\mathrm{id}+s\eta$.
\vupe{As in the proof of Proposition~\ref{p:homog-from-induction}, the idea is to define
\[
  \Psi(s) \coloneqq \varphi(F\circ\gamma_s),\qquad \Theta(s) \coloneqq \varphi(F)\circ\gamma_s,
\]
and prove that $\Psi(1)=\Theta(1)$. Again, using (2) and a Leibniz rule similar to Lemma~\ref{lem:Leibniz},
one can prove that} $D \coloneqq \Psi-\Theta$ solves the linear transport equation
\[
  \partial_s D=\frac{\eta}{\gamma_s'}\,\partial_t D,\qquad D(0,\cdot)=0.
\]

\vupe{However, since functions in $E$ do not necessarily have power series expansions,
we need to follow a different route from Proposition~\ref{p:homog-from-induction} to
conclude that $D(1,\cdot)=0$.
We}  integrate \vf{our PDE} in the coordinate $u  \coloneqq  \gamma_s(t)$. More precisely, for every
$s \in [0,s_0]$ one has the inequality $\gamma_s'(t) = (1-s) + s\gamma'(t) \ge 1 - s \geq 1 - s_0 > 0$,
so $\gamma_s \colon \mathbb{R} \to \mathbb{R}$ is a strictly increasing bijection;
we may thus set
\[
  \widetilde{D}(s,u) \; \coloneqq \; D\!\left(s,\, \gamma_s^{-1}(u)\right)
\]
and write $t = \gamma_s^{-1}(u)$.
Differentiating the identity $\gamma_s\!\left(\gamma_s^{-1}(u)\right) = u$ with respect
to $s$ gives
\[
 \gamma'_s(t) \cdot \partial_s \!\left(\gamma_s^{-1}(u)\right) + \partial_s \gamma_s(t)=0\,,
\]
so, using $\partial_s \gamma_s(t) = \eta(t)$,
$\partial_s\!\left(\gamma_s^{-1}(u)\right) = -\,\eta(t)/\gamma_s'(t)$.
Hence, by the chain rule and the transport equation,
\[
  \partial_s \widetilde{D}(s,u)
  \;=\; \partial_s D(s,t) - \frac{\eta(t)}{\gamma_s'(t)}\,\partial_t D(s,t)
  \;=\; 0 .
\]
Thus $s \mapsto \widetilde{D}(s,u)$ is constant on $[0,s_0]$ for each $u$. Since
$\gamma_0 = \mathrm{id}$ we have $\widetilde{D}(0,u) = D(0,u) = 0$. Therefore
$\widetilde{D} \equiv 0$, i.e.\ $D \equiv 0$ on $[0,s_0] \times \mathbb{R}$.
In particular $\Psi(s) = \Theta(s)$ for every $s \in [0,s_0]$. As $s_0 \in ]0,1[$ was
arbitrary, this holds for all $s \in [0,1[$, and path-continuity gives
$\Psi(1) = \Theta(1)$, namely $\varphi(F \circ \gamma) = \varphi(F) \circ \gamma$.
%
\end{proof}

%

\end{document}